\pdfoutput=1

\documentclass[12pt,a4paper,twoside]{amsart}

\usepackage{amsbsy}
\usepackage{amssymb}
\usepackage{amscd}
\usepackage{rotating}
\usepackage{amsmath}
\usepackage{bbding}
\usepackage{latexsym}
\usepackage{a4}
\usepackage[all]{xy}
\usepackage[english]{babel} 
\usepackage[latin1]{inputenc} 
\usepackage[T1]{fontenc}

\usepackage{array}
\usepackage{multirow}

\usepackage{listings}
\usepackage{stmaryrd}
\usepackage{cancel}
\usepackage{graphics}
\usepackage{graphicx}
\usepackage{tikz}
\usetikzlibrary{arrows,decorations.pathmorphing,decorations.footprints,backgrounds,positioning,fit,calc,through,shapes}
\usepackage{pgf}
\definecolor{MyLightMagenta}{cmyk}{0.1,0.8,0,0.1}

\usepackage[unicode=false,colorlinks=false]{hyperref}

\newcommand{\secA}[1]{\section{#1}}
\newcommand{\secB}[1]{\subsection{#1}}

\def\N{{\mathbb N}}
\def\Z{{\mathbb Z}}

\def\plus{\texttt{+}}
\def\minus{\texttt{-}}
\def\oo{\texttt{o}}

\def\coloneqq{\mathrel{\mathop:}=}

\newcommand{\lb}{\lbrack}
\newcommand{\rb}{\rbrack}

\def\Aut{\mathrm{Aut}}

\def\id{\mathrm{id\;}}

\def\0{{\bf 0}}
\def\1{{\bf 1}}

\newcommand{\defstil}[1]{\textup{\textbf{\/#1}}}

\newcommand{\names}[1]{#1}

\numberwithin{equation}{section}

\newtheorem{defn}{Definition}[section]

\newtheorem{ex*}{Example}
\newtheorem{lemma}[defn]{Lemma}
\newtheorem{prop}[defn]{Proposition}

\newtheorem{thm}[defn]{Theorem}

\newtheorem{conj}[defn]{Conjecture}

\linethickness{0.4pt}

\makeindex
\def\TubeL{
\begin{tikzpicture}
[scale=0.60,
inner sep=1.5,
vertex/.style={circle,draw=black!50,fill=black!20, very thick},
parta/.style={circle,draw=black!50,fill=black!0, very thick},
partb/.style={circle,draw=black!50,fill=black!100, very thick},
edgevertex/.style={circle,draw=blue!80,fill=blue!40, very thick},
normal/.style=,
post/.style={->, >=stealth'},
pre/.style={<-, >=stealth'}]
\node(a0) at (0.000000,0.000000) {};
\node(a1) at (0.000000,1.000000) {};
\node(a2) at (0.000000,2.000000) {};
\node[vertex](a3) at (2.300000,0.000000) {};
\node[vertex](a4) at (2.300000,1.000000) {};
\node[vertex](a5) at (2.300000,2.000000) {};
\node[vertex](a6) at (4.600000,0.000000) {};
\node[vertex](a7) at (4.600000,1.000000) {};
\node[vertex](a8) at (4.600000,2.000000) {};
\node[vertex](a9) at (6.900000,0.000000) {};
\node[vertex](a10) at (6.900000,1.000000) {};
\node[vertex](a11) at (6.900000,2.000000) {};
\node[vertex](a12) at (9.200000,0.000000) {};
\node[vertex](a13) at (9.200000,1.000000) {};
\node[vertex](a14) at (9.200000,2.000000) {};
\node[vertex](a15) at (11.500000,0.000000) {};
\node[vertex](a16) at (11.500000,1.000000) {};
\node[vertex](a17) at (11.500000,2.000000) {};
\node[vertex](a18) at (13.800000,0.000000) {};
\node[vertex](a19) at (13.800000,1.000000) {};
\node[vertex](a20) at (13.800000,2.000000) {};
\node[vertex](a21) at (16.100000,0.000000) {};
\node[vertex](a22) at (16.100000,1.000000) {};
\node[vertex](a23) at (16.100000,2.000000) {};
\node(a24) at (18.400000,0.000000) {};
\node(a25) at (18.400000,1.000000) {};
\node(a26) at (18.400000,2.000000) {};
\node(a27) at (2.300000,-1.000000) {\tiny -2};
\node(a28) at (4.600000,-1.000000) {\tiny -1};
\node(a29) at (6.900000,-1.000000) {\tiny  0};
\node(a30) at (9.200000,-1.000000) {\tiny  1};
\node(a31) at (11.500000,-1.000000) {\tiny  2};
\node(a32) at (13.800000,-1.000000) {\tiny  3};
\node(a33) at (16.100000,-1.000000) {\tiny  4};
\draw [post,black,thick] (a3) to (a6);
\draw [post,black,thick] (a3) to (a7);
\draw [post,black,thick] (a3) to (a8);
\draw [post,black,thick] (a4) to (a6);
\draw [post,black,thick] (a4) to (a7);
\draw [post,black,thick] (a4) to (a8);
\draw [post,black,thick] (a5) to (a6);
\draw [post,black,thick] (a5) to (a7);
\draw [post,black,thick] (a5) to (a8);
\draw [post,black,thick] (a6) to (a9);
\draw [post,black,thick] (a6) to (a10);
\draw [post,black,thick] (a7) to (a9);
\draw [post,black,thick] (a7) to (a11);
\draw [post,black,thick] (a8) to (a10);
\draw [post,black,thick] (a8) to (a11);
\draw [post,black,thick] (a9) to (a12);
\draw [post,black,thick] (a9) to (a13);
\draw [post,black,thick] (a9) to (a14);
\draw [post,black,thick] (a10) to (a12);
\draw [post,black,thick] (a10) to (a13);
\draw [post,black,thick] (a10) to (a14);
\draw [post,black,thick] (a11) to (a12);
\draw [post,black,thick] (a11) to (a13);
\draw [post,black,thick] (a11) to (a14);
\draw [post,black,thick] (a12) to (a15);
\draw [post,black,thick] (a12) to (a16);
\draw [post,black,thick] (a13) to (a15);
\draw [post,black,thick] (a13) to (a17);
\draw [post,black,thick] (a14) to (a16);
\draw [post,black,thick] (a14) to (a17);
\draw [post,black,thick] (a15) to (a18);
\draw [post,black,thick] (a15) to (a19);
\draw [post,black,thick] (a15) to (a20);
\draw [post,black,thick] (a16) to (a18);
\draw [post,black,thick] (a16) to (a19);
\draw [post,black,thick] (a16) to (a20);
\draw [post,black,thick] (a17) to (a18);
\draw [post,black,thick] (a17) to (a19);
\draw [post,black,thick] (a17) to (a20);
\draw [post,black,thick] (a18) to (a21);
\draw [post,black,thick] (a18) to (a22);
\draw [post,black,thick] (a19) to (a21);
\draw [post,black,thick] (a19) to (a23);
\draw [post,black,thick] (a20) to (a22);
\draw [post,black,thick] (a20) to (a23);
\draw [post,black,dotted] (a0) to (a3);
\draw [post,black,dotted] (a21) to (a24);
\draw [post,black,dotted] (a1) to (a4);
\draw [post,black,dotted] (a22) to (a25);
\draw [post,black,dotted] (a2) to (a5);
\draw [post,black,dotted] (a23) to (a26);
\end{tikzpicture}
}

\def\TubeD{
\begin{tikzpicture}
[scale=0.60,
inner sep=1.5,
vertex/.style={circle,draw=black!50,fill=black!20, very thick},
parta/.style={circle,draw=black!50,fill=black!0, very thick},
partb/.style={circle,draw=black!50,fill=black!100, very thick},
edgevertex/.style={circle,draw=blue!80,fill=blue!40, very thick},
normal/.style=,
post/.style={->, >=stealth'},
pre/.style={<-, >=stealth'}]
\node(a2) at (1.111140,1.662939) {};
\node(a5) at (1.111140,2.662939) {};
\node(a8) at (1.111140,3.662939) {};
\node[vertex](a11) at (3.411140,1.662939) {};
\node[vertex](a14) at (3.411140,2.662939) {};
\node[vertex](a17) at (3.411140,3.662939) {};
\node[vertex](a20) at (5.711140,1.662939) {};
\node[vertex](a23) at (5.711140,2.662939) {};
\node[vertex](a26) at (5.711140,3.662939) {};
\node[vertex](a29) at (8.011140,1.662939) {};
\node[vertex](a32) at (8.011140,2.662939) {};
\node[vertex](a35) at (8.011140,3.662939) {};
\node[vertex](a38) at (10.311140,1.662939) {};
\node[vertex](a41) at (10.311140,2.662939) {};
\node[vertex](a44) at (10.311140,3.662939) {};
\node[vertex](a47) at (12.611140,1.662939) {};
\node[vertex](a50) at (12.611140,2.662939) {};
\node[vertex](a53) at (12.611140,3.662939) {};
\node[vertex](a56) at (14.911140,1.662939) {};
\node[vertex](a59) at (14.911140,2.662939) {};
\node[vertex](a62) at (14.911140,3.662939) {};
\node[vertex](a65) at (17.211140,1.662939) {};
\node[vertex](a68) at (17.211140,2.662939) {};
\node[vertex](a71) at (17.211140,3.662939) {};
\node(a74) at (19.511140,1.662939) {};
\node(a77) at (19.511140,2.662939) {};
\node(a80) at (19.511140,3.662939) {};
\draw [post,black!30,thin] (a11) to (a20);
\draw [post,black!30,thin] (a11) to (a23);
\draw [post,black!30,thin] (a14) to (a20);
\draw [post,black!30,thin] (a14) to (a26);
\draw [post,black!30,thin] (a17) to (a23);
\draw [post,black!30,thin] (a20) to (a29);
\draw [post,black!30,thin] (a20) to (a32);
\draw [post,black!30,thin] (a20) to (a35);
\draw [post,black!30,thin] (a23) to (a29);
\draw [post,black!30,thin] (a23) to (a32);
\draw [post,black!30,thin] (a23) to (a35);
\draw [post,black!30,thin] (a26) to (a29);
\draw [post,black!30,thin] (a26) to (a32);
\draw [post,black!30,thin] (a29) to (a38);
\draw [post,black!30,thin] (a29) to (a41);
\draw [post,black!30,thin] (a32) to (a38);
\draw [post,black!30,thin] (a32) to (a44);
\draw [post,black!30,thin] (a35) to (a41);
\draw [post,black!30,thin] (a38) to (a47);
\draw [post,black!30,thin] (a38) to (a50);
\draw [post,black!30,thin] (a38) to (a53);
\draw [post,black!30,thin] (a41) to (a47);
\draw [post,black!30,thin] (a41) to (a50);
\draw [post,black!30,thin] (a41) to (a53);
\draw [post,black!30,thin] (a44) to (a47);
\draw [post,black!30,thin] (a44) to (a50);
\draw [post,black!30,thin] (a47) to (a56);
\draw [post,black!30,thin] (a47) to (a59);
\draw [post,black!30,thin] (a50) to (a56);
\draw [post,black!30,thin] (a50) to (a62);
\draw [post,black!30,thin] (a53) to (a59);
\draw [post,black!30,thin] (a56) to (a65);
\draw [post,black!30,thin] (a56) to (a68);
\draw [post,black!30,thin] (a56) to (a71);
\draw [post,black!30,thin] (a59) to (a65);
\draw [post,black!30,thin] (a59) to (a68);
\draw [post,black!30,thin] (a59) to (a71);
\draw [post,black!30,thin] (a62) to (a65);
\draw [post,black!30,thin] (a62) to (a68);
\draw [post,black,dotted] (a2) to (a11);
\draw [post,black,dotted] (a65) to (a74);
\draw [post,black,dotted] (a5) to (a14);
\draw [post,black,dotted] (a68) to (a77);
\draw [post,black,dotted] (a8) to (a17);
\draw [post,black,dotted] (a71) to (a80);
\node(a1) at (0.555570,0.831470) {};
\node(a4) at (0.555570,1.831470) {};
\node(a7) at (0.555570,2.831470) {};
\node[vertex](a10) at (2.855570,0.831470) {};
\node[vertex](a13) at (2.855570,1.831470) {};
\node[vertex](a16) at (2.855570,2.831470) {};
\node[vertex](a19) at (5.155570,0.831470) {};
\node[vertex](a22) at (5.155570,1.831470) {};
\node[vertex](a25) at (5.155570,2.831470) {};
\node[vertex](a28) at (7.455570,0.831470) {};
\node[vertex](a31) at (7.455570,1.831470) {};
\node[vertex](a34) at (7.455570,2.831470) {};
\node[vertex](a37) at (9.755570,0.831470) {};
\node[vertex](a40) at (9.755570,1.831470) {};
\node[vertex](a43) at (9.755570,2.831470) {};
\node[vertex](a46) at (12.055570,0.831470) {};
\node[vertex](a49) at (12.055570,1.831470) {};
\node[vertex](a52) at (12.055570,2.831470) {};
\node[vertex](a55) at (14.355570,0.831470) {};
\node[vertex](a58) at (14.355570,1.831470) {};
\node[vertex](a61) at (14.355570,2.831470) {};
\node[vertex](a64) at (16.655570,0.831470) {};
\node[vertex](a67) at (16.655570,1.831470) {};
\node[vertex](a70) at (16.655570,2.831470) {};
\node(a73) at (18.955570,0.831470) {};
\node(a76) at (18.955570,1.831470) {};
\node(a79) at (18.955570,2.831470) {};
\draw [post,black!30,thin] (a10) to (a19);
\draw [post,black!30,thin] (a10) to (a20);
\draw [post,black!30,thin] (a10) to (a22);
\draw [post,black!30,thin] (a10) to (a23);
\draw [post,black!30,thin] (a11) to (a19);
\draw [post,black!30,thin] (a11) to (a22);
\draw [post,black!30,thin] (a13) to (a19);
\draw [post,black!30,thin] (a13) to (a20);
\draw [post,black!30,thin] (a13) to (a25);
\draw [post,black!30,thin] (a13) to (a26);
\draw [post,black!30,thin] (a14) to (a19);
\draw [post,black!30,thin] (a14) to (a25);
\draw [post,black!30,thin] (a16) to (a22);
\draw [post,black!30,thin] (a16) to (a23);
\draw [post,black!30,thin] (a17) to (a22);
\draw [post,black!30,thin] (a19) to (a29);
\draw [post,black!30,thin] (a19) to (a32);
\draw [post,black!30,thin] (a19) to (a35);
\draw [post,black!30,thin] (a20) to (a28);
\draw [post,black!30,thin] (a20) to (a31);
\draw [post,black!30,thin] (a20) to (a34);
\draw [post,black!30,thin] (a22) to (a29);
\draw [post,black!30,thin] (a22) to (a32);
\draw [post,black!30,thin] (a22) to (a35);
\draw [post,black!30,thin] (a23) to (a28);
\draw [post,black!30,thin] (a23) to (a31);
\draw [post,black!30,thin] (a23) to (a34);
\draw [post,black!30,thin] (a25) to (a29);
\draw [post,black!30,thin] (a25) to (a32);
\draw [post,black!30,thin] (a26) to (a28);
\draw [post,black!30,thin] (a26) to (a31);
\draw [post,black!30,thin] (a28) to (a37);
\draw [post,black!30,thin] (a28) to (a38);
\draw [post,black!30,thin] (a28) to (a40);
\draw [post,black!30,thin] (a28) to (a41);
\draw [post,black!30,thin] (a29) to (a37);
\draw [post,black!30,thin] (a29) to (a40);
\draw [post,black!30,thin] (a31) to (a37);
\draw [post,black!30,thin] (a31) to (a38);
\draw [post,black!30,thin] (a31) to (a43);
\draw [post,black!30,thin] (a31) to (a44);
\draw [post,black!30,thin] (a32) to (a37);
\draw [post,black!30,thin] (a32) to (a43);
\draw [post,black!30,thin] (a34) to (a40);
\draw [post,black!30,thin] (a34) to (a41);
\draw [post,black!30,thin] (a35) to (a40);
\draw [post,black!30,thin] (a37) to (a47);
\draw [post,black!30,thin] (a37) to (a50);
\draw [post,black!30,thin] (a37) to (a53);
\draw [post,black!30,thin] (a38) to (a46);
\draw [post,black!30,thin] (a38) to (a49);
\draw [post,black!30,thin] (a38) to (a52);
\draw [post,black!30,thin] (a40) to (a47);
\draw [post,black!30,thin] (a40) to (a50);
\draw [post,black!30,thin] (a40) to (a53);
\draw [post,black!30,thin] (a41) to (a46);
\draw [post,black!30,thin] (a41) to (a49);
\draw [post,black!30,thin] (a41) to (a52);
\draw [post,black!30,thin] (a43) to (a47);
\draw [post,black!30,thin] (a43) to (a50);
\draw [post,black!30,thin] (a44) to (a46);
\draw [post,black!30,thin] (a44) to (a49);
\draw [post,black!30,thin] (a46) to (a55);
\draw [post,black!30,thin] (a46) to (a56);
\draw [post,black!30,thin] (a46) to (a58);
\draw [post,black!30,thin] (a46) to (a59);
\draw [post,black!30,thin] (a47) to (a55);
\draw [post,black!30,thin] (a47) to (a58);
\draw [post,black!30,thin] (a49) to (a55);
\draw [post,black!30,thin] (a49) to (a56);
\draw [post,black!30,thin] (a49) to (a61);
\draw [post,black!30,thin] (a49) to (a62);
\draw [post,black!30,thin] (a50) to (a55);
\draw [post,black!30,thin] (a50) to (a61);
\draw [post,black!30,thin] (a52) to (a58);
\draw [post,black!30,thin] (a52) to (a59);
\draw [post,black!30,thin] (a53) to (a58);
\draw [post,black!30,thin] (a55) to (a65);
\draw [post,black!30,thin] (a55) to (a68);
\draw [post,black!30,thin] (a55) to (a71);
\draw [post,black!30,thin] (a56) to (a64);
\draw [post,black!30,thin] (a56) to (a67);
\draw [post,black!30,thin] (a56) to (a70);
\draw [post,black!30,thin] (a58) to (a65);
\draw [post,black!30,thin] (a58) to (a68);
\draw [post,black!30,thin] (a58) to (a71);
\draw [post,black!30,thin] (a59) to (a64);
\draw [post,black!30,thin] (a59) to (a67);
\draw [post,black!30,thin] (a59) to (a70);
\draw [post,black!30,thin] (a61) to (a65);
\draw [post,black!30,thin] (a61) to (a68);
\draw [post,black!30,thin] (a62) to (a64);
\draw [post,black!30,thin] (a62) to (a67);
\draw [post,black,dotted] (a1) to (a10);
\draw [post,black,dotted] (a64) to (a73);
\draw [post,black,dotted] (a4) to (a13);
\draw [post,black,dotted] (a67) to (a76);
\draw [post,black,dotted] (a7) to (a16);
\draw [post,black,dotted] (a70) to (a79);
\node(a0) at (0.000000,0.000000) {};
\node(a3) at (0.000000,1.000000) {};
\node(a6) at (0.000000,2.000000) {};
\node[vertex](a9) at (2.300000,0.000000) {};
\node[vertex](a12) at (2.300000,1.000000) {};
\node[vertex](a15) at (2.300000,2.000000) {};
\node[vertex](a18) at (4.600000,0.000000) {};
\node[vertex](a21) at (4.600000,1.000000) {};
\node[vertex](a24) at (4.600000,2.000000) {};
\node[vertex](a27) at (6.900000,0.000000) {};
\node[vertex](a30) at (6.900000,1.000000) {};
\node[vertex](a33) at (6.900000,2.000000) {};
\node[vertex](a36) at (9.200000,0.000000) {};
\node[vertex](a39) at (9.200000,1.000000) {};
\node[vertex](a42) at (9.200000,2.000000) {};
\node[vertex](a45) at (11.500000,0.000000) {};
\node[vertex](a48) at (11.500000,1.000000) {};
\node[vertex](a51) at (11.500000,2.000000) {};
\node[vertex](a54) at (13.800000,0.000000) {};
\node[vertex](a57) at (13.800000,1.000000) {};
\node[vertex](a60) at (13.800000,2.000000) {};
\node[vertex](a63) at (16.100000,0.000000) {};
\node[vertex](a66) at (16.100000,1.000000) {};
\node[vertex](a69) at (16.100000,2.000000) {};
\node(a72) at (18.400000,0.000000) {};
\node(a75) at (18.400000,1.000000) {};
\node(a78) at (18.400000,2.000000) {};
\node(a81) at (2.300000,-1.000000) {\tiny -2};
\node(a82) at (4.600000,-1.000000) {\tiny -1};
\node(a83) at (6.900000,-1.000000) {\tiny  0};
\node(a84) at (9.200000,-1.000000) {\tiny  1};
\node(a85) at (11.500000,-1.000000) {\tiny  2};
\node(a86) at (13.800000,-1.000000) {\tiny  3};
\node(a87) at (16.100000,-1.000000) {\tiny  4};
\draw [post,black!30,thin] (a9) to (a19);
\draw [post,black!30,thin] (a9) to (a20);
\draw [post,black!30,thin] (a9) to (a22);
\draw [post,black!30,thin] (a9) to (a23);
\draw [post,black!30,thin] (a10) to (a18);
\draw [post,black!30,thin] (a10) to (a21);
\draw [post,black!30,thin] (a11) to (a18);
\draw [post,black!30,thin] (a11) to (a21);
\draw [post,black!30,thin] (a12) to (a19);
\draw [post,black!30,thin] (a12) to (a20);
\draw [post,black!30,thin] (a12) to (a25);
\draw [post,black!30,thin] (a12) to (a26);
\draw [post,black!30,thin] (a13) to (a18);
\draw [post,black!30,thin] (a13) to (a24);
\draw [post,black!30,thin] (a14) to (a18);
\draw [post,black!30,thin] (a14) to (a24);
\draw [post,black!30,thin] (a15) to (a22);
\draw [post,black!30,thin] (a15) to (a23);
\draw [post,black!30,thin] (a16) to (a21);
\draw [post,black!30,thin] (a17) to (a21);
\draw [post,black!30,thin] (a18) to (a28);
\draw [post,black!30,thin] (a18) to (a31);
\draw [post,black!30,thin] (a18) to (a34);
\draw [post,black!30,thin] (a19) to (a27);
\draw [post,black!30,thin] (a19) to (a30);
\draw [post,black!30,thin] (a19) to (a33);
\draw [post,black!30,thin] (a21) to (a28);
\draw [post,black!30,thin] (a21) to (a31);
\draw [post,black!30,thin] (a21) to (a34);
\draw [post,black!30,thin] (a22) to (a27);
\draw [post,black!30,thin] (a22) to (a30);
\draw [post,black!30,thin] (a22) to (a33);
\draw [post,black!30,thin] (a24) to (a28);
\draw [post,black!30,thin] (a24) to (a31);
\draw [post,black!30,thin] (a25) to (a27);
\draw [post,black!30,thin] (a25) to (a30);
\draw [post,black!30,thin] (a27) to (a37);
\draw [post,black!30,thin] (a27) to (a38);
\draw [post,black!30,thin] (a27) to (a40);
\draw [post,black!30,thin] (a27) to (a41);
\draw [post,black!30,thin] (a28) to (a36);
\draw [post,black!30,thin] (a28) to (a39);
\draw [post,black!30,thin] (a29) to (a36);
\draw [post,black!30,thin] (a29) to (a39);
\draw [post,black!30,thin] (a30) to (a37);
\draw [post,black!30,thin] (a30) to (a38);
\draw [post,black!30,thin] (a30) to (a43);
\draw [post,black!30,thin] (a30) to (a44);
\draw [post,black!30,thin] (a31) to (a36);
\draw [post,black!30,thin] (a31) to (a42);
\draw [post,black!30,thin] (a32) to (a36);
\draw [post,black!30,thin] (a32) to (a42);
\draw [post,black!30,thin] (a33) to (a40);
\draw [post,black!30,thin] (a33) to (a41);
\draw [post,black!30,thin] (a34) to (a39);
\draw [post,black!30,thin] (a35) to (a39);
\draw [post,black!30,thin] (a36) to (a46);
\draw [post,black!30,thin] (a36) to (a49);
\draw [post,black!30,thin] (a36) to (a52);
\draw [post,black!30,thin] (a37) to (a45);
\draw [post,black!30,thin] (a37) to (a48);
\draw [post,black!30,thin] (a37) to (a51);
\draw [post,black!30,thin] (a39) to (a46);
\draw [post,black!30,thin] (a39) to (a49);
\draw [post,black!30,thin] (a39) to (a52);
\draw [post,black!30,thin] (a40) to (a45);
\draw [post,black!30,thin] (a40) to (a48);
\draw [post,black!30,thin] (a40) to (a51);
\draw [post,black!30,thin] (a42) to (a46);
\draw [post,black!30,thin] (a42) to (a49);
\draw [post,black!30,thin] (a43) to (a45);
\draw [post,black!30,thin] (a43) to (a48);
\draw [post,black!30,thin] (a45) to (a55);
\draw [post,black!30,thin] (a45) to (a56);
\draw [post,black!30,thin] (a45) to (a58);
\draw [post,black!30,thin] (a45) to (a59);
\draw [post,black!30,thin] (a46) to (a54);
\draw [post,black!30,thin] (a46) to (a57);
\draw [post,black!30,thin] (a47) to (a54);
\draw [post,black!30,thin] (a47) to (a57);
\draw [post,black!30,thin] (a48) to (a55);
\draw [post,black!30,thin] (a48) to (a56);
\draw [post,black!30,thin] (a48) to (a61);
\draw [post,black!30,thin] (a48) to (a62);
\draw [post,black!30,thin] (a49) to (a54);
\draw [post,black!30,thin] (a49) to (a60);
\draw [post,black!30,thin] (a50) to (a54);
\draw [post,black!30,thin] (a50) to (a60);
\draw [post,black!30,thin] (a51) to (a58);
\draw [post,black!30,thin] (a51) to (a59);
\draw [post,black!30,thin] (a52) to (a57);
\draw [post,black!30,thin] (a53) to (a57);
\draw [post,black!30,thin] (a54) to (a64);
\draw [post,black!30,thin] (a54) to (a67);
\draw [post,black!30,thin] (a54) to (a70);
\draw [post,black!30,thin] (a55) to (a63);
\draw [post,black!30,thin] (a55) to (a66);
\draw [post,black!30,thin] (a55) to (a69);
\draw [post,black!30,thin] (a57) to (a64);
\draw [post,black!30,thin] (a57) to (a67);
\draw [post,black!30,thin] (a57) to (a70);
\draw [post,black!30,thin] (a58) to (a63);
\draw [post,black!30,thin] (a58) to (a66);
\draw [post,black!30,thin] (a58) to (a69);
\draw [post,black!30,thin] (a60) to (a64);
\draw [post,black!30,thin] (a60) to (a67);
\draw [post,black!30,thin] (a61) to (a63);
\draw [post,black!30,thin] (a61) to (a66);
\draw [post,black,dotted] (a0) to (a9);
\draw [post,black,dotted] (a63) to (a72);
\draw [post,black,dotted] (a3) to (a12);
\draw [post,black,dotted] (a66) to (a75);
\draw [post,black,dotted] (a6) to (a15);
\draw [post,black,dotted] (a69) to (a78);
\draw [post,black,very thick] (a48) to (a54);
\draw [post,black,very thick] (a48) to (a60);
\draw [post,black,very thick] (a51) to (a57);
\draw [post,black,very thick] (a51) to (a60);
\draw [post,black,very thick] (a51) to (a61);
\draw [post,black,very thick] (a51) to (a62);
\draw [post,black,very thick] (a52) to (a60);
\draw [post,black,very thick] (a53) to (a60);
\draw [post,black,very thick] (a54) to (a63);
\draw [post,black,very thick] (a54) to (a66);
\draw [post,black,very thick] (a54) to (a69);
\draw [post,black,very thick] (a57) to (a63);
\draw [post,black,very thick] (a57) to (a66);
\draw [post,black,very thick] (a57) to (a69);
\draw [post,black,very thick] (a60) to (a63);
\draw [post,black,very thick] (a60) to (a69);
\draw [post,black,very thick] (a60) to (a70);
\draw [post,black,very thick] (a60) to (a66);
\draw [post,black,very thick] (a61) to (a69);
\draw [post,black,very thick] (a16) to (a24);
\draw [post,black,very thick] (a17) to (a24);
\draw [post,black,very thick] (a18) to (a27);
\draw [post,black,very thick] (a18) to (a30);
\draw [post,black,very thick] (a18) to (a33);
\draw [post,black,very thick] (a21) to (a27);
\draw [post,black,very thick] (a21) to (a30);
\draw [post,black,very thick] (a21) to (a33);
\draw [post,black,very thick] (a24) to (a27);
\draw [post,black,very thick] (a24) to (a30);
\draw [post,black,very thick] (a24) to (a33);
\draw [post,black,very thick] (a24) to (a34);
\draw [post,black,very thick] (a25) to (a33);
\draw [post,black,very thick] (a27) to (a36);
\draw [post,black,very thick] (a27) to (a39);
\draw [post,black,very thick] (a30) to (a36);
\draw [post,black,very thick] (a30) to (a42);
\draw [post,black,very thick] (a33) to (a39);
\draw [post,black,very thick] (a33) to (a42);
\draw [post,black,very thick] (a33) to (a43);
\draw [post,black,very thick] (a33) to (a44);
\draw [post,black,very thick] (a34) to (a42);
\draw [post,black,very thick] (a35) to (a42);
\draw [post,black,very thick] (a36) to (a45);
\draw [post,black,very thick] (a36) to (a48);
\draw [post,black,very thick] (a36) to (a51);
\draw [post,black,very thick] (a39) to (a45);
\draw [post,black,very thick] (a39) to (a48);
\draw [post,black,very thick] (a39) to (a51);
\draw [post,black,very thick] (a42) to (a45);
\draw [post,black,very thick] (a42) to (a48);
\draw [post,black,very thick] (a42) to (a51);
\draw [post,black,very thick] (a42) to (a52);
\draw [post,black,very thick] (a43) to (a51);
\draw [post,black,very thick] (a45) to (a54);
\draw [post,black,very thick] (a45) to (a57);
\draw [post,black,very thick] (a53) to (a62);
\draw [post,black,very thick] (a17) to (a26);
\draw [post,black,very thick] (a26) to (a35);
\draw [post,black,very thick] (a35) to (a44);
\draw [post,black,very thick] (a44) to (a53);
\draw [post,black,very thick] (a62) to (a71);
\draw [post,black,very thick] (a16) to (a25);
\draw [post,black,very thick] (a16) to (a26);
\draw [post,black,very thick] (a17) to (a25);
\draw [post,black,very thick] (a25) to (a35);
\draw [post,black,very thick] (a26) to (a34);
\draw [post,black,very thick] (a34) to (a43);
\draw [post,black,very thick] (a34) to (a44);
\draw [post,black,very thick] (a35) to (a43);
\draw [post,black,very thick] (a43) to (a53);
\draw [post,black,very thick] (a44) to (a52);
\draw [post,black,very thick] (a52) to (a61);
\draw [post,black,very thick] (a52) to (a62);
\draw [post,black,very thick] (a53) to (a61);
\draw [post,black,very thick] (a61) to (a71);
\draw [post,black,very thick] (a62) to (a70);
\draw [post,black,very thick] (a9) to (a18);
\draw [post,black,very thick] (a9) to (a21);
\draw [post,black,very thick] (a12) to (a18);
\draw [post,black,very thick] (a12) to (a24);
\draw [post,black,very thick] (a15) to (a21);
\draw [post,black,very thick] (a15) to (a24);
\draw [post,black,very thick] (a15) to (a25);
\draw [post,black,very thick] (a15) to (a26);
\end{tikzpicture}
}

\def\TubeDeltaAAA{
\begin{tikzpicture}
[scale=0.60,
inner sep=1.5,
vertex/.style={circle,draw=black!50,fill=black!20, very thick},
parta/.style={circle,draw=black!50,fill=black!0, very thick},
partb/.style={circle,draw=black!50,fill=black!100, very thick},
edgevertex/.style={circle,draw=blue!80,fill=blue!40, very thick},
normal/.style=,
post/.style={->, >=stealth'},
pre/.style={<-, >=stealth'}]
\node[vertex](a2) at (1.111140,1.662939) {};
\node[vertex](a5) at (1.111140,2.662939) {};
\node[vertex](a8) at (1.111140,3.662939) {};
\node[vertex](a11) at (3.411140,1.662939) {};
\node[vertex](a14) at (3.411140,2.662939) {};
\node[vertex](a17) at (3.411140,3.662939) {};
\draw [post,black!30,thin] (a2) to (a11);
\draw [post,black!30,thin] (a2) to (a14);
\draw [post,black!30,thin] (a5) to (a11);
\draw [post,black!30,thin] (a5) to (a17);
\draw [post,black!30,thin] (a8) to (a14);
\draw [post,black!30,thin] (a8) to (a17);
\node[vertex](a1) at (0.555570,0.831470) {};
\node[vertex](a4) at (0.555570,1.831470) {};
\node[vertex](a7) at (0.555570,2.831470) {};
\node[vertex](a10) at (2.855570,0.831470) {};
\node[vertex](a13) at (2.855570,1.831470) {};
\node[vertex](a16) at (2.855570,2.831470) {};
\draw [post,black!30,thin] (a1) to (a10);
\draw [post,black!30,thin] (a1) to (a11);
\draw [post,black!30,thin] (a1) to (a13);
\draw [post,black!30,thin] (a1) to (a14);
\draw [post,black!30,thin] (a2) to (a10);
\draw [post,black!30,thin] (a2) to (a13);
\draw [post,black!30,thin] (a4) to (a10);
\draw [post,black!30,thin] (a4) to (a11);
\draw [post,black!30,thin] (a4) to (a16);
\draw [post,black!30,thin] (a4) to (a17);
\draw [post,black!30,thin] (a5) to (a10);
\draw [post,black!30,thin] (a5) to (a16);
\draw [post,black!30,thin] (a7) to (a13);
\draw [post,black!30,thin] (a7) to (a14);
\draw [post,black!30,thin] (a7) to (a16);
\draw [post,black!30,thin] (a7) to (a17);
\draw [post,black!30,thin] (a8) to (a13);
\draw [post,black!30,thin] (a8) to (a16);
\node[vertex](a0) at (0.000000,0.000000) {};
\node[vertex](a3) at (0.000000,1.000000) {};
\node[vertex](a6) at (0.000000,2.000000) {};
\node[vertex](a9) at (2.300000,0.000000) {};
\node[vertex](a12) at (2.300000,1.000000) {};
\node[vertex](a15) at (2.300000,2.000000) {};
\draw [post,black!30,thin] (a0) to (a10);
\draw [post,black!30,thin] (a0) to (a11);
\draw [post,black!30,thin] (a0) to (a13);
\draw [post,black!30,thin] (a0) to (a14);
\draw [post,black!30,thin] (a1) to (a9);
\draw [post,black!30,thin] (a1) to (a12);
\draw [post,black!30,thin] (a2) to (a9);
\draw [post,black!30,thin] (a2) to (a12);
\draw [post,black!30,thin] (a3) to (a10);
\draw [post,black!30,thin] (a3) to (a11);
\draw [post,black!30,thin] (a3) to (a16);
\draw [post,black!30,thin] (a3) to (a17);
\draw [post,black!30,thin] (a4) to (a9);
\draw [post,black!30,thin] (a4) to (a15);
\draw [post,black!30,thin] (a5) to (a9);
\draw [post,black!30,thin] (a5) to (a15);
\draw [post,black!30,thin] (a6) to (a13);
\draw [post,black!30,thin] (a6) to (a14);
\draw [post,black!30,thin] (a6) to (a16);
\draw [post,black!30,thin] (a6) to (a17);
\draw [post,black!30,thin] (a7) to (a12);
\draw [post,black!30,thin] (a7) to (a15);
\draw [post,black!30,thin] (a8) to (a12);
\draw [post,black!30,thin] (a8) to (a15);
\draw [post,black,very thick] (a0) to (a9);
\draw [post,black,very thick] (a0) to (a12);
\draw [post,black,very thick] (a3) to (a9);
\draw [post,black,very thick] (a3) to (a15);
\draw [post,black,very thick] (a6) to (a12);
\draw [post,black,very thick] (a6) to (a15);
\end{tikzpicture}
}

\def\TubeDeltaABA{
\begin{tikzpicture}
[scale=0.60,
inner sep=1.5,
vertex/.style={circle,draw=black!50,fill=black!20, very thick},
parta/.style={circle,draw=black!50,fill=black!0, very thick},
partb/.style={circle,draw=black!50,fill=black!100, very thick},
edgevertex/.style={circle,draw=blue!80,fill=blue!40, very thick},
normal/.style=,
post/.style={->, >=stealth'},
pre/.style={<-, >=stealth'}]
\node[vertex](a2) at (1.111140,1.662939) {};
\node[vertex](a5) at (1.111140,2.662939) {};
\node[vertex](a8) at (1.111140,3.662939) {};
\node[vertex](a11) at (3.411140,1.662939) {};
\node[vertex](a14) at (3.411140,2.662939) {};
\node[vertex](a17) at (3.411140,3.662939) {};
\draw [post,black!30,thin] (a2) to (a11);
\draw [post,black!30,thin] (a2) to (a14);
\draw [post,black!30,thin] (a5) to (a11);
\draw [post,black!30,thin] (a5) to (a17);
\draw [post,black!30,thin] (a8) to (a14);
\draw [post,black!30,thin] (a8) to (a17);
\node[vertex](a1) at (0.555570,0.831470) {};
\node[vertex](a4) at (0.555570,1.831470) {};
\node[vertex](a7) at (0.555570,2.831470) {};
\node[vertex](a10) at (2.855570,0.831470) {};
\node[vertex](a13) at (2.855570,1.831470) {};
\node[vertex](a16) at (2.855570,2.831470) {};
\draw [post,black!30,thin] (a1) to (a10);
\draw [post,black!30,thin] (a1) to (a11);
\draw [post,black!30,thin] (a1) to (a13);
\draw [post,black!30,thin] (a1) to (a14);
\draw [post,black!30,thin] (a2) to (a10);
\draw [post,black!30,thin] (a2) to (a13);
\draw [post,black!30,thin] (a4) to (a10);
\draw [post,black!30,thin] (a4) to (a11);
\draw [post,black!30,thin] (a4) to (a16);
\draw [post,black!30,thin] (a4) to (a17);
\draw [post,black!30,thin] (a5) to (a10);
\draw [post,black!30,thin] (a5) to (a16);
\draw [post,black!30,thin] (a7) to (a13);
\draw [post,black!30,thin] (a7) to (a14);
\draw [post,black!30,thin] (a7) to (a16);
\draw [post,black!30,thin] (a7) to (a17);
\draw [post,black!30,thin] (a8) to (a13);
\draw [post,black!30,thin] (a8) to (a16);
\node[vertex](a0) at (0.000000,0.000000) {};
\node[vertex](a3) at (0.000000,1.000000) {};
\node[vertex](a6) at (0.000000,2.000000) {};
\node[vertex](a9) at (2.300000,0.000000) {};
\node[vertex](a12) at (2.300000,1.000000) {};
\node[vertex](a15) at (2.300000,2.000000) {};
\draw [post,black!30,thin] (a0) to (a9);
\draw [post,black!30,thin] (a0) to (a11);
\draw [post,black!30,thin] (a0) to (a12);
\draw [post,black!30,thin] (a0) to (a14);
\draw [post,black!30,thin] (a1) to (a9);
\draw [post,black!30,thin] (a1) to (a12);
\draw [post,black!30,thin] (a2) to (a9);
\draw [post,black!30,thin] (a2) to (a12);
\draw [post,black!30,thin] (a3) to (a9);
\draw [post,black!30,thin] (a3) to (a11);
\draw [post,black!30,thin] (a3) to (a15);
\draw [post,black!30,thin] (a3) to (a17);
\draw [post,black!30,thin] (a4) to (a9);
\draw [post,black!30,thin] (a4) to (a15);
\draw [post,black!30,thin] (a5) to (a9);
\draw [post,black!30,thin] (a5) to (a15);
\draw [post,black!30,thin] (a6) to (a12);
\draw [post,black!30,thin] (a6) to (a14);
\draw [post,black!30,thin] (a6) to (a15);
\draw [post,black!30,thin] (a6) to (a17);
\draw [post,black!30,thin] (a7) to (a12);
\draw [post,black!30,thin] (a7) to (a15);
\draw [post,black!30,thin] (a8) to (a12);
\draw [post,black!30,thin] (a8) to (a15);
\draw [post,black,very thick] (a6) to (a13);
\draw [post,black,very thick] (a0) to (a10);
\draw [post,black,very thick] (a0) to (a13);
\draw [post,black,very thick] (a3) to (a10);
\draw [post,black,very thick] (a3) to (a16);
\draw [post,black,very thick] (a6) to (a16);
\end{tikzpicture}
}

\def\TubeDeltaACA{
\begin{tikzpicture}
[scale=0.60,
inner sep=1.5,
vertex/.style={circle,draw=black!50,fill=black!20, very thick},
parta/.style={circle,draw=black!50,fill=black!0, very thick},
partb/.style={circle,draw=black!50,fill=black!100, very thick},
edgevertex/.style={circle,draw=blue!80,fill=blue!40, very thick},
normal/.style=,
post/.style={->, >=stealth'},
pre/.style={<-, >=stealth'}]
\node[vertex](a2) at (1.111140,1.662939) {};
\node[vertex](a5) at (1.111140,2.662939) {};
\node[vertex](a8) at (1.111140,3.662939) {};
\node[vertex](a11) at (3.411140,1.662939) {};
\node[vertex](a14) at (3.411140,2.662939) {};
\node[vertex](a17) at (3.411140,3.662939) {};
\draw [post,black!30,thin] (a2) to (a11);
\draw [post,black!30,thin] (a2) to (a14);
\draw [post,black!30,thin] (a5) to (a11);
\draw [post,black!30,thin] (a5) to (a17);
\draw [post,black!30,thin] (a8) to (a14);
\draw [post,black!30,thin] (a8) to (a17);
\node[vertex](a1) at (0.555570,0.831470) {};
\node[vertex](a4) at (0.555570,1.831470) {};
\node[vertex](a7) at (0.555570,2.831470) {};
\node[vertex](a10) at (2.855570,0.831470) {};
\node[vertex](a13) at (2.855570,1.831470) {};
\node[vertex](a16) at (2.855570,2.831470) {};
\draw [post,black!30,thin] (a1) to (a10);
\draw [post,black!30,thin] (a1) to (a11);
\draw [post,black!30,thin] (a1) to (a13);
\draw [post,black!30,thin] (a1) to (a14);
\draw [post,black!30,thin] (a2) to (a10);
\draw [post,black!30,thin] (a2) to (a13);
\draw [post,black!30,thin] (a4) to (a10);
\draw [post,black!30,thin] (a4) to (a11);
\draw [post,black!30,thin] (a4) to (a16);
\draw [post,black!30,thin] (a4) to (a17);
\draw [post,black!30,thin] (a5) to (a10);
\draw [post,black!30,thin] (a5) to (a16);
\draw [post,black!30,thin] (a7) to (a13);
\draw [post,black!30,thin] (a7) to (a14);
\draw [post,black!30,thin] (a7) to (a16);
\draw [post,black!30,thin] (a7) to (a17);
\draw [post,black!30,thin] (a8) to (a13);
\draw [post,black!30,thin] (a8) to (a16);
\node[vertex](a0) at (0.000000,0.000000) {};
\node[vertex](a3) at (0.000000,1.000000) {};
\node[vertex](a6) at (0.000000,2.000000) {};
\node[vertex](a9) at (2.300000,0.000000) {};
\node[vertex](a12) at (2.300000,1.000000) {};
\node[vertex](a15) at (2.300000,2.000000) {};
\draw [post,black!30,thin] (a0) to (a9);
\draw [post,black!30,thin] (a0) to (a10);
\draw [post,black!30,thin] (a0) to (a12);
\draw [post,black!30,thin] (a0) to (a13);
\draw [post,black!30,thin] (a1) to (a9);
\draw [post,black!30,thin] (a1) to (a12);
\draw [post,black!30,thin] (a2) to (a9);
\draw [post,black!30,thin] (a2) to (a12);
\draw [post,black!30,thin] (a3) to (a9);
\draw [post,black!30,thin] (a3) to (a10);
\draw [post,black!30,thin] (a3) to (a15);
\draw [post,black!30,thin] (a3) to (a16);
\draw [post,black!30,thin] (a4) to (a9);
\draw [post,black!30,thin] (a4) to (a15);
\draw [post,black!30,thin] (a5) to (a9);
\draw [post,black!30,thin] (a5) to (a15);
\draw [post,black!30,thin] (a6) to (a12);
\draw [post,black!30,thin] (a6) to (a13);
\draw [post,black!30,thin] (a6) to (a15);
\draw [post,black!30,thin] (a6) to (a16);
\draw [post,black!30,thin] (a7) to (a12);
\draw [post,black!30,thin] (a7) to (a15);
\draw [post,black!30,thin] (a8) to (a12);
\draw [post,black!30,thin] (a8) to (a15);
\draw [post,black,very thick] (a6) to (a17);
\draw [post,black,very thick] (a6) to (a14);
\draw [post,black,very thick] (a0) to (a11);
\draw [post,black,very thick] (a0) to (a14);
\draw [post,black,very thick] (a3) to (a11);
\draw [post,black,very thick] (a3) to (a17);
\end{tikzpicture}
}

\def\TubeDeltaBAA{
\begin{tikzpicture}
[scale=0.60,
inner sep=1.5,
vertex/.style={circle,draw=black!50,fill=black!20, very thick},
parta/.style={circle,draw=black!50,fill=black!0, very thick},
partb/.style={circle,draw=black!50,fill=black!100, very thick},
edgevertex/.style={circle,draw=blue!80,fill=blue!40, very thick},
normal/.style=,
post/.style={->, >=stealth'},
pre/.style={<-, >=stealth'}]
\node[vertex](a2) at (1.111140,1.662939) {};
\node[vertex](a5) at (1.111140,2.662939) {};
\node[vertex](a8) at (1.111140,3.662939) {};
\node[vertex](a11) at (3.411140,1.662939) {};
\node[vertex](a14) at (3.411140,2.662939) {};
\node[vertex](a17) at (3.411140,3.662939) {};
\draw [post,black!30,thin] (a2) to (a11);
\draw [post,black!30,thin] (a2) to (a14);
\draw [post,black!30,thin] (a5) to (a11);
\draw [post,black!30,thin] (a5) to (a17);
\draw [post,black!30,thin] (a8) to (a14);
\draw [post,black!30,thin] (a8) to (a17);
\node[vertex](a1) at (0.555570,0.831470) {};
\node[vertex](a4) at (0.555570,1.831470) {};
\node[vertex](a7) at (0.555570,2.831470) {};
\node[vertex](a10) at (2.855570,0.831470) {};
\node[vertex](a13) at (2.855570,1.831470) {};
\node[vertex](a16) at (2.855570,2.831470) {};
\draw [post,black!30,thin] (a1) to (a10);
\draw [post,black!30,thin] (a1) to (a11);
\draw [post,black!30,thin] (a1) to (a13);
\draw [post,black!30,thin] (a1) to (a14);
\draw [post,black!30,thin] (a2) to (a10);
\draw [post,black!30,thin] (a2) to (a13);
\draw [post,black!30,thin] (a4) to (a10);
\draw [post,black!30,thin] (a4) to (a11);
\draw [post,black!30,thin] (a4) to (a16);
\draw [post,black!30,thin] (a4) to (a17);
\draw [post,black!30,thin] (a5) to (a10);
\draw [post,black!30,thin] (a5) to (a16);
\draw [post,black!30,thin] (a7) to (a13);
\draw [post,black!30,thin] (a7) to (a14);
\draw [post,black!30,thin] (a7) to (a16);
\draw [post,black!30,thin] (a7) to (a17);
\draw [post,black!30,thin] (a8) to (a13);
\draw [post,black!30,thin] (a8) to (a16);
\node[vertex](a0) at (0.000000,0.000000) {};
\node[vertex](a3) at (0.000000,1.000000) {};
\node[vertex](a6) at (0.000000,2.000000) {};
\node[vertex](a9) at (2.300000,0.000000) {};
\node[vertex](a12) at (2.300000,1.000000) {};
\node[vertex](a15) at (2.300000,2.000000) {};
\draw [post,black!30,thin] (a0) to (a9);
\draw [post,black!30,thin] (a0) to (a10);
\draw [post,black!30,thin] (a0) to (a11);
\draw [post,black!30,thin] (a0) to (a12);
\draw [post,black!30,thin] (a0) to (a13);
\draw [post,black!30,thin] (a0) to (a14);
\draw [post,black!30,thin] (a2) to (a9);
\draw [post,black!30,thin] (a2) to (a12);
\draw [post,black!30,thin] (a3) to (a9);
\draw [post,black!30,thin] (a3) to (a10);
\draw [post,black!30,thin] (a3) to (a11);
\draw [post,black!30,thin] (a3) to (a15);
\draw [post,black!30,thin] (a3) to (a16);
\draw [post,black!30,thin] (a3) to (a17);
\draw [post,black!30,thin] (a5) to (a9);
\draw [post,black!30,thin] (a5) to (a15);
\draw [post,black!30,thin] (a6) to (a12);
\draw [post,black!30,thin] (a6) to (a13);
\draw [post,black!30,thin] (a6) to (a14);
\draw [post,black!30,thin] (a6) to (a15);
\draw [post,black!30,thin] (a6) to (a16);
\draw [post,black!30,thin] (a6) to (a17);
\draw [post,black!30,thin] (a8) to (a12);
\draw [post,black!30,thin] (a8) to (a15);
\draw [post,black,very thick] (a7) to (a12);
\draw [post,black,very thick] (a7) to (a15);
\draw [post,black,very thick] (a1) to (a9);
\draw [post,black,very thick] (a1) to (a12);
\draw [post,black,very thick] (a4) to (a9);
\draw [post,black,very thick] (a4) to (a15);
\end{tikzpicture}
}

\def\TubeDeltaBBA{
\begin{tikzpicture}
[scale=0.60,
inner sep=1.5,
vertex/.style={circle,draw=black!50,fill=black!20, very thick},
parta/.style={circle,draw=black!50,fill=black!0, very thick},
partb/.style={circle,draw=black!50,fill=black!100, very thick},
edgevertex/.style={circle,draw=blue!80,fill=blue!40, very thick},
normal/.style=,
post/.style={->, >=stealth'},
pre/.style={<-, >=stealth'}]
\node[vertex](a2) at (1.111140,1.662939) {};
\node[vertex](a5) at (1.111140,2.662939) {};
\node[vertex](a8) at (1.111140,3.662939) {};
\node[vertex](a11) at (3.411140,1.662939) {};
\node[vertex](a14) at (3.411140,2.662939) {};
\node[vertex](a17) at (3.411140,3.662939) {};
\draw [post,black!30,thin] (a2) to (a11);
\draw [post,black!30,thin] (a2) to (a14);
\draw [post,black!30,thin] (a5) to (a11);
\draw [post,black!30,thin] (a5) to (a17);
\draw [post,black!30,thin] (a8) to (a14);
\draw [post,black!30,thin] (a8) to (a17);
\node[vertex](a1) at (0.555570,0.831470) {};
\node[vertex](a4) at (0.555570,1.831470) {};
\node[vertex](a7) at (0.555570,2.831470) {};
\node[vertex](a10) at (2.855570,0.831470) {};
\node[vertex](a13) at (2.855570,1.831470) {};
\node[vertex](a16) at (2.855570,2.831470) {};
\draw [post,black!30,thin] (a1) to (a11);
\draw [post,black!30,thin] (a1) to (a14);
\draw [post,black!30,thin] (a2) to (a10);
\draw [post,black!30,thin] (a2) to (a13);
\draw [post,black!30,thin] (a4) to (a11);
\draw [post,black!30,thin] (a4) to (a17);
\draw [post,black!30,thin] (a5) to (a10);
\draw [post,black!30,thin] (a5) to (a16);
\draw [post,black!30,thin] (a7) to (a14);
\draw [post,black!30,thin] (a7) to (a17);
\draw [post,black!30,thin] (a8) to (a13);
\draw [post,black!30,thin] (a8) to (a16);
\node[vertex](a0) at (0.000000,0.000000) {};
\node[vertex](a3) at (0.000000,1.000000) {};
\node[vertex](a6) at (0.000000,2.000000) {};
\node[vertex](a9) at (2.300000,0.000000) {};
\node[vertex](a12) at (2.300000,1.000000) {};
\node[vertex](a15) at (2.300000,2.000000) {};
\draw [post,black!30,thin] (a0) to (a9);
\draw [post,black!30,thin] (a0) to (a10);
\draw [post,black!30,thin] (a0) to (a11);
\draw [post,black!30,thin] (a0) to (a12);
\draw [post,black!30,thin] (a0) to (a13);
\draw [post,black!30,thin] (a0) to (a14);
\draw [post,black!30,thin] (a1) to (a9);
\draw [post,black!30,thin] (a1) to (a12);
\draw [post,black!30,thin] (a2) to (a9);
\draw [post,black!30,thin] (a2) to (a12);
\draw [post,black!30,thin] (a3) to (a9);
\draw [post,black!30,thin] (a3) to (a10);
\draw [post,black!30,thin] (a3) to (a11);
\draw [post,black!30,thin] (a3) to (a15);
\draw [post,black!30,thin] (a3) to (a16);
\draw [post,black!30,thin] (a3) to (a17);
\draw [post,black!30,thin] (a4) to (a9);
\draw [post,black!30,thin] (a4) to (a15);
\draw [post,black!30,thin] (a5) to (a9);
\draw [post,black!30,thin] (a5) to (a15);
\draw [post,black!30,thin] (a6) to (a12);
\draw [post,black!30,thin] (a6) to (a13);
\draw [post,black!30,thin] (a6) to (a14);
\draw [post,black!30,thin] (a6) to (a15);
\draw [post,black!30,thin] (a6) to (a16);
\draw [post,black!30,thin] (a6) to (a17);
\draw [post,black!30,thin] (a7) to (a12);
\draw [post,black!30,thin] (a7) to (a15);
\draw [post,black!30,thin] (a8) to (a12);
\draw [post,black!30,thin] (a8) to (a15);
\draw [post,black,very thick] (a1) to (a10);
\draw [post,black,very thick] (a1) to (a13);
\draw [post,black,very thick] (a4) to (a10);
\draw [post,black,very thick] (a4) to (a16);
\draw [post,black,very thick] (a7) to (a13);
\draw [post,black,very thick] (a7) to (a16);
\end{tikzpicture}
}

\def\TubeDeltaBCA{
\begin{tikzpicture}
[scale=0.60,
inner sep=1.5,
vertex/.style={circle,draw=black!50,fill=black!20, very thick},
parta/.style={circle,draw=black!50,fill=black!0, very thick},
partb/.style={circle,draw=black!50,fill=black!100, very thick},
edgevertex/.style={circle,draw=blue!80,fill=blue!40, very thick},
normal/.style=,
post/.style={->, >=stealth'},
pre/.style={<-, >=stealth'}]
\node[vertex](a2) at (1.111140,1.662939) {};
\node[vertex](a5) at (1.111140,2.662939) {};
\node[vertex](a8) at (1.111140,3.662939) {};
\node[vertex](a11) at (3.411140,1.662939) {};
\node[vertex](a14) at (3.411140,2.662939) {};
\node[vertex](a17) at (3.411140,3.662939) {};
\draw [post,black!30,thin] (a2) to (a11);
\draw [post,black!30,thin] (a2) to (a14);
\draw [post,black!30,thin] (a5) to (a11);
\draw [post,black!30,thin] (a5) to (a17);
\draw [post,black!30,thin] (a8) to (a14);
\draw [post,black!30,thin] (a8) to (a17);
\node[vertex](a1) at (0.555570,0.831470) {};
\node[vertex](a4) at (0.555570,1.831470) {};
\node[vertex](a7) at (0.555570,2.831470) {};
\node[vertex](a10) at (2.855570,0.831470) {};
\node[vertex](a13) at (2.855570,1.831470) {};
\node[vertex](a16) at (2.855570,2.831470) {};
\draw [post,black!30,thin] (a1) to (a10);
\draw [post,black!30,thin] (a1) to (a13);
\draw [post,black!30,thin] (a2) to (a10);
\draw [post,black!30,thin] (a2) to (a13);
\draw [post,black!30,thin] (a4) to (a10);
\draw [post,black!30,thin] (a4) to (a16);
\draw [post,black!30,thin] (a5) to (a10);
\draw [post,black!30,thin] (a5) to (a16);
\draw [post,black!30,thin] (a7) to (a13);
\draw [post,black!30,thin] (a7) to (a16);
\draw [post,black!30,thin] (a8) to (a13);
\draw [post,black!30,thin] (a8) to (a16);
\node[vertex](a0) at (0.000000,0.000000) {};
\node[vertex](a3) at (0.000000,1.000000) {};
\node[vertex](a6) at (0.000000,2.000000) {};
\node[vertex](a9) at (2.300000,0.000000) {};
\node[vertex](a12) at (2.300000,1.000000) {};
\node[vertex](a15) at (2.300000,2.000000) {};
\draw [post,black!30,thin] (a0) to (a9);
\draw [post,black!30,thin] (a0) to (a10);
\draw [post,black!30,thin] (a0) to (a11);
\draw [post,black!30,thin] (a0) to (a12);
\draw [post,black!30,thin] (a0) to (a13);
\draw [post,black!30,thin] (a0) to (a14);
\draw [post,black!30,thin] (a1) to (a9);
\draw [post,black!30,thin] (a1) to (a12);
\draw [post,black!30,thin] (a2) to (a9);
\draw [post,black!30,thin] (a2) to (a12);
\draw [post,black!30,thin] (a3) to (a9);
\draw [post,black!30,thin] (a3) to (a10);
\draw [post,black!30,thin] (a3) to (a11);
\draw [post,black!30,thin] (a3) to (a15);
\draw [post,black!30,thin] (a3) to (a16);
\draw [post,black!30,thin] (a3) to (a17);
\draw [post,black!30,thin] (a4) to (a9);
\draw [post,black!30,thin] (a4) to (a15);
\draw [post,black!30,thin] (a5) to (a9);
\draw [post,black!30,thin] (a5) to (a15);
\draw [post,black!30,thin] (a6) to (a12);
\draw [post,black!30,thin] (a6) to (a13);
\draw [post,black!30,thin] (a6) to (a14);
\draw [post,black!30,thin] (a6) to (a15);
\draw [post,black!30,thin] (a6) to (a16);
\draw [post,black!30,thin] (a6) to (a17);
\draw [post,black!30,thin] (a7) to (a12);
\draw [post,black!30,thin] (a7) to (a15);
\draw [post,black!30,thin] (a8) to (a12);
\draw [post,black!30,thin] (a8) to (a15);
\draw [post,black,very thick] (a1) to (a11);
\draw [post,black,very thick] (a1) to (a14);
\draw [post,black,very thick] (a4) to (a11);
\draw [post,black,very thick] (a4) to (a17);
\draw [post,black,very thick] (a7) to (a14);
\draw [post,black,very thick] (a7) to (a17);
\end{tikzpicture}
}

\def\TubeDeltaCAA{
\begin{tikzpicture}
[scale=0.60,
inner sep=1.5,
vertex/.style={circle,draw=black!50,fill=black!20, very thick},
parta/.style={circle,draw=black!50,fill=black!0, very thick},
partb/.style={circle,draw=black!50,fill=black!100, very thick},
edgevertex/.style={circle,draw=blue!80,fill=blue!40, very thick},
normal/.style=,
post/.style={->, >=stealth'},
pre/.style={<-, >=stealth'}]
\node[vertex](a2) at (1.111140,1.662939) {};
\node[vertex](a5) at (1.111140,2.662939) {};
\node[vertex](a8) at (1.111140,3.662939) {};
\node[vertex](a11) at (3.411140,1.662939) {};
\node[vertex](a14) at (3.411140,2.662939) {};
\node[vertex](a17) at (3.411140,3.662939) {};
\draw [post,black!30,thin] (a2) to (a11);
\draw [post,black!30,thin] (a2) to (a14);
\draw [post,black!30,thin] (a5) to (a11);
\draw [post,black!30,thin] (a5) to (a17);
\draw [post,black!30,thin] (a8) to (a14);
\draw [post,black!30,thin] (a8) to (a17);
\node[vertex](a1) at (0.555570,0.831470) {};
\node[vertex](a4) at (0.555570,1.831470) {};
\node[vertex](a7) at (0.555570,2.831470) {};
\node[vertex](a10) at (2.855570,0.831470) {};
\node[vertex](a13) at (2.855570,1.831470) {};
\node[vertex](a16) at (2.855570,2.831470) {};
\draw [post,black!30,thin] (a1) to (a10);
\draw [post,black!30,thin] (a1) to (a11);
\draw [post,black!30,thin] (a1) to (a13);
\draw [post,black!30,thin] (a1) to (a14);
\draw [post,black!30,thin] (a2) to (a10);
\draw [post,black!30,thin] (a2) to (a13);
\draw [post,black!30,thin] (a4) to (a10);
\draw [post,black!30,thin] (a4) to (a11);
\draw [post,black!30,thin] (a4) to (a16);
\draw [post,black!30,thin] (a4) to (a17);
\draw [post,black!30,thin] (a5) to (a10);
\draw [post,black!30,thin] (a5) to (a16);
\draw [post,black!30,thin] (a7) to (a13);
\draw [post,black!30,thin] (a7) to (a14);
\draw [post,black!30,thin] (a7) to (a16);
\draw [post,black!30,thin] (a7) to (a17);
\draw [post,black!30,thin] (a8) to (a13);
\draw [post,black!30,thin] (a8) to (a16);
\node[vertex](a0) at (0.000000,0.000000) {};
\node[vertex](a3) at (0.000000,1.000000) {};
\node[vertex](a6) at (0.000000,2.000000) {};
\node[vertex](a9) at (2.300000,0.000000) {};
\node[vertex](a12) at (2.300000,1.000000) {};
\node[vertex](a15) at (2.300000,2.000000) {};
\draw [post,black!30,thin] (a0) to (a9);
\draw [post,black!30,thin] (a0) to (a10);
\draw [post,black!30,thin] (a0) to (a11);
\draw [post,black!30,thin] (a0) to (a12);
\draw [post,black!30,thin] (a0) to (a13);
\draw [post,black!30,thin] (a0) to (a14);
\draw [post,black!30,thin] (a1) to (a9);
\draw [post,black!30,thin] (a1) to (a12);
\draw [post,black!30,thin] (a3) to (a9);
\draw [post,black!30,thin] (a3) to (a10);
\draw [post,black!30,thin] (a3) to (a11);
\draw [post,black!30,thin] (a3) to (a15);
\draw [post,black!30,thin] (a3) to (a16);
\draw [post,black!30,thin] (a3) to (a17);
\draw [post,black!30,thin] (a4) to (a9);
\draw [post,black!30,thin] (a4) to (a15);
\draw [post,black!30,thin] (a6) to (a12);
\draw [post,black!30,thin] (a6) to (a13);
\draw [post,black!30,thin] (a6) to (a14);
\draw [post,black!30,thin] (a6) to (a15);
\draw [post,black!30,thin] (a6) to (a16);
\draw [post,black!30,thin] (a6) to (a17);
\draw [post,black!30,thin] (a7) to (a12);
\draw [post,black!30,thin] (a7) to (a15);
\draw [post,black,very thick] (a2) to (a9);
\draw [post,black,very thick] (a2) to (a12);
\draw [post,black,very thick] (a8) to (a12);
\draw [post,black,very thick] (a8) to (a15);
\draw [post,black,very thick] (a5) to (a9);
\draw [post,black,very thick] (a5) to (a15);
\end{tikzpicture}
}

\def\TubeDeltaCBA{
\begin{tikzpicture}
[scale=0.60,
inner sep=1.5,
vertex/.style={circle,draw=black!50,fill=black!20, very thick},
parta/.style={circle,draw=black!50,fill=black!0, very thick},
partb/.style={circle,draw=black!50,fill=black!100, very thick},
edgevertex/.style={circle,draw=blue!80,fill=blue!40, very thick},
normal/.style=,
post/.style={->, >=stealth'},
pre/.style={<-, >=stealth'}]
\node[vertex](a2) at (1.111140,1.662939) {};
\node[vertex](a5) at (1.111140,2.662939) {};
\node[vertex](a8) at (1.111140,3.662939) {};
\node[vertex](a11) at (3.411140,1.662939) {};
\node[vertex](a14) at (3.411140,2.662939) {};
\node[vertex](a17) at (3.411140,3.662939) {};
\draw [post,black!30,thin] (a2) to (a11);
\draw [post,black!30,thin] (a2) to (a14);
\draw [post,black!30,thin] (a5) to (a11);
\draw [post,black!30,thin] (a5) to (a17);
\draw [post,black!30,thin] (a8) to (a14);
\draw [post,black!30,thin] (a8) to (a17);
\node[vertex](a1) at (0.555570,0.831470) {};
\node[vertex](a4) at (0.555570,1.831470) {};
\node[vertex](a7) at (0.555570,2.831470) {};
\node[vertex](a10) at (2.855570,0.831470) {};
\node[vertex](a13) at (2.855570,1.831470) {};
\node[vertex](a16) at (2.855570,2.831470) {};
\draw [post,black!30,thin] (a1) to (a10);
\draw [post,black!30,thin] (a1) to (a11);
\draw [post,black!30,thin] (a1) to (a13);
\draw [post,black!30,thin] (a1) to (a14);
\draw [post,black!30,thin] (a4) to (a10);
\draw [post,black!30,thin] (a4) to (a11);
\draw [post,black!30,thin] (a4) to (a16);
\draw [post,black!30,thin] (a4) to (a17);
\draw [post,black!30,thin] (a7) to (a13);
\draw [post,black!30,thin] (a7) to (a14);
\draw [post,black!30,thin] (a7) to (a16);
\draw [post,black!30,thin] (a7) to (a17);
\node[vertex](a0) at (0.000000,0.000000) {};
\node[vertex](a3) at (0.000000,1.000000) {};
\node[vertex](a6) at (0.000000,2.000000) {};
\node[vertex](a9) at (2.300000,0.000000) {};
\node[vertex](a12) at (2.300000,1.000000) {};
\node[vertex](a15) at (2.300000,2.000000) {};
\draw [post,black!30,thin] (a0) to (a9);
\draw [post,black!30,thin] (a0) to (a10);
\draw [post,black!30,thin] (a0) to (a11);
\draw [post,black!30,thin] (a0) to (a12);
\draw [post,black!30,thin] (a0) to (a13);
\draw [post,black!30,thin] (a0) to (a14);
\draw [post,black!30,thin] (a1) to (a9);
\draw [post,black!30,thin] (a1) to (a12);
\draw [post,black!30,thin] (a2) to (a9);
\draw [post,black!30,thin] (a2) to (a12);
\draw [post,black!30,thin] (a3) to (a9);
\draw [post,black!30,thin] (a3) to (a10);
\draw [post,black!30,thin] (a3) to (a11);
\draw [post,black!30,thin] (a3) to (a15);
\draw [post,black!30,thin] (a3) to (a16);
\draw [post,black!30,thin] (a3) to (a17);
\draw [post,black!30,thin] (a4) to (a9);
\draw [post,black!30,thin] (a4) to (a15);
\draw [post,black!30,thin] (a5) to (a9);
\draw [post,black!30,thin] (a5) to (a15);
\draw [post,black!30,thin] (a6) to (a12);
\draw [post,black!30,thin] (a6) to (a13);
\draw [post,black!30,thin] (a6) to (a14);
\draw [post,black!30,thin] (a6) to (a15);
\draw [post,black!30,thin] (a6) to (a16);
\draw [post,black!30,thin] (a6) to (a17);
\draw [post,black!30,thin] (a7) to (a12);
\draw [post,black!30,thin] (a7) to (a15);
\draw [post,black!30,thin] (a8) to (a12);
\draw [post,black!30,thin] (a8) to (a15);
\draw [post,black,very thick] (a5) to (a10);
\draw [post,black,very thick] (a5) to (a16);
\draw [post,black,very thick] (a2) to (a10);
\draw [post,black,very thick] (a2) to (a13);
\draw [post,black,very thick] (a8) to (a13);
\draw [post,black,very thick] (a8) to (a16);
\end{tikzpicture}
}

\def\TubeDeltaCCA{
\begin{tikzpicture}
[scale=0.60,
inner sep=1.5,
vertex/.style={circle,draw=black!50,fill=black!20, very thick},
parta/.style={circle,draw=black!50,fill=black!0, very thick},
partb/.style={circle,draw=black!50,fill=black!100, very thick},
edgevertex/.style={circle,draw=blue!80,fill=blue!40, very thick},
normal/.style=,
post/.style={->, >=stealth'},
pre/.style={<-, >=stealth'}]
\node[vertex](a2) at (1.111140,1.662939) {};
\node[vertex](a5) at (1.111140,2.662939) {};
\node[vertex](a8) at (1.111140,3.662939) {};
\node[vertex](a11) at (3.411140,1.662939) {};
\node[vertex](a14) at (3.411140,2.662939) {};
\node[vertex](a17) at (3.411140,3.662939) {};
\node[vertex](a1) at (0.555570,0.831470) {};
\node[vertex](a4) at (0.555570,1.831470) {};
\node[vertex](a7) at (0.555570,2.831470) {};
\node[vertex](a10) at (2.855570,0.831470) {};
\node[vertex](a13) at (2.855570,1.831470) {};
\node[vertex](a16) at (2.855570,2.831470) {};
\draw [post,black!30,thin] (a1) to (a10);
\draw [post,black!30,thin] (a1) to (a11);
\draw [post,black!30,thin] (a1) to (a13);
\draw [post,black!30,thin] (a1) to (a14);
\draw [post,black!30,thin] (a2) to (a10);
\draw [post,black!30,thin] (a2) to (a13);
\draw [post,black!30,thin] (a4) to (a10);
\draw [post,black!30,thin] (a4) to (a11);
\draw [post,black!30,thin] (a4) to (a16);
\draw [post,black!30,thin] (a4) to (a17);
\draw [post,black!30,thin] (a5) to (a10);
\draw [post,black!30,thin] (a5) to (a16);
\draw [post,black!30,thin] (a7) to (a13);
\draw [post,black!30,thin] (a7) to (a14);
\draw [post,black!30,thin] (a7) to (a16);
\draw [post,black!30,thin] (a7) to (a17);
\draw [post,black!30,thin] (a8) to (a13);
\draw [post,black!30,thin] (a8) to (a16);
\node[vertex](a0) at (0.000000,0.000000) {};
\node[vertex](a3) at (0.000000,1.000000) {};
\node[vertex](a6) at (0.000000,2.000000) {};
\node[vertex](a9) at (2.300000,0.000000) {};
\node[vertex](a12) at (2.300000,1.000000) {};
\node[vertex](a15) at (2.300000,2.000000) {};
\draw [post,black!30,thin] (a0) to (a9);
\draw [post,black!30,thin] (a0) to (a10);
\draw [post,black!30,thin] (a0) to (a11);
\draw [post,black!30,thin] (a0) to (a12);
\draw [post,black!30,thin] (a0) to (a13);
\draw [post,black!30,thin] (a0) to (a14);
\draw [post,black!30,thin] (a1) to (a9);
\draw [post,black!30,thin] (a1) to (a12);
\draw [post,black!30,thin] (a2) to (a9);
\draw [post,black!30,thin] (a2) to (a12);
\draw [post,black!30,thin] (a3) to (a9);
\draw [post,black!30,thin] (a3) to (a10);
\draw [post,black!30,thin] (a3) to (a11);
\draw [post,black!30,thin] (a3) to (a15);
\draw [post,black!30,thin] (a3) to (a16);
\draw [post,black!30,thin] (a3) to (a17);
\draw [post,black!30,thin] (a4) to (a9);
\draw [post,black!30,thin] (a4) to (a15);
\draw [post,black!30,thin] (a5) to (a9);
\draw [post,black!30,thin] (a5) to (a15);
\draw [post,black!30,thin] (a6) to (a12);
\draw [post,black!30,thin] (a6) to (a13);
\draw [post,black!30,thin] (a6) to (a14);
\draw [post,black!30,thin] (a6) to (a15);
\draw [post,black!30,thin] (a6) to (a16);
\draw [post,black!30,thin] (a6) to (a17);
\draw [post,black!30,thin] (a7) to (a12);
\draw [post,black!30,thin] (a7) to (a15);
\draw [post,black!30,thin] (a8) to (a12);
\draw [post,black!30,thin] (a8) to (a15);
\draw [post,black,very thick] (a2) to (a11);
\draw [post,black,very thick] (a2) to (a14);
\draw [post,black,very thick] (a5) to (a11);
\draw [post,black,very thick] (a5) to (a17);
\draw [post,black,very thick] (a8) to (a14);
\draw [post,black,very thick] (a8) to (a17);
\end{tikzpicture}
}

\def\TubeDeltaAAB{
\begin{tikzpicture}
[scale=0.60,
inner sep=1.5,
vertex/.style={circle,draw=black!50,fill=black!20, very thick},
parta/.style={circle,draw=black!50,fill=black!0, very thick},
partb/.style={circle,draw=black!50,fill=black!100, very thick},
edgevertex/.style={circle,draw=blue!80,fill=blue!40, very thick},
normal/.style=,
post/.style={->, >=stealth'},
pre/.style={<-, >=stealth'}]
\node[vertex](a2) at (1.111140,1.662939) {};
\node[vertex](a5) at (1.111140,2.662939) {};
\node[vertex](a8) at (1.111140,3.662939) {};
\node[vertex](a11) at (3.411140,1.662939) {};
\node[vertex](a14) at (3.411140,2.662939) {};
\node[vertex](a17) at (3.411140,3.662939) {};
\draw [post,black!30,thin] (a2) to (a14);
\draw [post,black!30,thin] (a5) to (a11);
\draw [post,black!30,thin] (a5) to (a17);
\draw [post,black!30,thin] (a8) to (a14);
\draw [post,black!30,thin] (a8) to (a17);
\node[vertex](a1) at (0.555570,0.831470) {};
\node[vertex](a4) at (0.555570,1.831470) {};
\node[vertex](a7) at (0.555570,2.831470) {};
\node[vertex](a10) at (2.855570,0.831470) {};
\node[vertex](a13) at (2.855570,1.831470) {};
\node[vertex](a16) at (2.855570,2.831470) {};
\draw [post,black!30,thin] (a1) to (a13);
\draw [post,black!30,thin] (a1) to (a14);
\draw [post,black!30,thin] (a2) to (a13);
\draw [post,black!30,thin] (a4) to (a10);
\draw [post,black!30,thin] (a4) to (a11);
\draw [post,black!30,thin] (a4) to (a16);
\draw [post,black!30,thin] (a4) to (a17);
\draw [post,black!30,thin] (a5) to (a10);
\draw [post,black!30,thin] (a5) to (a16);
\draw [post,black!30,thin] (a7) to (a13);
\draw [post,black!30,thin] (a7) to (a14);
\draw [post,black!30,thin] (a7) to (a16);
\draw [post,black!30,thin] (a7) to (a17);
\draw [post,black!30,thin] (a8) to (a13);
\draw [post,black!30,thin] (a8) to (a16);
\node[vertex](a0) at (0.000000,0.000000) {};
\node[vertex](a3) at (0.000000,1.000000) {};
\node[vertex](a6) at (0.000000,2.000000) {};
\node[vertex](a9) at (2.300000,0.000000) {};
\node[vertex](a12) at (2.300000,1.000000) {};
\node[vertex](a15) at (2.300000,2.000000) {};
\draw [post,black!30,thin] (a0) to (a12);
\draw [post,black!30,thin] (a0) to (a13);
\draw [post,black!30,thin] (a0) to (a14);
\draw [post,black!30,thin] (a1) to (a12);
\draw [post,black!30,thin] (a2) to (a12);
\draw [post,black!30,thin] (a3) to (a9);
\draw [post,black!30,thin] (a3) to (a10);
\draw [post,black!30,thin] (a3) to (a11);
\draw [post,black!30,thin] (a3) to (a15);
\draw [post,black!30,thin] (a3) to (a16);
\draw [post,black!30,thin] (a3) to (a17);
\draw [post,black!30,thin] (a4) to (a9);
\draw [post,black!30,thin] (a4) to (a15);
\draw [post,black!30,thin] (a5) to (a9);
\draw [post,black!30,thin] (a5) to (a15);
\draw [post,black!30,thin] (a6) to (a12);
\draw [post,black!30,thin] (a6) to (a13);
\draw [post,black!30,thin] (a6) to (a14);
\draw [post,black!30,thin] (a6) to (a15);
\draw [post,black!30,thin] (a6) to (a16);
\draw [post,black!30,thin] (a6) to (a17);
\draw [post,black!30,thin] (a7) to (a12);
\draw [post,black!30,thin] (a7) to (a15);
\draw [post,black!30,thin] (a8) to (a12);
\draw [post,black!30,thin] (a8) to (a15);
\draw [post,black,very thick] (a2) to (a11);
\draw [post,black,very thick] (a1) to (a10);
\draw [post,black,very thick] (a1) to (a11);
\draw [post,black,very thick] (a2) to (a10);
\draw [post,black,very thick] (a0) to (a9);
\draw [post,black,very thick] (a0) to (a10);
\draw [post,black,very thick] (a0) to (a11);
\draw [post,black,very thick] (a1) to (a9);
\draw [post,black,very thick] (a2) to (a9);
\end{tikzpicture}
}

\def\TubeDeltaABB{
\begin{tikzpicture}
[scale=0.60,
inner sep=1.5,
vertex/.style={circle,draw=black!50,fill=black!20, very thick},
parta/.style={circle,draw=black!50,fill=black!0, very thick},
partb/.style={circle,draw=black!50,fill=black!100, very thick},
edgevertex/.style={circle,draw=blue!80,fill=blue!40, very thick},
normal/.style=,
post/.style={->, >=stealth'},
pre/.style={<-, >=stealth'}]
\node[vertex](a2) at (1.111140,1.662939) {};
\node[vertex](a5) at (1.111140,2.662939) {};
\node[vertex](a8) at (1.111140,3.662939) {};
\node[vertex](a11) at (3.411140,1.662939) {};
\node[vertex](a14) at (3.411140,2.662939) {};
\node[vertex](a17) at (3.411140,3.662939) {};
\draw [post,black!30,thin] (a2) to (a11);
\draw [post,black!30,thin] (a5) to (a11);
\draw [post,black!30,thin] (a5) to (a17);
\draw [post,black!30,thin] (a8) to (a14);
\draw [post,black!30,thin] (a8) to (a17);
\node[vertex](a1) at (0.555570,0.831470) {};
\node[vertex](a4) at (0.555570,1.831470) {};
\node[vertex](a7) at (0.555570,2.831470) {};
\node[vertex](a10) at (2.855570,0.831470) {};
\node[vertex](a13) at (2.855570,1.831470) {};
\node[vertex](a16) at (2.855570,2.831470) {};
\draw [post,black!30,thin] (a1) to (a10);
\draw [post,black!30,thin] (a1) to (a11);
\draw [post,black!30,thin] (a2) to (a10);
\draw [post,black!30,thin] (a4) to (a10);
\draw [post,black!30,thin] (a4) to (a11);
\draw [post,black!30,thin] (a4) to (a16);
\draw [post,black!30,thin] (a4) to (a17);
\draw [post,black!30,thin] (a5) to (a10);
\draw [post,black!30,thin] (a5) to (a16);
\draw [post,black!30,thin] (a7) to (a13);
\draw [post,black!30,thin] (a7) to (a14);
\draw [post,black!30,thin] (a7) to (a16);
\draw [post,black!30,thin] (a7) to (a17);
\draw [post,black!30,thin] (a8) to (a13);
\draw [post,black!30,thin] (a8) to (a16);
\node[vertex](a0) at (0.000000,0.000000) {};
\node[vertex](a3) at (0.000000,1.000000) {};
\node[vertex](a6) at (0.000000,2.000000) {};
\node[vertex](a9) at (2.300000,0.000000) {};
\node[vertex](a12) at (2.300000,1.000000) {};
\node[vertex](a15) at (2.300000,2.000000) {};
\draw [post,black!30,thin] (a0) to (a9);
\draw [post,black!30,thin] (a0) to (a10);
\draw [post,black!30,thin] (a0) to (a11);
\draw [post,black!30,thin] (a1) to (a9);
\draw [post,black!30,thin] (a2) to (a9);
\draw [post,black!30,thin] (a3) to (a9);
\draw [post,black!30,thin] (a3) to (a10);
\draw [post,black!30,thin] (a3) to (a11);
\draw [post,black!30,thin] (a3) to (a15);
\draw [post,black!30,thin] (a3) to (a16);
\draw [post,black!30,thin] (a3) to (a17);
\draw [post,black!30,thin] (a4) to (a9);
\draw [post,black!30,thin] (a4) to (a15);
\draw [post,black!30,thin] (a5) to (a9);
\draw [post,black!30,thin] (a5) to (a15);
\draw [post,black!30,thin] (a6) to (a12);
\draw [post,black!30,thin] (a6) to (a13);
\draw [post,black!30,thin] (a6) to (a14);
\draw [post,black!30,thin] (a6) to (a15);
\draw [post,black!30,thin] (a6) to (a16);
\draw [post,black!30,thin] (a6) to (a17);
\draw [post,black!30,thin] (a7) to (a12);
\draw [post,black!30,thin] (a7) to (a15);
\draw [post,black!30,thin] (a8) to (a12);
\draw [post,black!30,thin] (a8) to (a15);
\draw [post,black,very thick] (a2) to (a14);
\draw [post,black,very thick] (a1) to (a13);
\draw [post,black,very thick] (a1) to (a14);
\draw [post,black,very thick] (a2) to (a13);
\draw [post,black,very thick] (a0) to (a12);
\draw [post,black,very thick] (a0) to (a13);
\draw [post,black,very thick] (a0) to (a14);
\draw [post,black,very thick] (a1) to (a12);
\draw [post,black,very thick] (a2) to (a12);
\end{tikzpicture}
}

\def\TubeDeltaACB{
\begin{tikzpicture}
[scale=0.60,
inner sep=1.5,
vertex/.style={circle,draw=black!50,fill=black!20, very thick},
parta/.style={circle,draw=black!50,fill=black!0, very thick},
partb/.style={circle,draw=black!50,fill=black!100, very thick},
edgevertex/.style={circle,draw=blue!80,fill=blue!40, very thick},
normal/.style=,
post/.style={->, >=stealth'},
pre/.style={<-, >=stealth'}]
\node[vertex](a2) at (1.111140,1.662939) {};
\node[vertex](a5) at (1.111140,2.662939) {};
\node[vertex](a8) at (1.111140,3.662939) {};
\node[vertex](a11) at (3.411140,1.662939) {};
\node[vertex](a14) at (3.411140,2.662939) {};
\node[vertex](a17) at (3.411140,3.662939) {};
\draw [post,black!30,thin] (a2) to (a11);
\draw [post,black!30,thin] (a2) to (a14);
\draw [post,black!30,thin] (a5) to (a11);
\draw [post,black!30,thin] (a5) to (a17);
\draw [post,black!30,thin] (a8) to (a14);
\draw [post,black!30,thin] (a8) to (a17);
\node[vertex](a1) at (0.555570,0.831470) {};
\node[vertex](a4) at (0.555570,1.831470) {};
\node[vertex](a7) at (0.555570,2.831470) {};
\node[vertex](a10) at (2.855570,0.831470) {};
\node[vertex](a13) at (2.855570,1.831470) {};
\node[vertex](a16) at (2.855570,2.831470) {};
\draw [post,black!30,thin] (a1) to (a10);
\draw [post,black!30,thin] (a1) to (a11);
\draw [post,black!30,thin] (a1) to (a13);
\draw [post,black!30,thin] (a1) to (a14);
\draw [post,black!30,thin] (a2) to (a10);
\draw [post,black!30,thin] (a2) to (a13);
\draw [post,black!30,thin] (a4) to (a10);
\draw [post,black!30,thin] (a4) to (a11);
\draw [post,black!30,thin] (a4) to (a16);
\draw [post,black!30,thin] (a4) to (a17);
\draw [post,black!30,thin] (a5) to (a10);
\draw [post,black!30,thin] (a5) to (a16);
\draw [post,black!30,thin] (a7) to (a13);
\draw [post,black!30,thin] (a7) to (a14);
\draw [post,black!30,thin] (a7) to (a16);
\draw [post,black!30,thin] (a7) to (a17);
\draw [post,black!30,thin] (a8) to (a13);
\draw [post,black!30,thin] (a8) to (a16);
\node[vertex](a0) at (0.000000,0.000000) {};
\node[vertex](a3) at (0.000000,1.000000) {};
\node[vertex](a6) at (0.000000,2.000000) {};
\node[vertex](a9) at (2.300000,0.000000) {};
\node[vertex](a12) at (2.300000,1.000000) {};
\node[vertex](a15) at (2.300000,2.000000) {};
\draw [post,black!30,thin] (a0) to (a9);
\draw [post,black!30,thin] (a0) to (a10);
\draw [post,black!30,thin] (a0) to (a11);
\draw [post,black!30,thin] (a0) to (a12);
\draw [post,black!30,thin] (a0) to (a13);
\draw [post,black!30,thin] (a0) to (a14);
\draw [post,black!30,thin] (a1) to (a9);
\draw [post,black!30,thin] (a1) to (a12);
\draw [post,black!30,thin] (a2) to (a9);
\draw [post,black!30,thin] (a2) to (a12);
\draw [post,black!30,thin] (a3) to (a9);
\draw [post,black!30,thin] (a3) to (a10);
\draw [post,black!30,thin] (a3) to (a11);
\draw [post,black!30,thin] (a3) to (a15);
\draw [post,black!30,thin] (a3) to (a16);
\draw [post,black!30,thin] (a3) to (a17);
\draw [post,black!30,thin] (a4) to (a9);
\draw [post,black!30,thin] (a4) to (a15);
\draw [post,black!30,thin] (a5) to (a9);
\draw [post,black!30,thin] (a5) to (a15);
\draw [post,black!30,thin] (a6) to (a12);
\draw [post,black!30,thin] (a6) to (a13);
\draw [post,black!30,thin] (a6) to (a14);
\draw [post,black!30,thin] (a6) to (a15);
\draw [post,black!30,thin] (a6) to (a16);
\draw [post,black!30,thin] (a6) to (a17);
\draw [post,black!30,thin] (a7) to (a12);
\draw [post,black!30,thin] (a7) to (a15);
\draw [post,black!30,thin] (a8) to (a12);
\draw [post,black!30,thin] (a8) to (a15);
\end{tikzpicture}
}

\def\TubeDeltaBAB{
\begin{tikzpicture}
[scale=0.60,
inner sep=1.5,
vertex/.style={circle,draw=black!50,fill=black!20, very thick},
parta/.style={circle,draw=black!50,fill=black!0, very thick},
partb/.style={circle,draw=black!50,fill=black!100, very thick},
edgevertex/.style={circle,draw=blue!80,fill=blue!40, very thick},
normal/.style=,
post/.style={->, >=stealth'},
pre/.style={<-, >=stealth'}]
\node[vertex](a2) at (1.111140,1.662939) {};
\node[vertex](a5) at (1.111140,2.662939) {};
\node[vertex](a8) at (1.111140,3.662939) {};
\node[vertex](a11) at (3.411140,1.662939) {};
\node[vertex](a14) at (3.411140,2.662939) {};
\node[vertex](a17) at (3.411140,3.662939) {};
\draw [post,black!30,thin] (a2) to (a11);
\draw [post,black!30,thin] (a2) to (a14);
\draw [post,black!30,thin] (a5) to (a17);
\draw [post,black!30,thin] (a8) to (a14);
\draw [post,black!30,thin] (a8) to (a17);
\node[vertex](a1) at (0.555570,0.831470) {};
\node[vertex](a4) at (0.555570,1.831470) {};
\node[vertex](a7) at (0.555570,2.831470) {};
\node[vertex](a10) at (2.855570,0.831470) {};
\node[vertex](a13) at (2.855570,1.831470) {};
\node[vertex](a16) at (2.855570,2.831470) {};
\draw [post,black!30,thin] (a1) to (a10);
\draw [post,black!30,thin] (a1) to (a11);
\draw [post,black!30,thin] (a1) to (a13);
\draw [post,black!30,thin] (a1) to (a14);
\draw [post,black!30,thin] (a2) to (a10);
\draw [post,black!30,thin] (a2) to (a13);
\draw [post,black!30,thin] (a4) to (a16);
\draw [post,black!30,thin] (a4) to (a17);
\draw [post,black!30,thin] (a5) to (a16);
\draw [post,black!30,thin] (a7) to (a13);
\draw [post,black!30,thin] (a7) to (a14);
\draw [post,black!30,thin] (a7) to (a16);
\draw [post,black!30,thin] (a7) to (a17);
\draw [post,black!30,thin] (a8) to (a13);
\draw [post,black!30,thin] (a8) to (a16);
\node[vertex](a0) at (0.000000,0.000000) {};
\node[vertex](a3) at (0.000000,1.000000) {};
\node[vertex](a6) at (0.000000,2.000000) {};
\node[vertex](a9) at (2.300000,0.000000) {};
\node[vertex](a12) at (2.300000,1.000000) {};
\node[vertex](a15) at (2.300000,2.000000) {};
\draw [post,black!30,thin] (a0) to (a9);
\draw [post,black!30,thin] (a0) to (a10);
\draw [post,black!30,thin] (a0) to (a11);
\draw [post,black!30,thin] (a0) to (a12);
\draw [post,black!30,thin] (a0) to (a13);
\draw [post,black!30,thin] (a0) to (a14);
\draw [post,black!30,thin] (a1) to (a9);
\draw [post,black!30,thin] (a1) to (a12);
\draw [post,black!30,thin] (a2) to (a9);
\draw [post,black!30,thin] (a2) to (a12);
\draw [post,black!30,thin] (a3) to (a15);
\draw [post,black!30,thin] (a3) to (a16);
\draw [post,black!30,thin] (a3) to (a17);
\draw [post,black!30,thin] (a4) to (a15);
\draw [post,black!30,thin] (a5) to (a15);
\draw [post,black!30,thin] (a6) to (a12);
\draw [post,black!30,thin] (a6) to (a13);
\draw [post,black!30,thin] (a6) to (a14);
\draw [post,black!30,thin] (a6) to (a15);
\draw [post,black!30,thin] (a6) to (a16);
\draw [post,black!30,thin] (a6) to (a17);
\draw [post,black!30,thin] (a7) to (a12);
\draw [post,black!30,thin] (a7) to (a15);
\draw [post,black!30,thin] (a8) to (a12);
\draw [post,black!30,thin] (a8) to (a15);
\draw [post,black,very thick] (a5) to (a9);
\draw [post,black,very thick] (a5) to (a10);
\draw [post,black,very thick] (a4) to (a9);
\draw [post,black,very thick] (a5) to (a11);
\draw [post,black,very thick] (a4) to (a10);
\draw [post,black,very thick] (a4) to (a11);
\draw [post,black,very thick] (a3) to (a9);
\draw [post,black,very thick] (a3) to (a10);
\draw [post,black,very thick] (a3) to (a11);
\end{tikzpicture}
}

\def\TubeDeltaBBB{
\begin{tikzpicture}
[scale=0.60,
inner sep=1.5,
vertex/.style={circle,draw=black!50,fill=black!20, very thick},
parta/.style={circle,draw=black!50,fill=black!0, very thick},
partb/.style={circle,draw=black!50,fill=black!100, very thick},
edgevertex/.style={circle,draw=blue!80,fill=blue!40, very thick},
normal/.style=,
post/.style={->, >=stealth'},
pre/.style={<-, >=stealth'}]
\node[vertex](a2) at (1.111140,1.662939) {};
\node[vertex](a5) at (1.111140,2.662939) {};
\node[vertex](a8) at (1.111140,3.662939) {};
\node[vertex](a11) at (3.411140,1.662939) {};
\node[vertex](a14) at (3.411140,2.662939) {};
\node[vertex](a17) at (3.411140,3.662939) {};
\draw [post,black!30,thin] (a2) to (a11);
\draw [post,black!30,thin] (a2) to (a14);
\draw [post,black!30,thin] (a5) to (a11);
\draw [post,black!30,thin] (a5) to (a17);
\draw [post,black!30,thin] (a8) to (a14);
\draw [post,black!30,thin] (a8) to (a17);
\node[vertex](a1) at (0.555570,0.831470) {};
\node[vertex](a4) at (0.555570,1.831470) {};
\node[vertex](a7) at (0.555570,2.831470) {};
\node[vertex](a10) at (2.855570,0.831470) {};
\node[vertex](a13) at (2.855570,1.831470) {};
\node[vertex](a16) at (2.855570,2.831470) {};
\draw [post,black!30,thin] (a1) to (a10);
\draw [post,black!30,thin] (a1) to (a11);
\draw [post,black!30,thin] (a1) to (a13);
\draw [post,black!30,thin] (a1) to (a14);
\draw [post,black!30,thin] (a2) to (a10);
\draw [post,black!30,thin] (a2) to (a13);
\draw [post,black!30,thin] (a4) to (a10);
\draw [post,black!30,thin] (a4) to (a11);
\draw [post,black!30,thin] (a4) to (a16);
\draw [post,black!30,thin] (a4) to (a17);
\draw [post,black!30,thin] (a5) to (a10);
\draw [post,black!30,thin] (a5) to (a16);
\draw [post,black!30,thin] (a7) to (a13);
\draw [post,black!30,thin] (a7) to (a14);
\draw [post,black!30,thin] (a7) to (a16);
\draw [post,black!30,thin] (a7) to (a17);
\draw [post,black!30,thin] (a8) to (a13);
\draw [post,black!30,thin] (a8) to (a16);
\node[vertex](a0) at (0.000000,0.000000) {};
\node[vertex](a3) at (0.000000,1.000000) {};
\node[vertex](a6) at (0.000000,2.000000) {};
\node[vertex](a9) at (2.300000,0.000000) {};
\node[vertex](a12) at (2.300000,1.000000) {};
\node[vertex](a15) at (2.300000,2.000000) {};
\draw [post,black!30,thin] (a0) to (a9);
\draw [post,black!30,thin] (a0) to (a10);
\draw [post,black!30,thin] (a0) to (a11);
\draw [post,black!30,thin] (a0) to (a12);
\draw [post,black!30,thin] (a0) to (a13);
\draw [post,black!30,thin] (a0) to (a14);
\draw [post,black!30,thin] (a1) to (a9);
\draw [post,black!30,thin] (a1) to (a12);
\draw [post,black!30,thin] (a2) to (a9);
\draw [post,black!30,thin] (a2) to (a12);
\draw [post,black!30,thin] (a3) to (a9);
\draw [post,black!30,thin] (a3) to (a10);
\draw [post,black!30,thin] (a3) to (a11);
\draw [post,black!30,thin] (a3) to (a15);
\draw [post,black!30,thin] (a3) to (a16);
\draw [post,black!30,thin] (a3) to (a17);
\draw [post,black!30,thin] (a4) to (a9);
\draw [post,black!30,thin] (a4) to (a15);
\draw [post,black!30,thin] (a5) to (a9);
\draw [post,black!30,thin] (a5) to (a15);
\draw [post,black!30,thin] (a6) to (a12);
\draw [post,black!30,thin] (a6) to (a13);
\draw [post,black!30,thin] (a6) to (a14);
\draw [post,black!30,thin] (a6) to (a15);
\draw [post,black!30,thin] (a6) to (a16);
\draw [post,black!30,thin] (a6) to (a17);
\draw [post,black!30,thin] (a7) to (a12);
\draw [post,black!30,thin] (a7) to (a15);
\draw [post,black!30,thin] (a8) to (a12);
\draw [post,black!30,thin] (a8) to (a15);
\end{tikzpicture}
}

\def\TubeDeltaBCB{
\begin{tikzpicture}
[scale=0.60,
inner sep=1.5,
vertex/.style={circle,draw=black!50,fill=black!20, very thick},
parta/.style={circle,draw=black!50,fill=black!0, very thick},
partb/.style={circle,draw=black!50,fill=black!100, very thick},
edgevertex/.style={circle,draw=blue!80,fill=blue!40, very thick},
normal/.style=,
post/.style={->, >=stealth'},
pre/.style={<-, >=stealth'}]
\node[vertex](a2) at (1.111140,1.662939) {};
\node[vertex](a5) at (1.111140,2.662939) {};
\node[vertex](a8) at (1.111140,3.662939) {};
\node[vertex](a11) at (3.411140,1.662939) {};
\node[vertex](a14) at (3.411140,2.662939) {};
\node[vertex](a17) at (3.411140,3.662939) {};
\draw [post,black!30,thin] (a2) to (a11);
\draw [post,black!30,thin] (a2) to (a14);
\draw [post,black!30,thin] (a5) to (a11);
\draw [post,black!30,thin] (a8) to (a14);
\draw [post,black!30,thin] (a8) to (a17);
\node[vertex](a1) at (0.555570,0.831470) {};
\node[vertex](a4) at (0.555570,1.831470) {};
\node[vertex](a7) at (0.555570,2.831470) {};
\node[vertex](a10) at (2.855570,0.831470) {};
\node[vertex](a13) at (2.855570,1.831470) {};
\node[vertex](a16) at (2.855570,2.831470) {};
\draw [post,black!30,thin] (a1) to (a10);
\draw [post,black!30,thin] (a1) to (a11);
\draw [post,black!30,thin] (a1) to (a13);
\draw [post,black!30,thin] (a1) to (a14);
\draw [post,black!30,thin] (a2) to (a10);
\draw [post,black!30,thin] (a2) to (a13);
\draw [post,black!30,thin] (a4) to (a10);
\draw [post,black!30,thin] (a4) to (a11);
\draw [post,black!30,thin] (a5) to (a10);
\draw [post,black!30,thin] (a7) to (a13);
\draw [post,black!30,thin] (a7) to (a14);
\draw [post,black!30,thin] (a7) to (a16);
\draw [post,black!30,thin] (a7) to (a17);
\draw [post,black!30,thin] (a8) to (a13);
\draw [post,black!30,thin] (a8) to (a16);
\node[vertex](a0) at (0.000000,0.000000) {};
\node[vertex](a3) at (0.000000,1.000000) {};
\node[vertex](a6) at (0.000000,2.000000) {};
\node[vertex](a9) at (2.300000,0.000000) {};
\node[vertex](a12) at (2.300000,1.000000) {};
\node[vertex](a15) at (2.300000,2.000000) {};
\draw [post,black!30,thin] (a0) to (a9);
\draw [post,black!30,thin] (a0) to (a10);
\draw [post,black!30,thin] (a0) to (a11);
\draw [post,black!30,thin] (a0) to (a12);
\draw [post,black!30,thin] (a0) to (a13);
\draw [post,black!30,thin] (a0) to (a14);
\draw [post,black!30,thin] (a1) to (a9);
\draw [post,black!30,thin] (a1) to (a12);
\draw [post,black!30,thin] (a2) to (a9);
\draw [post,black!30,thin] (a2) to (a12);
\draw [post,black!30,thin] (a3) to (a9);
\draw [post,black!30,thin] (a3) to (a10);
\draw [post,black!30,thin] (a3) to (a11);
\draw [post,black!30,thin] (a4) to (a9);
\draw [post,black!30,thin] (a5) to (a9);
\draw [post,black!30,thin] (a6) to (a12);
\draw [post,black!30,thin] (a6) to (a13);
\draw [post,black!30,thin] (a6) to (a14);
\draw [post,black!30,thin] (a6) to (a15);
\draw [post,black!30,thin] (a6) to (a16);
\draw [post,black!30,thin] (a6) to (a17);
\draw [post,black!30,thin] (a7) to (a12);
\draw [post,black!30,thin] (a7) to (a15);
\draw [post,black!30,thin] (a8) to (a12);
\draw [post,black!30,thin] (a8) to (a15);
\draw [post,black,very thick] (a5) to (a17);
\draw [post,black,very thick] (a4) to (a16);
\draw [post,black,very thick] (a4) to (a17);
\draw [post,black,very thick] (a5) to (a16);
\draw [post,black,very thick] (a3) to (a15);
\draw [post,black,very thick] (a3) to (a16);
\draw [post,black,very thick] (a3) to (a17);
\draw [post,black,very thick] (a4) to (a15);
\draw [post,black,very thick] (a5) to (a15);
\end{tikzpicture}
}

\def\TubeDeltaCAB{
\begin{tikzpicture}
[scale=0.60,
inner sep=1.5,
vertex/.style={circle,draw=black!50,fill=black!20, very thick},
parta/.style={circle,draw=black!50,fill=black!0, very thick},
partb/.style={circle,draw=black!50,fill=black!100, very thick},
edgevertex/.style={circle,draw=blue!80,fill=blue!40, very thick},
normal/.style=,
post/.style={->, >=stealth'},
pre/.style={<-, >=stealth'}]
\node[vertex](a2) at (1.111140,1.662939) {};
\node[vertex](a5) at (1.111140,2.662939) {};
\node[vertex](a8) at (1.111140,3.662939) {};
\node[vertex](a11) at (3.411140,1.662939) {};
\node[vertex](a14) at (3.411140,2.662939) {};
\node[vertex](a17) at (3.411140,3.662939) {};
\draw [post,black!30,thin] (a2) to (a11);
\draw [post,black!30,thin] (a2) to (a14);
\draw [post,black!30,thin] (a5) to (a11);
\draw [post,black!30,thin] (a5) to (a17);
\draw [post,black!30,thin] (a8) to (a14);
\draw [post,black!30,thin] (a8) to (a17);
\node[vertex](a1) at (0.555570,0.831470) {};
\node[vertex](a4) at (0.555570,1.831470) {};
\node[vertex](a7) at (0.555570,2.831470) {};
\node[vertex](a10) at (2.855570,0.831470) {};
\node[vertex](a13) at (2.855570,1.831470) {};
\node[vertex](a16) at (2.855570,2.831470) {};
\draw [post,black!30,thin] (a1) to (a10);
\draw [post,black!30,thin] (a1) to (a11);
\draw [post,black!30,thin] (a1) to (a13);
\draw [post,black!30,thin] (a1) to (a14);
\draw [post,black!30,thin] (a2) to (a10);
\draw [post,black!30,thin] (a2) to (a13);
\draw [post,black!30,thin] (a4) to (a10);
\draw [post,black!30,thin] (a4) to (a11);
\draw [post,black!30,thin] (a4) to (a16);
\draw [post,black!30,thin] (a4) to (a17);
\draw [post,black!30,thin] (a5) to (a10);
\draw [post,black!30,thin] (a5) to (a16);
\draw [post,black!30,thin] (a7) to (a13);
\draw [post,black!30,thin] (a7) to (a14);
\draw [post,black!30,thin] (a7) to (a16);
\draw [post,black!30,thin] (a7) to (a17);
\draw [post,black!30,thin] (a8) to (a13);
\draw [post,black!30,thin] (a8) to (a16);
\node[vertex](a0) at (0.000000,0.000000) {};
\node[vertex](a3) at (0.000000,1.000000) {};
\node[vertex](a6) at (0.000000,2.000000) {};
\node[vertex](a9) at (2.300000,0.000000) {};
\node[vertex](a12) at (2.300000,1.000000) {};
\node[vertex](a15) at (2.300000,2.000000) {};
\draw [post,black!30,thin] (a0) to (a9);
\draw [post,black!30,thin] (a0) to (a10);
\draw [post,black!30,thin] (a0) to (a11);
\draw [post,black!30,thin] (a0) to (a12);
\draw [post,black!30,thin] (a0) to (a13);
\draw [post,black!30,thin] (a0) to (a14);
\draw [post,black!30,thin] (a1) to (a9);
\draw [post,black!30,thin] (a1) to (a12);
\draw [post,black!30,thin] (a2) to (a9);
\draw [post,black!30,thin] (a2) to (a12);
\draw [post,black!30,thin] (a3) to (a9);
\draw [post,black!30,thin] (a3) to (a10);
\draw [post,black!30,thin] (a3) to (a11);
\draw [post,black!30,thin] (a3) to (a15);
\draw [post,black!30,thin] (a3) to (a16);
\draw [post,black!30,thin] (a3) to (a17);
\draw [post,black!30,thin] (a4) to (a9);
\draw [post,black!30,thin] (a4) to (a15);
\draw [post,black!30,thin] (a5) to (a9);
\draw [post,black!30,thin] (a5) to (a15);
\draw [post,black!30,thin] (a6) to (a12);
\draw [post,black!30,thin] (a6) to (a13);
\draw [post,black!30,thin] (a6) to (a14);
\draw [post,black!30,thin] (a6) to (a15);
\draw [post,black!30,thin] (a6) to (a16);
\draw [post,black!30,thin] (a6) to (a17);
\draw [post,black!30,thin] (a7) to (a12);
\draw [post,black!30,thin] (a7) to (a15);
\draw [post,black!30,thin] (a8) to (a12);
\draw [post,black!30,thin] (a8) to (a15);
\end{tikzpicture}
}

\def\TubeDeltaCBB{
\begin{tikzpicture}
[scale=0.60,
inner sep=1.5,
vertex/.style={circle,draw=black!50,fill=black!20, very thick},
parta/.style={circle,draw=black!50,fill=black!0, very thick},
partb/.style={circle,draw=black!50,fill=black!100, very thick},
edgevertex/.style={circle,draw=blue!80,fill=blue!40, very thick},
normal/.style=,
post/.style={->, >=stealth'},
pre/.style={<-, >=stealth'}]
\node[vertex](a2) at (1.111140,1.662939) {};
\node[vertex](a5) at (1.111140,2.662939) {};
\node[vertex](a8) at (1.111140,3.662939) {};
\node[vertex](a11) at (3.411140,1.662939) {};
\node[vertex](a14) at (3.411140,2.662939) {};
\node[vertex](a17) at (3.411140,3.662939) {};
\draw [post,black!30,thin] (a2) to (a11);
\draw [post,black!30,thin] (a2) to (a14);
\draw [post,black!30,thin] (a5) to (a11);
\draw [post,black!30,thin] (a5) to (a17);
\draw [post,black!30,thin] (a8) to (a17);
\node[vertex](a1) at (0.555570,0.831470) {};
\node[vertex](a4) at (0.555570,1.831470) {};
\node[vertex](a7) at (0.555570,2.831470) {};
\node[vertex](a10) at (2.855570,0.831470) {};
\node[vertex](a13) at (2.855570,1.831470) {};
\node[vertex](a16) at (2.855570,2.831470) {};
\draw [post,black!30,thin] (a1) to (a10);
\draw [post,black!30,thin] (a1) to (a11);
\draw [post,black!30,thin] (a1) to (a13);
\draw [post,black!30,thin] (a1) to (a14);
\draw [post,black!30,thin] (a2) to (a10);
\draw [post,black!30,thin] (a2) to (a13);
\draw [post,black!30,thin] (a4) to (a10);
\draw [post,black!30,thin] (a4) to (a11);
\draw [post,black!30,thin] (a4) to (a16);
\draw [post,black!30,thin] (a4) to (a17);
\draw [post,black!30,thin] (a5) to (a10);
\draw [post,black!30,thin] (a5) to (a16);
\draw [post,black!30,thin] (a7) to (a16);
\draw [post,black!30,thin] (a7) to (a17);
\draw [post,black!30,thin] (a8) to (a16);
\node[vertex](a0) at (0.000000,0.000000) {};
\node[vertex](a3) at (0.000000,1.000000) {};
\node[vertex](a6) at (0.000000,2.000000) {};
\node[vertex](a9) at (2.300000,0.000000) {};
\node[vertex](a12) at (2.300000,1.000000) {};
\node[vertex](a15) at (2.300000,2.000000) {};
\draw [post,black!30,thin] (a0) to (a9);
\draw [post,black!30,thin] (a0) to (a10);
\draw [post,black!30,thin] (a0) to (a11);
\draw [post,black!30,thin] (a0) to (a12);
\draw [post,black!30,thin] (a0) to (a13);
\draw [post,black!30,thin] (a0) to (a14);
\draw [post,black!30,thin] (a1) to (a9);
\draw [post,black!30,thin] (a1) to (a12);
\draw [post,black!30,thin] (a2) to (a9);
\draw [post,black!30,thin] (a2) to (a12);
\draw [post,black!30,thin] (a3) to (a9);
\draw [post,black!30,thin] (a3) to (a10);
\draw [post,black!30,thin] (a3) to (a11);
\draw [post,black!30,thin] (a3) to (a15);
\draw [post,black!30,thin] (a3) to (a16);
\draw [post,black!30,thin] (a3) to (a17);
\draw [post,black!30,thin] (a4) to (a9);
\draw [post,black!30,thin] (a4) to (a15);
\draw [post,black!30,thin] (a5) to (a9);
\draw [post,black!30,thin] (a5) to (a15);
\draw [post,black!30,thin] (a6) to (a15);
\draw [post,black!30,thin] (a6) to (a16);
\draw [post,black!30,thin] (a6) to (a17);
\draw [post,black!30,thin] (a7) to (a15);
\draw [post,black!30,thin] (a8) to (a15);
\draw [post,black,very thick] (a7) to (a12);
\draw [post,black,very thick] (a8) to (a12);
\draw [post,black,very thick] (a8) to (a14);
\draw [post,black,very thick] (a7) to (a13);
\draw [post,black,very thick] (a7) to (a14);
\draw [post,black,very thick] (a8) to (a13);
\draw [post,black,very thick] (a6) to (a12);
\draw [post,black,very thick] (a6) to (a13);
\draw [post,black,very thick] (a6) to (a14);
\end{tikzpicture}
}

\def\TubeDeltaCCB{
\begin{tikzpicture}
[scale=0.60,
inner sep=1.5,
vertex/.style={circle,draw=black!50,fill=black!20, very thick},
parta/.style={circle,draw=black!50,fill=black!0, very thick},
partb/.style={circle,draw=black!50,fill=black!100, very thick},
edgevertex/.style={circle,draw=blue!80,fill=blue!40, very thick},
normal/.style=,
post/.style={->, >=stealth'},
pre/.style={<-, >=stealth'}]
\node[vertex](a2) at (1.111140,1.662939) {};
\node[vertex](a5) at (1.111140,2.662939) {};
\node[vertex](a8) at (1.111140,3.662939) {};
\node[vertex](a11) at (3.411140,1.662939) {};
\node[vertex](a14) at (3.411140,2.662939) {};
\node[vertex](a17) at (3.411140,3.662939) {};
\draw [post,black!30,thin] (a2) to (a11);
\draw [post,black!30,thin] (a2) to (a14);
\draw [post,black!30,thin] (a5) to (a11);
\draw [post,black!30,thin] (a5) to (a17);
\draw [post,black!30,thin] (a8) to (a14);
\node[vertex](a1) at (0.555570,0.831470) {};
\node[vertex](a4) at (0.555570,1.831470) {};
\node[vertex](a7) at (0.555570,2.831470) {};
\node[vertex](a10) at (2.855570,0.831470) {};
\node[vertex](a13) at (2.855570,1.831470) {};
\node[vertex](a16) at (2.855570,2.831470) {};
\draw [post,black!30,thin] (a1) to (a10);
\draw [post,black!30,thin] (a1) to (a11);
\draw [post,black!30,thin] (a1) to (a13);
\draw [post,black!30,thin] (a1) to (a14);
\draw [post,black!30,thin] (a2) to (a10);
\draw [post,black!30,thin] (a2) to (a13);
\draw [post,black!30,thin] (a4) to (a10);
\draw [post,black!30,thin] (a4) to (a11);
\draw [post,black!30,thin] (a4) to (a16);
\draw [post,black!30,thin] (a4) to (a17);
\draw [post,black!30,thin] (a5) to (a10);
\draw [post,black!30,thin] (a5) to (a16);
\draw [post,black!30,thin] (a7) to (a13);
\draw [post,black!30,thin] (a7) to (a14);
\draw [post,black!30,thin] (a8) to (a13);
\node[vertex](a0) at (0.000000,0.000000) {};
\node[vertex](a3) at (0.000000,1.000000) {};
\node[vertex](a6) at (0.000000,2.000000) {};
\node[vertex](a9) at (2.300000,0.000000) {};
\node[vertex](a12) at (2.300000,1.000000) {};
\node[vertex](a15) at (2.300000,2.000000) {};
\draw [post,black!30,thin] (a0) to (a9);
\draw [post,black!30,thin] (a0) to (a10);
\draw [post,black!30,thin] (a0) to (a11);
\draw [post,black!30,thin] (a0) to (a12);
\draw [post,black!30,thin] (a0) to (a13);
\draw [post,black!30,thin] (a0) to (a14);
\draw [post,black!30,thin] (a1) to (a9);
\draw [post,black!30,thin] (a1) to (a12);
\draw [post,black!30,thin] (a2) to (a9);
\draw [post,black!30,thin] (a2) to (a12);
\draw [post,black!30,thin] (a3) to (a9);
\draw [post,black!30,thin] (a3) to (a10);
\draw [post,black!30,thin] (a3) to (a11);
\draw [post,black!30,thin] (a3) to (a15);
\draw [post,black!30,thin] (a3) to (a16);
\draw [post,black!30,thin] (a3) to (a17);
\draw [post,black!30,thin] (a4) to (a9);
\draw [post,black!30,thin] (a4) to (a15);
\draw [post,black!30,thin] (a5) to (a9);
\draw [post,black!30,thin] (a5) to (a15);
\draw [post,black!30,thin] (a6) to (a12);
\draw [post,black!30,thin] (a6) to (a13);
\draw [post,black!30,thin] (a6) to (a14);
\draw [post,black!30,thin] (a7) to (a12);
\draw [post,black!30,thin] (a8) to (a12);
\draw [post,black,very thick] (a8) to (a17);
\draw [post,black,very thick] (a7) to (a15);
\draw [post,black,very thick] (a8) to (a15);
\draw [post,black,very thick] (a7) to (a16);
\draw [post,black,very thick] (a7) to (a17);
\draw [post,black,very thick] (a8) to (a16);
\draw [post,black,very thick] (a6) to (a15);
\draw [post,black,very thick] (a6) to (a16);
\draw [post,black,very thick] (a6) to (a17);
\end{tikzpicture}
}

\def\TubeAZ{
\begin{tikzpicture}
[scale=0.60,
inner sep=1.5,
vertex/.style={circle,draw=black!50,fill=black!20, very thick},
parta/.style={circle,draw=black!50,fill=black!0, very thick},
partb/.style={circle,draw=black!50,fill=black!100, very thick},
edgevertex/.style={circle,draw=blue!80,fill=blue!40, very thick},
normal/.style=,
post/.style={->, >=stealth'},
pre/.style={<-, >=stealth'}]
\node(a2) at (1.111140,1.662939) {};
\node(a5) at (1.111140,2.662939) {};
\node(a8) at (1.111140,3.662939) {};
\node[vertex](a11) at (3.411140,1.662939) {};
\node[vertex](a14) at (3.411140,2.662939) {};
\node[vertex](a17) at (3.411140,3.662939) {};
\node[vertex](a20) at (5.711140,1.662939) {};
\node[vertex](a23) at (5.711140,2.662939) {};
\node[vertex](a26) at (5.711140,3.662939) {};
\node[vertex](a29) at (8.011140,1.662939) {};
\node[vertex](a32) at (8.011140,2.662939) {};
\node[vertex](a35) at (8.011140,3.662939) {};
\node[vertex](a38) at (10.311140,1.662939) {};
\node[vertex](a41) at (10.311140,2.662939) {};
\node[vertex](a44) at (10.311140,3.662939) {};
\node[vertex](a47) at (12.611140,1.662939) {};
\node[vertex](a50) at (12.611140,2.662939) {};
\node[vertex](a53) at (12.611140,3.662939) {};
\node[vertex](a56) at (14.911140,1.662939) {};
\node[vertex](a59) at (14.911140,2.662939) {};
\node[vertex](a62) at (14.911140,3.662939) {};
\node[vertex](a65) at (17.211140,1.662939) {};
\node[vertex](a68) at (17.211140,2.662939) {};
\node[vertex](a71) at (17.211140,3.662939) {};
\node(a74) at (19.511140,1.662939) {};
\node(a77) at (19.511140,2.662939) {};
\node(a80) at (19.511140,3.662939) {};
\draw [post,black!30,thin] (a11) to (a20);
\draw [post,black!30,thin] (a11) to (a23);
\draw [post,black!30,thin] (a14) to (a20);
\draw [post,black!30,thin] (a14) to (a26);
\draw [post,black!30,thin] (a17) to (a23);
\draw [post,black!30,thin] (a17) to (a26);
\draw [post,black!30,thin] (a20) to (a29);
\draw [post,black!30,thin] (a20) to (a32);
\draw [post,black!30,thin] (a20) to (a35);
\draw [post,black!30,thin] (a23) to (a29);
\draw [post,black!30,thin] (a23) to (a32);
\draw [post,black!30,thin] (a23) to (a35);
\draw [post,black!30,thin] (a26) to (a29);
\draw [post,black!30,thin] (a26) to (a32);
\draw [post,black!30,thin] (a26) to (a35);
\draw [post,black!30,thin] (a29) to (a38);
\draw [post,black!30,thin] (a29) to (a41);
\draw [post,black!30,thin] (a32) to (a38);
\draw [post,black!30,thin] (a32) to (a44);
\draw [post,black!30,thin] (a35) to (a41);
\draw [post,black!30,thin] (a35) to (a44);
\draw [post,black!30,thin] (a38) to (a47);
\draw [post,black!30,thin] (a38) to (a50);
\draw [post,black!30,thin] (a38) to (a53);
\draw [post,black!30,thin] (a41) to (a47);
\draw [post,black!30,thin] (a41) to (a50);
\draw [post,black!30,thin] (a41) to (a53);
\draw [post,black!30,thin] (a44) to (a47);
\draw [post,black!30,thin] (a44) to (a50);
\draw [post,black!30,thin] (a44) to (a53);
\draw [post,black!30,thin] (a47) to (a56);
\draw [post,black!30,thin] (a47) to (a59);
\draw [post,black!30,thin] (a50) to (a56);
\draw [post,black!30,thin] (a50) to (a62);
\draw [post,black!30,thin] (a53) to (a59);
\draw [post,black!30,thin] (a53) to (a62);
\draw [post,black!30,thin] (a56) to (a65);
\draw [post,black!30,thin] (a56) to (a68);
\draw [post,black!30,thin] (a56) to (a71);
\draw [post,black!30,thin] (a59) to (a65);
\draw [post,black!30,thin] (a59) to (a68);
\draw [post,black!30,thin] (a59) to (a71);
\draw [post,black!30,thin] (a62) to (a65);
\draw [post,black!30,thin] (a62) to (a68);
\draw [post,black!30,thin] (a62) to (a71);
\draw [post,black,dotted] (a2) to (a11);
\draw [post,black,dotted] (a65) to (a74);
\draw [post,black,dotted] (a5) to (a14);
\draw [post,black,dotted] (a68) to (a77);
\draw [post,black,dotted] (a8) to (a17);
\draw [post,black,dotted] (a71) to (a80);
\node(a1) at (0.555570,0.831470) {};
\node(a4) at (0.555570,1.831470) {};
\node(a7) at (0.555570,2.831470) {};
\node[vertex](a10) at (2.855570,0.831470) {};
\node[vertex](a13) at (2.855570,1.831470) {};
\node[vertex](a16) at (2.855570,2.831470) {};
\node[vertex](a19) at (5.155570,0.831470) {};
\node[vertex](a22) at (5.155570,1.831470) {};
\node[vertex](a25) at (5.155570,2.831470) {};
\node[vertex](a28) at (7.455570,0.831470) {};
\node[vertex](a31) at (7.455570,1.831470) {};
\node[vertex](a34) at (7.455570,2.831470) {};
\node[vertex](a37) at (9.755570,0.831470) {};
\node[vertex](a40) at (9.755570,1.831470) {};
\node[vertex](a43) at (9.755570,2.831470) {};
\node[vertex](a46) at (12.055570,0.831470) {};
\node[vertex](a49) at (12.055570,1.831470) {};
\node[vertex](a52) at (12.055570,2.831470) {};
\node[vertex](a55) at (14.355570,0.831470) {};
\node[vertex](a58) at (14.355570,1.831470) {};
\node[vertex](a61) at (14.355570,2.831470) {};
\node[vertex](a64) at (16.655570,0.831470) {};
\node[vertex](a67) at (16.655570,1.831470) {};
\node[vertex](a70) at (16.655570,2.831470) {};
\node(a73) at (18.955570,0.831470) {};
\node(a76) at (18.955570,1.831470) {};
\node(a79) at (18.955570,2.831470) {};
\draw [post,black!30,thin] (a10) to (a19);
\draw [post,black!30,thin] (a10) to (a20);
\draw [post,black!30,thin] (a10) to (a22);
\draw [post,black!30,thin] (a10) to (a23);
\draw [post,black!30,thin] (a11) to (a19);
\draw [post,black!30,thin] (a11) to (a22);
\draw [post,black!30,thin] (a13) to (a19);
\draw [post,black!30,thin] (a13) to (a20);
\draw [post,black!30,thin] (a13) to (a25);
\draw [post,black!30,thin] (a13) to (a26);
\draw [post,black!30,thin] (a14) to (a19);
\draw [post,black!30,thin] (a14) to (a25);
\draw [post,black!30,thin] (a16) to (a22);
\draw [post,black!30,thin] (a16) to (a23);
\draw [post,black!30,thin] (a16) to (a25);
\draw [post,black!30,thin] (a16) to (a26);
\draw [post,black!30,thin] (a17) to (a22);
\draw [post,black!30,thin] (a17) to (a25);
\draw [post,black!30,thin] (a19) to (a29);
\draw [post,black!30,thin] (a19) to (a32);
\draw [post,black!30,thin] (a19) to (a35);
\draw [post,black!30,thin] (a20) to (a28);
\draw [post,black!30,thin] (a20) to (a31);
\draw [post,black!30,thin] (a20) to (a34);
\draw [post,black!30,thin] (a22) to (a29);
\draw [post,black!30,thin] (a22) to (a32);
\draw [post,black!30,thin] (a22) to (a35);
\draw [post,black!30,thin] (a23) to (a28);
\draw [post,black!30,thin] (a23) to (a31);
\draw [post,black!30,thin] (a23) to (a34);
\draw [post,black!30,thin] (a25) to (a29);
\draw [post,black!30,thin] (a25) to (a32);
\draw [post,black!30,thin] (a25) to (a35);
\draw [post,black!30,thin] (a26) to (a28);
\draw [post,black!30,thin] (a26) to (a31);
\draw [post,black!30,thin] (a26) to (a34);
\draw [post,black!30,thin] (a28) to (a37);
\draw [post,black!30,thin] (a28) to (a38);
\draw [post,black!30,thin] (a28) to (a40);
\draw [post,black!30,thin] (a28) to (a41);
\draw [post,black!30,thin] (a29) to (a37);
\draw [post,black!30,thin] (a29) to (a40);
\draw [post,black!30,thin] (a31) to (a37);
\draw [post,black!30,thin] (a31) to (a38);
\draw [post,black!30,thin] (a31) to (a43);
\draw [post,black!30,thin] (a31) to (a44);
\draw [post,black!30,thin] (a32) to (a37);
\draw [post,black!30,thin] (a32) to (a43);
\draw [post,black!30,thin] (a34) to (a40);
\draw [post,black!30,thin] (a34) to (a41);
\draw [post,black!30,thin] (a34) to (a43);
\draw [post,black!30,thin] (a34) to (a44);
\draw [post,black!30,thin] (a35) to (a40);
\draw [post,black!30,thin] (a35) to (a43);
\draw [post,black!30,thin] (a37) to (a47);
\draw [post,black!30,thin] (a37) to (a50);
\draw [post,black!30,thin] (a37) to (a53);
\draw [post,black!30,thin] (a38) to (a46);
\draw [post,black!30,thin] (a38) to (a49);
\draw [post,black!30,thin] (a38) to (a52);
\draw [post,black!30,thin] (a40) to (a47);
\draw [post,black!30,thin] (a40) to (a50);
\draw [post,black!30,thin] (a40) to (a53);
\draw [post,black!30,thin] (a41) to (a46);
\draw [post,black!30,thin] (a41) to (a49);
\draw [post,black!30,thin] (a41) to (a52);
\draw [post,black!30,thin] (a43) to (a47);
\draw [post,black!30,thin] (a43) to (a50);
\draw [post,black!30,thin] (a43) to (a53);
\draw [post,black!30,thin] (a44) to (a46);
\draw [post,black!30,thin] (a44) to (a49);
\draw [post,black!30,thin] (a44) to (a52);
\draw [post,black!30,thin] (a46) to (a55);
\draw [post,black!30,thin] (a46) to (a56);
\draw [post,black!30,thin] (a46) to (a58);
\draw [post,black!30,thin] (a46) to (a59);
\draw [post,black!30,thin] (a47) to (a55);
\draw [post,black!30,thin] (a47) to (a58);
\draw [post,black!30,thin] (a49) to (a55);
\draw [post,black!30,thin] (a49) to (a56);
\draw [post,black!30,thin] (a49) to (a61);
\draw [post,black!30,thin] (a49) to (a62);
\draw [post,black!30,thin] (a50) to (a55);
\draw [post,black!30,thin] (a50) to (a61);
\draw [post,black!30,thin] (a52) to (a58);
\draw [post,black!30,thin] (a52) to (a59);
\draw [post,black!30,thin] (a52) to (a61);
\draw [post,black!30,thin] (a52) to (a62);
\draw [post,black!30,thin] (a53) to (a58);
\draw [post,black!30,thin] (a53) to (a61);
\draw [post,black!30,thin] (a55) to (a65);
\draw [post,black!30,thin] (a55) to (a68);
\draw [post,black!30,thin] (a55) to (a71);
\draw [post,black!30,thin] (a56) to (a64);
\draw [post,black!30,thin] (a56) to (a67);
\draw [post,black!30,thin] (a56) to (a70);
\draw [post,black!30,thin] (a58) to (a65);
\draw [post,black!30,thin] (a58) to (a68);
\draw [post,black!30,thin] (a58) to (a71);
\draw [post,black!30,thin] (a59) to (a64);
\draw [post,black!30,thin] (a59) to (a67);
\draw [post,black!30,thin] (a59) to (a70);
\draw [post,black!30,thin] (a61) to (a65);
\draw [post,black!30,thin] (a61) to (a68);
\draw [post,black!30,thin] (a61) to (a71);
\draw [post,black!30,thin] (a62) to (a64);
\draw [post,black!30,thin] (a62) to (a67);
\draw [post,black!30,thin] (a62) to (a70);
\draw [post,black,dotted] (a1) to (a10);
\draw [post,black,dotted] (a64) to (a73);
\draw [post,black,dotted] (a4) to (a13);
\draw [post,black,dotted] (a67) to (a76);
\draw [post,black,dotted] (a7) to (a16);
\draw [post,black,dotted] (a70) to (a79);
\node(a0) at (0.000000,0.000000) {};
\node(a3) at (0.000000,1.000000) {};
\node(a6) at (0.000000,2.000000) {};
\node[vertex](a9) at (2.300000,0.000000) {};
\node[vertex](a12) at (2.300000,1.000000) {};
\node[vertex](a15) at (2.300000,2.000000) {};
\node[vertex](a18) at (4.600000,0.000000) {};
\node[vertex](a21) at (4.600000,1.000000) {};
\node[vertex](a24) at (4.600000,2.000000) {};
\node[vertex](a27) at (6.900000,0.000000) {};
\node[vertex](a30) at (6.900000,1.000000) {};
\node[vertex](a33) at (6.900000,2.000000) {};
\node[vertex](a36) at (9.200000,0.000000) {};
\node[vertex](a39) at (9.200000,1.000000) {};
\node[vertex](a42) at (9.200000,2.000000) {};
\node[vertex](a45) at (11.500000,0.000000) {};
\node[vertex](a48) at (11.500000,1.000000) {};
\node[vertex](a51) at (11.500000,2.000000) {};
\node[vertex](a54) at (13.800000,0.000000) {};
\node[vertex](a57) at (13.800000,1.000000) {};
\node[vertex](a60) at (13.800000,2.000000) {};
\node[vertex](a63) at (16.100000,0.000000) {};
\node[vertex](a66) at (16.100000,1.000000) {};
\node[vertex](a69) at (16.100000,2.000000) {};
\node(a72) at (18.400000,0.000000) {};
\node(a75) at (18.400000,1.000000) {};
\node(a78) at (18.400000,2.000000) {};
\draw [post,black!30,thin] (a9) to (a18);
\draw [post,black!30,thin] (a9) to (a19);
\draw [post,black!30,thin] (a9) to (a20);
\draw [post,black!30,thin] (a9) to (a21);
\draw [post,black!30,thin] (a9) to (a22);
\draw [post,black!30,thin] (a9) to (a23);
\draw [post,black!30,thin] (a10) to (a18);
\draw [post,black!30,thin] (a10) to (a21);
\draw [post,black!30,thin] (a11) to (a18);
\draw [post,black!30,thin] (a11) to (a21);
\draw [post,black!30,thin] (a12) to (a18);
\draw [post,black!30,thin] (a12) to (a19);
\draw [post,black!30,thin] (a12) to (a20);
\draw [post,black!30,thin] (a12) to (a24);
\draw [post,black!30,thin] (a12) to (a25);
\draw [post,black!30,thin] (a12) to (a26);
\draw [post,black!30,thin] (a13) to (a18);
\draw [post,black!30,thin] (a13) to (a24);
\draw [post,black!30,thin] (a14) to (a18);
\draw [post,black!30,thin] (a14) to (a24);
\draw [post,black!30,thin] (a15) to (a21);
\draw [post,black!30,thin] (a15) to (a22);
\draw [post,black!30,thin] (a15) to (a23);
\draw [post,black!30,thin] (a15) to (a25);
\draw [post,black!30,thin] (a15) to (a26);
\draw [post,black!30,thin] (a16) to (a21);
\draw [post,black!30,thin] (a16) to (a24);
\draw [post,black!30,thin] (a17) to (a21);
\draw [post,black!30,thin] (a17) to (a24);
\draw [post,black!30,thin] (a18) to (a27);
\draw [post,black!30,thin] (a18) to (a28);
\draw [post,black!30,thin] (a18) to (a30);
\draw [post,black!30,thin] (a18) to (a31);
\draw [post,black!30,thin] (a18) to (a33);
\draw [post,black!30,thin] (a18) to (a34);
\draw [post,black!30,thin] (a19) to (a27);
\draw [post,black!30,thin] (a19) to (a30);
\draw [post,black!30,thin] (a19) to (a33);
\draw [post,black!30,thin] (a21) to (a27);
\draw [post,black!30,thin] (a21) to (a28);
\draw [post,black!30,thin] (a21) to (a30);
\draw [post,black!30,thin] (a21) to (a31);
\draw [post,black!30,thin] (a21) to (a33);
\draw [post,black!30,thin] (a21) to (a34);
\draw [post,black!30,thin] (a22) to (a27);
\draw [post,black!30,thin] (a22) to (a30);
\draw [post,black!30,thin] (a22) to (a33);
\draw [post,black!30,thin] (a24) to (a27);
\draw [post,black!30,thin] (a24) to (a28);
\draw [post,black!30,thin] (a24) to (a30);
\draw [post,black!30,thin] (a24) to (a31);
\draw [post,black!30,thin] (a24) to (a34);
\draw [post,black!30,thin] (a25) to (a27);
\draw [post,black!30,thin] (a25) to (a30);
\draw [post,black!30,thin] (a25) to (a33);
\draw [post,black!30,thin] (a27) to (a37);
\draw [post,black!30,thin] (a27) to (a38);
\draw [post,black!30,thin] (a27) to (a39);
\draw [post,black!30,thin] (a27) to (a40);
\draw [post,black!30,thin] (a27) to (a41);
\draw [post,black!30,thin] (a28) to (a36);
\draw [post,black!30,thin] (a28) to (a39);
\draw [post,black!30,thin] (a29) to (a36);
\draw [post,black!30,thin] (a29) to (a39);
\draw [post,black!30,thin] (a30) to (a36);
\draw [post,black!30,thin] (a30) to (a37);
\draw [post,black!30,thin] (a30) to (a38);
\draw [post,black!30,thin] (a30) to (a42);
\draw [post,black!30,thin] (a30) to (a43);
\draw [post,black!30,thin] (a30) to (a44);
\draw [post,black!30,thin] (a31) to (a36);
\draw [post,black!30,thin] (a31) to (a42);
\draw [post,black!30,thin] (a32) to (a36);
\draw [post,black!30,thin] (a32) to (a42);
\draw [post,black!30,thin] (a33) to (a39);
\draw [post,black!30,thin] (a33) to (a40);
\draw [post,black!30,thin] (a33) to (a41);
\draw [post,black!30,thin] (a33) to (a43);
\draw [post,black!30,thin] (a33) to (a44);
\draw [post,black!30,thin] (a34) to (a39);
\draw [post,black!30,thin] (a34) to (a42);
\draw [post,black!30,thin] (a35) to (a39);
\draw [post,black!30,thin] (a35) to (a42);
\draw [post,black!30,thin] (a36) to (a45);
\draw [post,black!30,thin] (a36) to (a46);
\draw [post,black!30,thin] (a36) to (a49);
\draw [post,black!30,thin] (a36) to (a51);
\draw [post,black!30,thin] (a36) to (a52);
\draw [post,black!30,thin] (a37) to (a45);
\draw [post,black!30,thin] (a37) to (a48);
\draw [post,black!30,thin] (a37) to (a51);
\draw [post,black!30,thin] (a39) to (a45);
\draw [post,black!30,thin] (a39) to (a46);
\draw [post,black!30,thin] (a39) to (a48);
\draw [post,black!30,thin] (a39) to (a49);
\draw [post,black!30,thin] (a39) to (a51);
\draw [post,black!30,thin] (a39) to (a52);
\draw [post,black!30,thin] (a40) to (a45);
\draw [post,black!30,thin] (a40) to (a48);
\draw [post,black!30,thin] (a40) to (a51);
\draw [post,black!30,thin] (a42) to (a45);
\draw [post,black!30,thin] (a42) to (a46);
\draw [post,black!30,thin] (a42) to (a48);
\draw [post,black!30,thin] (a42) to (a49);
\draw [post,black!30,thin] (a42) to (a52);
\draw [post,black!30,thin] (a43) to (a45);
\draw [post,black!30,thin] (a43) to (a48);
\draw [post,black!30,thin] (a43) to (a51);
\draw [post,black!30,thin] (a45) to (a54);
\draw [post,black!30,thin] (a45) to (a55);
\draw [post,black!30,thin] (a45) to (a56);
\draw [post,black!30,thin] (a45) to (a57);
\draw [post,black!30,thin] (a45) to (a58);
\draw [post,black!30,thin] (a45) to (a59);
\draw [post,black!30,thin] (a46) to (a54);
\draw [post,black!30,thin] (a46) to (a57);
\draw [post,black!30,thin] (a47) to (a54);
\draw [post,black!30,thin] (a47) to (a57);
\draw [post,black!30,thin] (a48) to (a54);
\draw [post,black!30,thin] (a48) to (a56);
\draw [post,black!30,thin] (a48) to (a60);
\draw [post,black!30,thin] (a48) to (a61);
\draw [post,black!30,thin] (a48) to (a62);
\draw [post,black!30,thin] (a49) to (a54);
\draw [post,black!30,thin] (a49) to (a60);
\draw [post,black!30,thin] (a50) to (a54);
\draw [post,black!30,thin] (a50) to (a60);
\draw [post,black!30,thin] (a51) to (a57);
\draw [post,black!30,thin] (a51) to (a58);
\draw [post,black!30,thin] (a51) to (a59);
\draw [post,black!30,thin] (a51) to (a60);
\draw [post,black!30,thin] (a51) to (a61);
\draw [post,black!30,thin] (a51) to (a62);
\draw [post,black!30,thin] (a52) to (a57);
\draw [post,black!30,thin] (a52) to (a60);
\draw [post,black!30,thin] (a53) to (a57);
\draw [post,black!30,thin] (a53) to (a60);
\draw [post,black!30,thin] (a54) to (a63);
\draw [post,black!30,thin] (a54) to (a64);
\draw [post,black!30,thin] (a54) to (a66);
\draw [post,black!30,thin] (a54) to (a67);
\draw [post,black!30,thin] (a54) to (a69);
\draw [post,black!30,thin] (a54) to (a70);
\draw [post,black!30,thin] (a55) to (a63);
\draw [post,black!30,thin] (a55) to (a69);
\draw [post,black!30,thin] (a57) to (a63);
\draw [post,black!30,thin] (a57) to (a64);
\draw [post,black!30,thin] (a57) to (a66);
\draw [post,black!30,thin] (a57) to (a67);
\draw [post,black!30,thin] (a57) to (a69);
\draw [post,black!30,thin] (a57) to (a70);
\draw [post,black!30,thin] (a58) to (a63);
\draw [post,black!30,thin] (a58) to (a66);
\draw [post,black!30,thin] (a58) to (a69);
\draw [post,black!30,thin] (a60) to (a63);
\draw [post,black!30,thin] (a60) to (a64);
\draw [post,black!30,thin] (a60) to (a66);
\draw [post,black!30,thin] (a60) to (a67);
\draw [post,black!30,thin] (a60) to (a69);
\draw [post,black!30,thin] (a60) to (a70);
\draw [post,black!30,thin] (a61) to (a63);
\draw [post,black!30,thin] (a61) to (a66);
\draw [post,black!30,thin] (a61) to (a69);
\draw [post,black,dotted] (a0) to (a9);
\draw [post,black,dotted] (a63) to (a72);
\draw [post,black,dotted] (a3) to (a12);
\draw [post,black,dotted] (a66) to (a75);
\draw [post,black,dotted] (a6) to (a15);
\draw [post,black,dotted] (a69) to (a78);
\draw [post,black,very thick] (a15) to (a24);
\draw [post,black,very thick] (a24) to (a33);
\draw [post,black,very thick] (a27) to (a36);
\draw [post,black,very thick] (a33) to (a42);
\draw [post,black,very thick] (a36) to (a48);
\draw [post,black,very thick] (a42) to (a51);
\draw [post,black,very thick] (a48) to (a55);
\draw [post,black,very thick] (a55) to (a66);
\end{tikzpicture}
}

\def\TubeTol{
\begin{tikzpicture}
[scale=0.50,
inner sep=1.5,font=\ttfamily,
movi/.style={rectangle,draw=black,minimum size=4mm},
vertex/.style={circle,draw=black!50,fill=black!20, very thick},
parta/.style={circle,draw=black!50,fill=black!0, very thick},
partb/.style={circle,draw=black!50,fill=black!100, very thick},
edgevertex/.style={circle,draw=blue!80,fill=blue!40, very thick},
curly/.style={decorate,decoration={snake,amplitude=0.2mm,segment length=1mm,post length=1mm}},
normal/.style=,
post/.style={->, >=stealth'},
pre/.style={<-, >=stealth'}]
\node(a0) at (1.000000,0.000000) {\footnotesize };
\node(a1) at (1.000000,1.000000) {\footnotesize };
\node(a2) at (1.000000,2.000000) {\footnotesize };
\node(a3) at (2.300000,0.000000) {\footnotesize +};
\node(a4) at (2.300000,1.000000) {\footnotesize o};
\node(a5) at (2.300000,2.000000) {\footnotesize -};
\node(a6) at (4.600000,0.000000) {\footnotesize +};
\node(a7) at (4.600000,1.000000) {\footnotesize o};
\node(a8) at (4.600000,2.000000) {\footnotesize -};
\node(a9) at (6.900000,0.000000) {\footnotesize +};
\node(a10) at (6.900000,1.000000) [movi] {\large\textbf{o}};
\node(a11) at (6.900000,2.000000) {\footnotesize -};
\node(a12) at (9.200000,0.000000) {\footnotesize +};
\node(a13) at (9.200000,1.000000) {\footnotesize o};
\node(a14) at (9.200000,2.000000) {\footnotesize -};
\node(a15) at (11.500000,0.000000) {\footnotesize +};
\node(a16) at (11.500000,1.000000) {\footnotesize o};
\node(a17) at (11.500000,2.000000) {\footnotesize -};
\node(a18) at (12.800000,0.000000) {\footnotesize };
\node(a19) at (12.800000,1.000000) {\footnotesize };
\node(a20) at (12.800000,2.000000) {\footnotesize };
\node(a21) at (2.300000,-1.000000) {\footnotesize -2};
\node(a22) at (4.600000,-1.000000) {\footnotesize -1};
\node(a23) at (6.900000,-1.000000) {\footnotesize  0};
\node(a24) at (9.200000,-1.000000) {\footnotesize  1};
\node(a25) at (11.500000,-1.000000) {\footnotesize  2};
\draw [post,black,dashed] (a3) to (a6);
\draw [post,black,dashed] (a3) to (a7);
\draw [post,black!30,thin] (a3) to (a8);
\draw [post,black,thick] (a4) to (a6);
\draw [post,black,thick] (a4) to (a7);
\draw [post,black!30,thin] (a4) to (a8);
\draw [post,black,curly] (a5) to (a6);
\draw [post,black,curly] (a5) to (a7);
\draw [post,black!30,thin] (a5) to (a8);
\draw [post,black,curly] (a6) to (a9);
\draw [post,black,dashed] (a6) to (a10);
\draw [post,black,curly] (a7) to (a9);
\draw [post,black,dashed] (a7) to (a11);
\draw [post,black,thick] (a8) to (a10);
\draw [post,black,thick] (a8) to (a11);
\draw [post,black!30,thin] (a9) to (a12);
\draw [post,black!30,thin] (a9) to (a13);
\draw [post,black!30,thin] (a9) to (a14);
\draw [post,black,thick] (a10) to (a12);
\draw [post,black,curly] (a10) to (a13);
\draw [post,black,dashed] (a10) to (a14);
\draw [post,black,thick] (a11) to (a12);
\draw [post,black,curly] (a11) to (a13);
\draw [post,black,dashed] (a11) to (a14);
\draw [post,black!30,thin] (a12) to (a15);
\draw [post,black!30,thin] (a12) to (a16);
\draw [post,black!30,thin] (a13) to (a15);
\draw [post,black!30,thin] (a13) to (a17);
\draw [post,black!30,thin] (a14) to (a16);
\draw [post,black!30,thin] (a14) to (a17);
\draw [post,black!30,dotted] (a0) to (a3);
\draw [post,black!30,dotted] (a15) to (a18);
\draw [post,black!30,dotted] (a1) to (a4);
\draw [post,black!30,dotted] (a16) to (a19);
\draw [post,black!30,dotted] (a2) to (a5);
\draw [post,black!30,dotted] (a17) to (a20);
\end{tikzpicture}
}

\def\TubeTor{
\begin{tikzpicture}
[scale=0.50,
inner sep=1.5,font=\ttfamily,
movi/.style={rectangle,draw=black,minimum size=4mm},
vertex/.style={circle,draw=black!50,fill=black!20, very thick},
parta/.style={circle,draw=black!50,fill=black!0, very thick},
partb/.style={circle,draw=black!50,fill=black!100, very thick},
edgevertex/.style={circle,draw=blue!80,fill=blue!40, very thick},
curly/.style={decorate,decoration={snake,amplitude=0.2mm,segment length=1mm,post length=1mm}},
normal/.style=,
post/.style={->, >=stealth'},
pre/.style={<-, >=stealth'}]
\node(a0) at (1.000000,0.000000) {\footnotesize };
\node(a1) at (1.000000,1.000000) {\footnotesize };
\node(a2) at (1.000000,2.000000) {\footnotesize };
\node(a3) at (2.300000,0.000000) {\footnotesize +};
\node(a4) at (2.300000,1.000000) {\footnotesize o};
\node(a5) at (2.300000,2.000000) {\footnotesize -};
\node(a6) at (4.600000,1.000000) [movi] {\large\textbf{+}};
\node(a7) at (4.600000,0.000000) [movi] {\large\textbf{o}};
\node(a8) at (4.600000,2.000000) {\footnotesize -};
\node(a9) at (6.900000,0.000000) {\footnotesize +};
\node(a10) at (6.900000,2.000000) [movi] {\large\textbf{o}};
\node(a11) at (6.900000,1.000000) [movi] {\large\textbf{-}};
\node(a12) at (9.200000,0.000000) {\footnotesize +};
\node(a13) at (9.200000,1.000000) {\footnotesize o};
\node(a14) at (9.200000,2.000000) {\footnotesize -};
\node(a15) at (11.500000,0.000000) {\footnotesize +};
\node(a16) at (11.500000,1.000000) {\footnotesize o};
\node(a17) at (11.500000,2.000000) {\footnotesize -};
\node(a18) at (12.800000,0.000000) {\footnotesize };
\node(a19) at (12.800000,1.000000) {\footnotesize };
\node(a20) at (12.800000,2.000000) {\footnotesize };
\node(a21) at (2.300000,-1.000000) {\footnotesize -2};
\node(a22) at (4.600000,-1.000000) {\footnotesize -1};
\node(a23) at (6.900000,-1.000000) {\footnotesize  0};
\node(a24) at (9.200000,-1.000000) {\footnotesize  1};
\node(a25) at (11.500000,-1.000000) {\footnotesize  2};
\draw [post,black,dashed] (a3) to (a6);
\draw [post,black,dashed] (a3) to (a7);
\draw [post,black!30,thin] (a3) to (a8);
\draw [post,black,thick] (a4) to (a6);
\draw [post,black,thick] (a4) to (a7);
\draw [post,black!30,thin] (a4) to (a8);
\draw [post,black,curly] (a5) to (a6);
\draw [post,black,curly] (a5) to (a7);
\draw [post,black!30,thin] (a5) to (a8);
\draw [post,black,curly] (a6) to (a9);
\draw [post,black,dashed] (a6) to (a10);
\draw [post,black,curly] (a7) to (a9);
\draw [post,black,dashed] (a7) to (a11);
\draw [post,black,thick] (a8) to (a10);
\draw [post,black,thick] (a8) to (a11);
\draw [post,black!30,thin] (a9) to (a12);
\draw [post,black!30,thin] (a9) to (a13);
\draw [post,black!30,thin] (a9) to (a14);
\draw [post,black,thick] (a10) to (a12);
\draw [post,black,curly] (a10) to (a13);
\draw [post,black,dashed] (a10) to (a14);
\draw [post,black,thick] (a11) to (a12);
\draw [post,black,curly] (a11) to (a13);
\draw [post,black,dashed] (a11) to (a14);
\draw [post,black!30,thin] (a12) to (a15);
\draw [post,black!30,thin] (a12) to (a16);
\draw [post,black!30,thin] (a13) to (a15);
\draw [post,black!30,thin] (a13) to (a17);
\draw [post,black!30,thin] (a14) to (a16);
\draw [post,black!30,thin] (a14) to (a17);
\draw [post,black!30,dotted] (a0) to (a3);
\draw [post,black!30,dotted] (a15) to (a18);
\draw [post,black!30,dotted] (a1) to (a4);
\draw [post,black!30,dotted] (a16) to (a19);
\draw [post,black!30,dotted] (a2) to (a5);
\draw [post,black!30,dotted] (a17) to (a20);
\end{tikzpicture}
}

\def\TubeTpl{
\begin{tikzpicture}
[scale=0.50,
inner sep=1.5,font=\ttfamily,
movi/.style={rectangle,draw=black,minimum size=4mm},
vertex/.style={circle,draw=black!50,fill=black!20, very thick},
parta/.style={circle,draw=black!50,fill=black!0, very thick},
partb/.style={circle,draw=black!50,fill=black!100, very thick},
edgevertex/.style={circle,draw=blue!80,fill=blue!40, very thick},
curly/.style={decorate,decoration={snake,amplitude=0.2mm,segment length=1mm,post length=1mm}},
normal/.style=,
post/.style={->, >=stealth'},
pre/.style={<-, >=stealth'}]
\node(a0) at (1.000000,0.000000) {\footnotesize };
\node(a1) at (1.000000,1.000000) {\footnotesize };
\node(a2) at (1.000000,2.000000) {\footnotesize };
\node(a3) at (2.300000,0.000000) {\footnotesize +};
\node(a4) at (2.300000,1.000000) {\footnotesize o};
\node(a5) at (2.300000,2.000000) {\footnotesize -};
\node(a6) at (4.600000,0.000000) {\footnotesize +};
\node(a7) at (4.600000,1.000000) {\footnotesize o};
\node(a8) at (4.600000,2.000000) {\footnotesize -};
\node(a9) at (6.900000,0.000000) [movi] {\large\textbf{+}};
\node(a10) at (6.900000,1.000000) {\footnotesize o};
\node(a11) at (6.900000,2.000000) {\footnotesize -};
\node(a12) at (9.200000,0.000000) {\footnotesize +};
\node(a13) at (9.200000,1.000000) {\footnotesize o};
\node(a14) at (9.200000,2.000000) {\footnotesize -};
\node(a15) at (11.500000,0.000000) {\footnotesize +};
\node(a16) at (11.500000,1.000000) {\footnotesize o};
\node(a17) at (11.500000,2.000000) {\footnotesize -};
\node(a18) at (12.800000,0.000000) {\footnotesize };
\node(a19) at (12.800000,1.000000) {\footnotesize };
\node(a20) at (12.800000,2.000000) {\footnotesize };
\node(a21) at (2.300000,-1.000000) {\footnotesize -2};
\node(a22) at (4.600000,-1.000000) {\footnotesize -1};
\node(a23) at (6.900000,-1.000000) {\footnotesize  0};
\node(a24) at (9.200000,-1.000000) {\footnotesize  1};
\node(a25) at (11.500000,-1.000000) {\footnotesize  2};
\draw [post,black,curly] (a3) to (a6);
\draw [post,black!30,thin] (a3) to (a7);
\draw [post,black,curly] (a3) to (a8);
\draw [post,black,dashed] (a4) to (a6);
\draw [post,black!30,thin] (a4) to (a7);
\draw [post,black,dashed] (a4) to (a8);
\draw [post,black,thick] (a5) to (a6);
\draw [post,black!30,thin] (a5) to (a7);
\draw [post,black,thick] (a5) to (a8);
\draw [post,black,thick] (a6) to (a9);
\draw [post,black,dashed] (a6) to (a10);
\draw [post,black,curly] (a7) to (a9);
\draw [post,black,curly] (a7) to (a11);
\draw [post,black,dashed] (a8) to (a10);
\draw [post,black,thick] (a8) to (a11);
\draw [post,black,thick] (a9) to (a12);
\draw [post,black,curly] (a9) to (a13);
\draw [post,black,dashed] (a9) to (a14);
\draw [post,black!30,thin] (a10) to (a12);
\draw [post,black!30,thin] (a10) to (a13);
\draw [post,black!30,thin] (a10) to (a14);
\draw [post,black,thick] (a11) to (a12);
\draw [post,black,curly] (a11) to (a13);
\draw [post,black,dashed] (a11) to (a14);
\draw [post,black!30,thin] (a12) to (a15);
\draw [post,black!30,thin] (a12) to (a16);
\draw [post,black!30,thin] (a13) to (a15);
\draw [post,black!30,thin] (a13) to (a17);
\draw [post,black!30,thin] (a14) to (a16);
\draw [post,black!30,thin] (a14) to (a17);
\draw [post,black!30,dotted] (a0) to (a3);
\draw [post,black!30,dotted] (a15) to (a18);
\draw [post,black!30,dotted] (a1) to (a4);
\draw [post,black!30,dotted] (a16) to (a19);
\draw [post,black!30,dotted] (a2) to (a5);
\draw [post,black!30,dotted] (a17) to (a20);
\end{tikzpicture}
}

\def\TubeTpr{
\begin{tikzpicture}
[scale=0.50,
inner sep=1.5,font=\ttfamily,
movi/.style={rectangle,draw=black,minimum size=4mm},
vertex/.style={circle,draw=black!50,fill=black!20, very thick},
parta/.style={circle,draw=black!50,fill=black!0, very thick},
partb/.style={circle,draw=black!50,fill=black!100, very thick},
edgevertex/.style={circle,draw=blue!80,fill=blue!40, very thick},
curly/.style={decorate,decoration={snake,amplitude=0.2mm,segment length=1mm,post length=1mm}},
normal/.style=,
post/.style={->, >=stealth'},
pre/.style={<-, >=stealth'}]
\node(a0) at (1.000000,0.000000) {\footnotesize };
\node(a1) at (1.000000,1.000000) {\footnotesize };
\node(a2) at (1.000000,2.000000) {\footnotesize };
\node(a3) at (2.300000,0.000000) {\footnotesize +};
\node(a4) at (2.300000,1.000000) {\footnotesize o};
\node(a5) at (2.300000,2.000000) {\footnotesize -};
\node(a6) at (4.600000,2.000000) [movi] {\large\textbf{+}};
\node(a7) at (4.600000,1.000000) {\footnotesize o};
\node(a8) at (4.600000,0.000000) [movi] {\large\textbf{-}};
\node(a9) at (6.900000,2.000000) [movi] {\large\textbf{+}};
\node(a10) at (6.900000,1.000000) {\footnotesize o};
\node(a11) at (6.900000,0.000000) [movi] {\large\textbf{-}};
\node(a12) at (9.200000,0.000000) {\footnotesize +};
\node(a13) at (9.200000,1.000000) {\footnotesize o};
\node(a14) at (9.200000,2.000000) {\footnotesize -};
\node(a15) at (11.500000,0.000000) {\footnotesize +};
\node(a16) at (11.500000,1.000000) {\footnotesize o};
\node(a17) at (11.500000,2.000000) {\footnotesize -};
\node(a18) at (12.800000,0.000000) {\footnotesize };
\node(a19) at (12.800000,1.000000) {\footnotesize };
\node(a20) at (12.800000,2.000000) {\footnotesize };
\node(a21) at (2.300000,-1.000000) {\footnotesize -2};
\node(a22) at (4.600000,-1.000000) {\footnotesize -1};
\node(a23) at (6.900000,-1.000000) {\footnotesize  0};
\node(a24) at (9.200000,-1.000000) {\footnotesize  1};
\node(a25) at (11.500000,-1.000000) {\footnotesize  2};
\draw [post,black,curly] (a3) to (a6);
\draw [post,black!30,thin] (a3) to (a7);
\draw [post,black,curly] (a3) to (a8);
\draw [post,black,dashed] (a4) to (a6);
\draw [post,black!30,thin] (a4) to (a7);
\draw [post,black,dashed] (a4) to (a8);
\draw [post,black,thick] (a5) to (a6);
\draw [post,black!30,thin] (a5) to (a7);
\draw [post,black,thick] (a5) to (a8);
\draw [post,black,thick] (a6) to (a9);
\draw [post,black,dashed] (a6) to (a10);
\draw [post,black,curly] (a7) to (a9);
\draw [post,black,curly] (a7) to (a11);
\draw [post,black,dashed] (a8) to (a10);
\draw [post,black,thick] (a8) to (a11);
\draw [post,black,thick] (a9) to (a12);
\draw [post,black,curly] (a9) to (a13);
\draw [post,black,dashed] (a9) to (a14);
\draw [post,black!30,thin] (a10) to (a12);
\draw [post,black!30,thin] (a10) to (a13);
\draw [post,black!30,thin] (a10) to (a14);
\draw [post,black,thick] (a11) to (a12);
\draw [post,black,curly] (a11) to (a13);
\draw [post,black,dashed] (a11) to (a14);
\draw [post,black!30,thin] (a12) to (a15);
\draw [post,black!30,thin] (a12) to (a16);
\draw [post,black!30,thin] (a13) to (a15);
\draw [post,black!30,thin] (a13) to (a17);
\draw [post,black!30,thin] (a14) to (a16);
\draw [post,black!30,thin] (a14) to (a17);
\draw [post,black!30,dotted] (a0) to (a3);
\draw [post,black!30,dotted] (a15) to (a18);
\draw [post,black!30,dotted] (a1) to (a4);
\draw [post,black!30,dotted] (a16) to (a19);
\draw [post,black!30,dotted] (a2) to (a5);
\draw [post,black!30,dotted] (a17) to (a20);
\end{tikzpicture}
}

\def\TubeTOl{
\begin{tikzpicture}
[scale=0.50,
inner sep=1.5,font=\ttfamily,
movi/.style={rectangle,draw=black,minimum size=4mm},
vertex/.style={circle,draw=black!50,fill=black!20, very thick},
parta/.style={circle,draw=black!50,fill=black!0, very thick},
partb/.style={circle,draw=black!50,fill=black!100, very thick},
edgevertex/.style={circle,draw=blue!80,fill=blue!40, very thick},
curly/.style={decorate,decoration={snake,amplitude=0.2mm,segment length=1mm,post length=1mm}},
normal/.style=,
post/.style={->, >=stealth'},
pre/.style={<-, >=stealth'}]
\node(a0) at (1.000000,0.000000) {\footnotesize };
\node(a1) at (1.000000,1.000000) {\footnotesize };
\node(a2) at (1.000000,2.000000) {\footnotesize };
\node(a3) at (2.300000,0.000000) {\footnotesize +};
\node(a4) at (2.300000,1.000000) {\footnotesize o};
\node(a5) at (2.300000,2.000000) {\footnotesize -};
\node(a6) at (4.600000,0.000000) {\footnotesize +};
\node(a7) at (4.600000,1.000000) {\footnotesize o};
\node(a8) at (4.600000,2.000000) {\footnotesize -};
\node(a9) at (6.900000,0.000000) {\footnotesize +};
\node(a10) at (6.900000,1.000000) [movi] {\large\textbf{o}};
\node(a11) at (6.900000,2.000000) {\footnotesize -};
\node(a12) at (9.200000,0.000000) {\footnotesize +};
\node(a13) at (9.200000,1.000000) {\footnotesize o};
\node(a14) at (9.200000,2.000000) {\footnotesize -};
\node(a15) at (11.500000,0.000000) {\footnotesize +};
\node(a16) at (11.500000,1.000000) {\footnotesize o};
\node(a17) at (11.500000,2.000000) {\footnotesize -};
\node(a18) at (12.800000,0.000000) {\footnotesize };
\node(a19) at (12.800000,1.000000) {\footnotesize };
\node(a20) at (12.800000,2.000000) {\footnotesize };
\node(a21) at (2.300000,-1.000000) {\footnotesize -2};
\node(a22) at (4.600000,-1.000000) {\footnotesize -1};
\node(a23) at (6.900000,-1.000000) {\footnotesize  0};
\node(a24) at (9.200000,-1.000000) {\footnotesize  1};
\node(a25) at (11.500000,-1.000000) {\footnotesize  2};
\draw [post,black!30,thin] (a3) to (a6);
\draw [post,black!30,thin] (a3) to (a7);
\draw [post,black!30,thin] (a4) to (a6);
\draw [post,black!30,thin] (a4) to (a8);
\draw [post,black!30,thin] (a5) to (a7);
\draw [post,black!30,thin] (a5) to (a8);
\draw [post,black!30,thin] (a6) to (a9);
\draw [post,black,curly] (a6) to (a10);
\draw [post,black,curly] (a6) to (a11);
\draw [post,black!30,thin] (a7) to (a9);
\draw [post,black,dashed] (a7) to (a10);
\draw [post,black,dashed] (a7) to (a11);
\draw [post,black!30,thin] (a8) to (a9);
\draw [post,black,thick] (a8) to (a10);
\draw [post,black,thick] (a8) to (a11);
\draw [post,black,thick] (a9) to (a12);
\draw [post,black,thick] (a9) to (a13);
\draw [post,black,curly] (a10) to (a12);
\draw [post,black,dashed] (a10) to (a14);
\draw [post,black,curly] (a11) to (a13);
\draw [post,black,dashed] (a11) to (a14);
\draw [post,black,thick] (a12) to (a15);
\draw [post,black,dashed] (a12) to (a16);
\draw [post,black,curly] (a12) to (a17);
\draw [post,black,thick] (a13) to (a15);
\draw [post,black,dashed] (a13) to (a16);
\draw [post,black,curly] (a13) to (a17);
\draw [post,black!30,thin] (a14) to (a15);
\draw [post,black!30,thin] (a14) to (a16);
\draw [post,black!30,thin] (a14) to (a17);
\draw [post,black!30,dotted] (a0) to (a3);
\draw [post,black!30,dotted] (a15) to (a18);
\draw [post,black!30,dotted] (a1) to (a4);
\draw [post,black!30,dotted] (a16) to (a19);
\draw [post,black!30,dotted] (a2) to (a5);
\draw [post,black!30,dotted] (a17) to (a20);
\end{tikzpicture}
}

\def\TubeTOr{
\begin{tikzpicture}
[scale=0.50,
inner sep=1.5,font=\ttfamily,
movi/.style={rectangle,draw=black,minimum size=4mm},
vertex/.style={circle,draw=black!50,fill=black!20, very thick},
parta/.style={circle,draw=black!50,fill=black!0, very thick},
partb/.style={circle,draw=black!50,fill=black!100, very thick},
edgevertex/.style={circle,draw=blue!80,fill=blue!40, very thick},
curly/.style={decorate,decoration={snake,amplitude=0.2mm,segment length=1mm,post length=1mm}},
normal/.style=,
post/.style={->, >=stealth'},
pre/.style={<-, >=stealth'}]
\node(a0) at (1.000000,0.000000) {\footnotesize };
\node(a1) at (1.000000,1.000000) {\footnotesize };
\node(a2) at (1.000000,2.000000) {\footnotesize };
\node(a3) at (2.300000,0.000000) {\footnotesize +};
\node(a4) at (2.300000,1.000000) {\footnotesize o};
\node(a5) at (2.300000,2.000000) {\footnotesize -};
\node(a6) at (4.600000,0.000000) {\footnotesize +};
\node(a7) at (4.600000,1.000000) {\footnotesize o};
\node(a8) at (4.600000,2.000000) {\footnotesize -};
\node(a9) at (6.900000,0.000000) {\footnotesize +};
\node(a10) at (6.900000,2.000000) [movi] {\large\textbf{o}};
\node(a11) at (6.900000,1.000000) [movi] {\large\textbf{-}};
\node(a12) at (9.200000,1.000000) [movi] {\large\textbf{+}};
\node(a13) at (9.200000,0.000000) [movi] {\large\textbf{o}};
\node(a14) at (9.200000,2.000000) {\footnotesize -};
\node(a15) at (11.500000,0.000000) {\footnotesize +};
\node(a16) at (11.500000,1.000000) {\footnotesize o};
\node(a17) at (11.500000,2.000000) {\footnotesize -};
\node(a18) at (12.800000,0.000000) {\footnotesize };
\node(a19) at (12.800000,1.000000) {\footnotesize };
\node(a20) at (12.800000,2.000000) {\footnotesize };
\node(a21) at (2.300000,-1.000000) {\footnotesize -2};
\node(a22) at (4.600000,-1.000000) {\footnotesize -1};
\node(a23) at (6.900000,-1.000000) {\footnotesize  0};
\node(a24) at (9.200000,-1.000000) {\footnotesize  1};
\node(a25) at (11.500000,-1.000000) {\footnotesize  2};
\draw [post,black!30,thin] (a3) to (a6);
\draw [post,black!30,thin] (a3) to (a7);
\draw [post,black!30,thin] (a4) to (a6);
\draw [post,black!30,thin] (a4) to (a8);
\draw [post,black!30,thin] (a5) to (a7);
\draw [post,black!30,thin] (a5) to (a8);
\draw [post,black!30,thin] (a6) to (a9);
\draw [post,black,curly] (a6) to (a10);
\draw [post,black,curly] (a6) to (a11);
\draw [post,black!30,thin] (a7) to (a9);
\draw [post,black,dashed] (a7) to (a10);
\draw [post,black,dashed] (a7) to (a11);
\draw [post,black!30,thin] (a8) to (a9);
\draw [post,black,thick] (a8) to (a10);
\draw [post,black,thick] (a8) to (a11);
\draw [post,black,thick] (a9) to (a12);
\draw [post,black,thick] (a9) to (a13);
\draw [post,black,curly] (a10) to (a12);
\draw [post,black,dashed] (a10) to (a14);
\draw [post,black,curly] (a11) to (a13);
\draw [post,black,dashed] (a11) to (a14);
\draw [post,black,thick] (a12) to (a15);
\draw [post,black,dashed] (a12) to (a16);
\draw [post,black,curly] (a12) to (a17);
\draw [post,black,thick] (a13) to (a15);
\draw [post,black,dashed] (a13) to (a16);
\draw [post,black,curly] (a13) to (a17);
\draw [post,black!30,thin] (a14) to (a15);
\draw [post,black!30,thin] (a14) to (a16);
\draw [post,black!30,thin] (a14) to (a17);
\draw [post,black!30,dotted] (a0) to (a3);
\draw [post,black!30,dotted] (a15) to (a18);
\draw [post,black!30,dotted] (a1) to (a4);
\draw [post,black!30,dotted] (a16) to (a19);
\draw [post,black!30,dotted] (a2) to (a5);
\draw [post,black!30,dotted] (a17) to (a20);
\end{tikzpicture}
}

\def\TubeTPl{
\begin{tikzpicture}
[scale=0.50,
inner sep=1.5,font=\ttfamily,
movi/.style={rectangle,draw=black,minimum size=4mm},
vertex/.style={circle,draw=black!50,fill=black!20, very thick},
parta/.style={circle,draw=black!50,fill=black!0, very thick},
partb/.style={circle,draw=black!50,fill=black!100, very thick},
edgevertex/.style={circle,draw=blue!80,fill=blue!40, very thick},
curly/.style={decorate,decoration={snake,amplitude=0.2mm,segment length=1mm,post length=1mm}},
normal/.style=,
post/.style={->, >=stealth'},
pre/.style={<-, >=stealth'}]
\node(a0) at (1.000000,0.000000) {\footnotesize };
\node(a1) at (1.000000,1.000000) {\footnotesize };
\node(a2) at (1.000000,2.000000) {\footnotesize };
\node(a3) at (2.300000,0.000000) {\footnotesize +};
\node(a4) at (2.300000,1.000000) {\footnotesize o};
\node(a5) at (2.300000,2.000000) {\footnotesize -};
\node(a6) at (4.600000,0.000000) {\footnotesize +};
\node(a7) at (4.600000,1.000000) {\footnotesize o};
\node(a8) at (4.600000,2.000000) {\footnotesize -};
\node(a9) at (6.900000,0.000000) [movi] {\large\textbf{+}};
\node(a10) at (6.900000,1.000000) {\footnotesize o};
\node(a11) at (6.900000,2.000000) {\footnotesize -};
\node(a12) at (9.200000,0.000000) {\footnotesize +};
\node(a13) at (9.200000,1.000000) {\footnotesize o};
\node(a14) at (9.200000,2.000000) {\footnotesize -};
\node(a15) at (11.500000,0.000000) {\footnotesize +};
\node(a16) at (11.500000,1.000000) {\footnotesize o};
\node(a17) at (11.500000,2.000000) {\footnotesize -};
\node(a18) at (12.800000,0.000000) {\footnotesize };
\node(a19) at (12.800000,1.000000) {\footnotesize };
\node(a20) at (12.800000,2.000000) {\footnotesize };
\node(a21) at (2.300000,-1.000000) {\footnotesize -2};
\node(a22) at (4.600000,-1.000000) {\footnotesize -1};
\node(a23) at (6.900000,-1.000000) {\footnotesize  0};
\node(a24) at (9.200000,-1.000000) {\footnotesize  1};
\node(a25) at (11.500000,-1.000000) {\footnotesize  2};
\draw [post,black!30,thin] (a3) to (a6);
\draw [post,black!30,thin] (a3) to (a7);
\draw [post,black!30,thin] (a4) to (a6);
\draw [post,black!30,thin] (a4) to (a8);
\draw [post,black!30,thin] (a5) to (a7);
\draw [post,black!30,thin] (a5) to (a8);
\draw [post,black,dashed] (a6) to (a9);
\draw [post,black!30,thin] (a6) to (a10);
\draw [post,black,dashed] (a6) to (a11);
\draw [post,black,curly] (a7) to (a9);
\draw [post,black!30,thin] (a7) to (a10);
\draw [post,black,curly] (a7) to (a11);
\draw [post,black,thick] (a8) to (a9);
\draw [post,black!30,thin] (a8) to (a10);
\draw [post,black,thick] (a8) to (a11);
\draw [post,black,thick] (a9) to (a12);
\draw [post,black,dashed] (a9) to (a13);
\draw [post,black,curly] (a10) to (a12);
\draw [post,black,curly] (a10) to (a14);
\draw [post,black,dashed] (a11) to (a13);
\draw [post,black,thick] (a11) to (a14);
\draw [post,black,thick] (a12) to (a15);
\draw [post,black,dashed] (a12) to (a16);
\draw [post,black,curly] (a12) to (a17);
\draw [post,black!30,thin] (a13) to (a15);
\draw [post,black!30,thin] (a13) to (a16);
\draw [post,black!30,thin] (a13) to (a17);
\draw [post,black,thick] (a14) to (a15);
\draw [post,black,dashed] (a14) to (a16);
\draw [post,black,curly] (a14) to (a17);
\draw [post,black!30,dotted] (a0) to (a3);
\draw [post,black!30,dotted] (a15) to (a18);
\draw [post,black!30,dotted] (a1) to (a4);
\draw [post,black!30,dotted] (a16) to (a19);
\draw [post,black!30,dotted] (a2) to (a5);
\draw [post,black!30,dotted] (a17) to (a20);
\end{tikzpicture}
}

\def\TubeTPr{
\begin{tikzpicture}
[scale=0.50,
inner sep=1.5,font=\ttfamily,
movi/.style={rectangle,draw=black,minimum size=4mm},
vertex/.style={circle,draw=black!50,fill=black!20, very thick},
parta/.style={circle,draw=black!50,fill=black!0, very thick},
partb/.style={circle,draw=black!50,fill=black!100, very thick},
edgevertex/.style={circle,draw=blue!80,fill=blue!40, very thick},
curly/.style={decorate,decoration={snake,amplitude=0.2mm,segment length=1mm,post length=1mm}},
normal/.style=,
post/.style={->, >=stealth'},
pre/.style={<-, >=stealth'}]
\node(a0) at (1.000000,0.000000) {\footnotesize };
\node(a1) at (1.000000,1.000000) {\footnotesize };
\node(a2) at (1.000000,2.000000) {\footnotesize };
\node(a3) at (2.300000,0.000000) {\footnotesize +};
\node(a4) at (2.300000,1.000000) {\footnotesize o};
\node(a5) at (2.300000,2.000000) {\footnotesize -};
\node(a6) at (4.600000,0.000000) {\footnotesize +};
\node(a7) at (4.600000,1.000000) {\footnotesize o};
\node(a8) at (4.600000,2.000000) {\footnotesize -};
\node(a9) at (6.900000,2.000000) [movi] {\large\textbf{+}};
\node(a10) at (6.900000,1.000000) {\footnotesize o};
\node(a11) at (6.900000,0.000000) [movi] {\large\textbf{-}};
\node(a12) at (9.200000,2.000000) [movi] {\large\textbf{+}};
\node(a13) at (9.200000,1.000000) {\footnotesize o};
\node(a14) at (9.200000,0.000000) [movi] {\large\textbf{-}};
\node(a15) at (11.500000,0.000000) {\footnotesize +};
\node(a16) at (11.500000,1.000000) {\footnotesize o};
\node(a17) at (11.500000,2.000000) {\footnotesize -};
\node(a18) at (12.800000,0.000000) {\footnotesize };
\node(a19) at (12.800000,1.000000) {\footnotesize };
\node(a20) at (12.800000,2.000000) {\footnotesize };
\node(a21) at (2.300000,-1.000000) {\footnotesize -2};
\node(a22) at (4.600000,-1.000000) {\footnotesize -1};
\node(a23) at (6.900000,-1.000000) {\footnotesize  0};
\node(a24) at (9.200000,-1.000000) {\footnotesize  1};
\node(a25) at (11.500000,-1.000000) {\footnotesize  2};
\draw [post,black!30,thin] (a3) to (a6);
\draw [post,black!30,thin] (a3) to (a7);
\draw [post,black!30,thin] (a4) to (a6);
\draw [post,black!30,thin] (a4) to (a8);
\draw [post,black!30,thin] (a5) to (a7);
\draw [post,black!30,thin] (a5) to (a8);
\draw [post,black,dashed] (a6) to (a9);
\draw [post,black!30,thin] (a6) to (a10);
\draw [post,black,dashed] (a6) to (a11);
\draw [post,black,curly] (a7) to (a9);
\draw [post,black!30,thin] (a7) to (a10);
\draw [post,black,curly] (a7) to (a11);
\draw [post,black,thick] (a8) to (a9);
\draw [post,black!30,thin] (a8) to (a10);
\draw [post,black,thick] (a8) to (a11);
\draw [post,black,thick] (a9) to (a12);
\draw [post,black,dashed] (a9) to (a13);
\draw [post,black,curly] (a10) to (a12);
\draw [post,black,curly] (a10) to (a14);
\draw [post,black,dashed] (a11) to (a13);
\draw [post,black,thick] (a11) to (a14);
\draw [post,black,thick] (a12) to (a15);
\draw [post,black,dashed] (a12) to (a16);
\draw [post,black,curly] (a12) to (a17);
\draw [post,black!30,thin] (a13) to (a15);
\draw [post,black!30,thin] (a13) to (a16);
\draw [post,black!30,thin] (a13) to (a17);
\draw [post,black,thick] (a14) to (a15);
\draw [post,black,dashed] (a14) to (a16);
\draw [post,black,curly] (a14) to (a17);
\draw [post,black!30,dotted] (a0) to (a3);
\draw [post,black!30,dotted] (a15) to (a18);
\draw [post,black!30,dotted] (a1) to (a4);
\draw [post,black!30,dotted] (a16) to (a19);
\draw [post,black!30,dotted] (a2) to (a5);
\draw [post,black!30,dotted] (a17) to (a20);
\end{tikzpicture}
}

\def\TubeAol{
\begin{tikzpicture}
[scale=0.50,
inner sep=1.5,font=\ttfamily,
movi/.style={rectangle,draw=black,minimum size=4mm},
vertex/.style={circle,draw=black!50,fill=black!20, very thick},
parta/.style={circle,draw=black!50,fill=black!0, very thick},
partb/.style={circle,draw=black!50,fill=black!100, very thick},
edgevertex/.style={circle,draw=blue!80,fill=blue!40, very thick},
curly/.style={decorate,decoration={snake,amplitude=0.2mm,segment length=1mm,post length=1mm}},
normal/.style=,
post/.style={->, >=stealth'},
pre/.style={<-, >=stealth'}]
\node(a0) at (1.000000,0.000000) {\footnotesize };
\node(a1) at (1.000000,1.000000) {\footnotesize };
\node(a2) at (1.000000,2.000000) {\footnotesize };
\node(a3) at (2.300000,0.000000) {\footnotesize +};
\node(a4) at (2.300000,1.000000) {\footnotesize o};
\node(a5) at (2.300000,2.000000) {\footnotesize -};
\node(a6) at (4.600000,0.000000) {\footnotesize +};
\node(a7) at (4.600000,1.000000) {\footnotesize o};
\node(a8) at (4.600000,2.000000) {\footnotesize -};
\node(a9) at (6.900000,0.000000) {\footnotesize +};
\node(a10) at (6.900000,1.000000) [movi] {\large\textbf{o}};
\node(a11) at (6.900000,2.000000) {\footnotesize -};
\node(a12) at (9.200000,0.000000) {\footnotesize +};
\node(a13) at (9.200000,1.000000) {\footnotesize o};
\node(a14) at (9.200000,2.000000) {\footnotesize -};
\node(a15) at (11.500000,0.000000) {\footnotesize +};
\node(a16) at (11.500000,1.000000) {\footnotesize o};
\node(a17) at (11.500000,2.000000) {\footnotesize -};
\node(a18) at (12.800000,0.000000) {\footnotesize };
\node(a19) at (12.800000,1.000000) {\footnotesize };
\node(a20) at (12.800000,2.000000) {\footnotesize };
\node(a21) at (2.300000,-1.000000) {\footnotesize -1};
\node(a22) at (4.600000,-1.000000) {\footnotesize 0};
\node(a23) at (6.900000,-1.000000) {\footnotesize  1};
\node(a24) at (9.200000,-1.000000) {\footnotesize  2};
\node(a25) at (11.500000,-1.000000) {\footnotesize  3};
\draw [post,black!30,thin] (a3) to (a6);
\draw [post,black!30,thin] (a3) to (a7);
\draw [post,black!30,thin] (a4) to (a6);
\draw [post,black!30,thin] (a4) to (a8);
\draw [post,black!30,thin] (a5) to (a7);
\draw [post,black,ultra thick] (a5) to (a8);
\draw [post,black!30,thin] (a6) to (a9);
\draw [post,black,dashed] (a6) to (a10);
\draw [post,black,dashed] (a6) to (a11);
\draw [post,black!30,thin] (a7) to (a9);
\draw [post,black,curly] (a7) to (a10);
\draw [post,black,curly] (a7) to (a11);
\draw [post,black!30,thin] (a8) to (a9);
\draw [post,black,ultra thick] (a8) to (a10);
\draw [post,black,thick] (a8) to (a11);
\draw [post,black,dashed] (a9) to (a12);
\draw [post,black,dashed] (a9) to (a13);
\draw [post,black,curly] (a10) to (a12);
\draw [post,black,thick] (a10) to (a14);
\draw [post,black,curly] (a11) to (a13);
\draw [post,black,thick] (a11) to (a14);
\draw [post,black,thick] (a12) to (a15);
\draw [post,black,curly] (a12) to (a16);
\draw [post,black,dashed] (a12) to (a17);
\draw [post,black,thick] (a13) to (a15);
\draw [post,black,curly] (a13) to (a16);
\draw [post,black,dashed] (a13) to (a17);
\draw [post,black!30,thin] (a14) to (a15);
\draw [post,black!30,thin] (a14) to (a16);
\draw [post,black!30,thin] (a14) to (a17);
\draw [post,black!30,dotted] (a0) to (a3);
\draw [post,black!30,dotted] (a15) to (a18);
\draw [post,black!30,dotted] (a1) to (a4);
\draw [post,black!30,dotted] (a16) to (a19);
\draw [post,black!30,dotted] (a2) to (a5);
\draw [post,black!30,dotted] (a17) to (a20);
\end{tikzpicture}
}

\def\TubeAor{
\begin{tikzpicture}
[scale=0.50,
inner sep=1.5,font=\ttfamily,
movi/.style={rectangle,draw=black,minimum size=4mm},
vertex/.style={circle,draw=black!50,fill=black!20, very thick},
parta/.style={circle,draw=black!50,fill=black!0, very thick},
partb/.style={circle,draw=black!50,fill=black!100, very thick},
edgevertex/.style={circle,draw=blue!80,fill=blue!40, very thick},
curly/.style={decorate,decoration={snake,amplitude=0.2mm,segment length=1mm,post length=1mm}},
normal/.style=,
post/.style={->, >=stealth'},
pre/.style={<-, >=stealth'}]
\node(a0) at (1.000000,0.000000) {\footnotesize };
\node(a1) at (1.000000,1.000000) {\footnotesize };
\node(a2) at (1.000000,2.000000) {\footnotesize };
\node(a3) at (2.300000,0.000000) {\footnotesize +};
\node(a4) at (2.300000,1.000000) {\footnotesize o};
\node(a5) at (2.300000,2.000000) {\footnotesize -};
\node(a6) at (4.600000,0.000000) {\footnotesize +};
\node(a7) at (4.600000,1.000000) {\footnotesize o};
\node(a8) at (4.600000,2.000000) {\footnotesize -};
\node(a9) at (6.900000,0.000000) {\footnotesize +};
\node(a10) at (6.900000,2.000000) [movi] {\large\textbf{o}};
\node(a11) at (6.900000,1.000000) [movi] {\large\textbf{-}};
\node(a12) at (9.200000,1.000000) [movi] {\large\textbf{+}};
\node(a13) at (9.200000,0.000000) [movi] {\large\textbf{o}};
\node(a14) at (9.200000,2.000000) {\footnotesize -};
\node(a15) at (11.500000,0.000000) {\footnotesize +};
\node(a16) at (11.500000,1.000000) {\footnotesize o};
\node(a17) at (11.500000,2.000000) {\footnotesize -};
\node(a18) at (12.800000,0.000000) {\footnotesize };
\node(a19) at (12.800000,1.000000) {\footnotesize };
\node(a20) at (12.800000,2.000000) {\footnotesize };
\node(a21) at (2.300000,-1.000000) {\footnotesize -1};
\node(a22) at (4.600000,-1.000000) {\footnotesize 0};
\node(a23) at (6.900000,-1.000000) {\footnotesize  1};
\node(a24) at (9.200000,-1.000000) {\footnotesize  2};
\node(a25) at (11.500000,-1.000000) {\footnotesize  3};
\draw [post,black!30,thin] (a3) to (a6);
\draw [post,black!30,thin] (a3) to (a7);
\draw [post,black!30,thin] (a4) to (a6);
\draw [post,black!30,thin] (a4) to (a8);
\draw [post,black!30,thin] (a5) to (a7);
\draw [post,black,ultra thick] (a5) to (a8);
\draw [post,black!30,thin] (a6) to (a9);
\draw [post,black,dashed] (a6) to (a10);
\draw [post,black,dashed] (a6) to (a11);
\draw [post,black!30,thin] (a7) to (a9);
\draw [post,black,curly] (a7) to (a10);
\draw [post,black,curly] (a7) to (a11);
\draw [post,black!30,thin] (a8) to (a9);
\draw [post,black,ultra thick] (a8) to (a10);
\draw [post,black,thick] (a8) to (a11);
\draw [post,black,dashed] (a9) to (a12);
\draw [post,black,dashed] (a9) to (a13);
\draw [post,black,curly] (a10) to (a12);
\draw [post,black,thick] (a10) to (a14);
\draw [post,black,curly] (a11) to (a13);
\draw [post,black,thick] (a11) to (a14);
\draw [post,black,thick] (a12) to (a15);
\draw [post,black,curly] (a12) to (a16);
\draw [post,black,dashed] (a12) to (a17);
\draw [post,black,thick] (a13) to (a15);
\draw [post,black,curly] (a13) to (a16);
\draw [post,black,dashed] (a13) to (a17);
\draw [post,black!30,thin] (a14) to (a15);
\draw [post,black!30,thin] (a14) to (a16);
\draw [post,black!30,thin] (a14) to (a17);
\draw [post,black!30,dotted] (a0) to (a3);
\draw [post,black!30,dotted] (a15) to (a18);
\draw [post,black!30,dotted] (a1) to (a4);
\draw [post,black!30,dotted] (a16) to (a19);
\draw [post,black!30,dotted] (a2) to (a5);
\draw [post,black!30,dotted] (a17) to (a20);
\end{tikzpicture}
}

\def\TubeApl{
\begin{tikzpicture}
[scale=0.50,
inner sep=1.5,font=\ttfamily,
movi/.style={rectangle,draw=black,minimum size=4mm},
vertex/.style={circle,draw=black!50,fill=black!20, very thick},
parta/.style={circle,draw=black!50,fill=black!0, very thick},
partb/.style={circle,draw=black!50,fill=black!100, very thick},
edgevertex/.style={circle,draw=blue!80,fill=blue!40, very thick},
curly/.style={decorate,decoration={snake,amplitude=0.2mm,segment length=1mm,post length=1mm}},
normal/.style=,
post/.style={->, >=stealth'},
pre/.style={<-, >=stealth'}]
\node(a0) at (1.000000,0.000000) {\footnotesize };
\node(a1) at (1.000000,1.000000) {\footnotesize };
\node(a2) at (1.000000,2.000000) {\footnotesize };
\node(a3) at (2.300000,0.000000) {\footnotesize +};
\node(a4) at (2.300000,1.000000) {\footnotesize o};
\node(a5) at (2.300000,2.000000) {\footnotesize -};
\node(a6) at (4.600000,0.000000) {\footnotesize +};
\node(a7) at (4.600000,1.000000) {\footnotesize o};
\node(a8) at (4.600000,2.000000) {\footnotesize -};
\node(a9) at (6.900000,0.000000) [movi] {\large\textbf{+}};
\node(a10) at (6.900000,1.000000) {\footnotesize o};
\node(a11) at (6.900000,2.000000) {\footnotesize -};
\node(a12) at (9.200000,0.000000) {\footnotesize +};
\node(a13) at (9.200000,1.000000) {\footnotesize o};
\node(a14) at (9.200000,2.000000) {\footnotesize -};
\node(a15) at (11.500000,0.000000) {\footnotesize +};
\node(a16) at (11.500000,1.000000) {\footnotesize o};
\node(a17) at (11.500000,2.000000) {\footnotesize -};
\node(a18) at (12.800000,0.000000) {\footnotesize };
\node(a19) at (12.800000,1.000000) {\footnotesize };
\node(a20) at (12.800000,2.000000) {\footnotesize };
\node(a21) at (2.300000,-1.000000) {\footnotesize -1};
\node(a22) at (4.600000,-1.000000) {\footnotesize 0};
\node(a23) at (6.900000,-1.000000) {\footnotesize  1};
\node(a24) at (9.200000,-1.000000) {\footnotesize  2};
\node(a25) at (11.500000,-1.000000) {\footnotesize  3};
\draw [post,black!30,thin] (a3) to (a6);
\draw [post,black!30,thin] (a3) to (a7);
\draw [post,black!30,thin] (a4) to (a6);
\draw [post,black!30,thin] (a4) to (a8);
\draw [post,black!30,thin] (a5) to (a7);
\draw [post,black,ultra thick] (a5) to (a8);
\draw [post,black,dashed] (a6) to (a9);
\draw [post,black!30,thin] (a6) to (a10);
\draw [post,black,dashed] (a6) to (a11);
\draw [post,black,curly] (a7) to (a9);
\draw [post,black!30,thin] (a7) to (a10);
\draw [post,black,curly] (a7) to (a11);
\draw [post,black,ultra thick] (a8) to (a9);
\draw [post,black!30,thin] (a8) to (a10);
\draw [post,black,thick] (a8) to (a11);
\draw [post,black,thick] (a9) to (a12);
\draw [post,black,curly] (a9) to (a13);
\draw [post,black,dashed] (a10) to (a12);
\draw [post,black,dashed] (a10) to (a14);
\draw [post,black,curly] (a11) to (a13);
\draw [post,black,thick] (a11) to (a14);
\draw [post,black,thick] (a12) to (a15);
\draw [post,black,curly] (a12) to (a16);
\draw [post,black,dashed] (a12) to (a17);
\draw [post,black!30,thin] (a13) to (a15);
\draw [post,black!30,thin] (a13) to (a16);
\draw [post,black!30,thin] (a13) to (a17);
\draw [post,black,thick] (a14) to (a15);
\draw [post,black,curly] (a14) to (a16);
\draw [post,black,dashed] (a14) to (a17);
\draw [post,black!30,dotted] (a0) to (a3);
\draw [post,black!30,dotted] (a15) to (a18);
\draw [post,black!30,dotted] (a1) to (a4);
\draw [post,black!30,dotted] (a16) to (a19);
\draw [post,black!30,dotted] (a2) to (a5);
\draw [post,black!30,dotted] (a17) to (a20);
\end{tikzpicture}
}

\def\TubeApr{
\begin{tikzpicture}
[scale=0.50,
inner sep=1.5,font=\ttfamily,
movi/.style={rectangle,draw=black,minimum size=4mm},
vertex/.style={circle,draw=black!50,fill=black!20, very thick},
parta/.style={circle,draw=black!50,fill=black!0, very thick},
partb/.style={circle,draw=black!50,fill=black!100, very thick},
edgevertex/.style={circle,draw=blue!80,fill=blue!40, very thick},
curly/.style={decorate,decoration={snake,amplitude=0.2mm,segment length=1mm,post length=1mm}},
normal/.style=,
post/.style={->, >=stealth'},
pre/.style={<-, >=stealth'}]
\node(a0) at (1.000000,0.000000) {\footnotesize };
\node(a1) at (1.000000,1.000000) {\footnotesize };
\node(a2) at (1.000000,2.000000) {\footnotesize };
\node(a3) at (2.300000,0.000000) {\footnotesize +};
\node(a4) at (2.300000,1.000000) {\footnotesize o};
\node(a5) at (2.300000,2.000000) {\footnotesize -};
\node(a6) at (4.600000,0.000000) {\footnotesize +};
\node(a7) at (4.600000,1.000000) {\footnotesize o};
\node(a8) at (4.600000,2.000000) {\footnotesize -};
\node(a9) at (6.900000,2.000000) [movi] {\large\textbf{+}};
\node(a10) at (6.900000,1.000000) {\footnotesize o};
\node(a11) at (6.900000,0.000000) [movi] {\large\textbf{-}};
\node(a12) at (9.200000,2.000000) [movi] {\large\textbf{+}};
\node(a13) at (9.200000,1.000000) {\footnotesize o};
\node(a14) at (9.200000,0.000000) [movi] {\large\textbf{-}};
\node(a15) at (11.500000,0.000000) {\footnotesize +};
\node(a16) at (11.500000,1.000000) {\footnotesize o};
\node(a17) at (11.500000,2.000000) {\footnotesize -};
\node(a18) at (12.800000,0.000000) {\footnotesize };
\node(a19) at (12.800000,1.000000) {\footnotesize };
\node(a20) at (12.800000,2.000000) {\footnotesize };
\node(a21) at (2.300000,-1.000000) {\footnotesize -1};
\node(a22) at (4.600000,-1.000000) {\footnotesize 0};
\node(a23) at (6.900000,-1.000000) {\footnotesize  1};
\node(a24) at (9.200000,-1.000000) {\footnotesize  2};
\node(a25) at (11.500000,-1.000000) {\footnotesize  3};
\draw [post,black!30,thin] (a3) to (a6);
\draw [post,black!30,thin] (a3) to (a7);
\draw [post,black!30,thin] (a4) to (a6);
\draw [post,black!30,thin] (a4) to (a8);
\draw [post,black!30,thin] (a5) to (a7);
\draw [post,black,ultra thick] (a5) to (a8);
\draw [post,black,dashed] (a6) to (a9);
\draw [post,black!30,thin] (a6) to (a10);
\draw [post,black,dashed] (a6) to (a11);
\draw [post,black,curly] (a7) to (a9);
\draw [post,black!30,thin] (a7) to (a10);
\draw [post,black,curly] (a7) to (a11);
\draw [post,black,ultra thick] (a8) to (a9);
\draw [post,black!30,thin] (a8) to (a10);
\draw [post,black,thick] (a8) to (a11);
\draw [post,black,thick] (a9) to (a12);
\draw [post,black,curly] (a9) to (a13);
\draw [post,black,dashed] (a10) to (a12);
\draw [post,black,dashed] (a10) to (a14);
\draw [post,black,curly] (a11) to (a13);
\draw [post,black,thick] (a11) to (a14);
\draw [post,black,thick] (a12) to (a15);
\draw [post,black,curly] (a12) to (a16);
\draw [post,black,dashed] (a12) to (a17);
\draw [post,black!30,thin] (a13) to (a15);
\draw [post,black!30,thin] (a13) to (a16);
\draw [post,black!30,thin] (a13) to (a17);
\draw [post,black,thick] (a14) to (a15);
\draw [post,black,curly] (a14) to (a16);
\draw [post,black,dashed] (a14) to (a17);
\draw [post,black!30,dotted] (a0) to (a3);
\draw [post,black!30,dotted] (a15) to (a18);
\draw [post,black!30,dotted] (a1) to (a4);
\draw [post,black!30,dotted] (a16) to (a19);
\draw [post,black!30,dotted] (a2) to (a5);
\draw [post,black!30,dotted] (a17) to (a20);
\end{tikzpicture}
}

\def\TubeAOl{
\begin{tikzpicture}
[scale=0.50,
inner sep=1.5,font=\ttfamily,
movi/.style={rectangle,draw=black,minimum size=4mm},
vertex/.style={circle,draw=black!50,fill=black!20, very thick},
parta/.style={circle,draw=black!50,fill=black!0, very thick},
partb/.style={circle,draw=black!50,fill=black!100, very thick},
edgevertex/.style={circle,draw=blue!80,fill=blue!40, very thick},
curly/.style={decorate,decoration={snake,amplitude=0.2mm,segment length=1mm,post length=1mm}},
normal/.style=,
post/.style={->, >=stealth'},
pre/.style={<-, >=stealth'}]
\node(a0) at (1.000000,0.000000) {\footnotesize };
\node(a1) at (1.000000,1.000000) {\footnotesize };
\node(a2) at (1.000000,2.000000) {\footnotesize };
\node(a3) at (2.300000,0.000000) {\footnotesize +};
\node(a4) at (2.300000,1.000000) {\footnotesize o};
\node(a5) at (2.300000,2.000000) {\footnotesize -};
\node(a6) at (4.600000,0.000000) {\footnotesize +};
\node(a7) at (4.600000,1.000000) {\footnotesize o};
\node(a8) at (4.600000,2.000000) {\footnotesize -};
\node(a9) at (6.900000,0.000000) {\footnotesize +};
\node(a10) at (6.900000,1.000000) {\footnotesize o};
\node(a11) at (6.900000,2.000000) {\footnotesize -};
\node(a12) at (9.200000,0.000000) {\footnotesize +};
\node(a13) at (9.200000,1.000000) [movi] {\large\textbf{o}};
\node(a14) at (9.200000,2.000000) {\footnotesize -};
\node(a15) at (11.500000,0.000000) {\footnotesize +};
\node(a16) at (11.500000,1.000000) {\footnotesize o};
\node(a17) at (11.500000,2.000000) {\footnotesize -};
\node(a18) at (12.800000,0.000000) {\footnotesize };
\node(a19) at (12.800000,1.000000) {\footnotesize };
\node(a20) at (12.800000,2.000000) {\footnotesize };
\node(a21) at (2.300000,-1.000000) {\footnotesize -2};
\node(a22) at (4.600000,-1.000000) {\footnotesize -1};
\node(a23) at (6.900000,-1.000000) {\footnotesize  0};
\node(a24) at (9.200000,-1.000000) {\footnotesize  1};
\node(a25) at (11.500000,-1.000000) {\footnotesize  2};
\draw [post,black!30,thin] (a3) to (a6);
\draw [post,black!30,thin] (a3) to (a7);
\draw [post,black!30,thin] (a4) to (a6);
\draw [post,black!30,thin] (a4) to (a8);
\draw [post,black!30,thin] (a5) to (a7);
\draw [post,black!30,thin] (a5) to (a8);
\draw [post,black,thick] (a6) to (a9);
\draw [post,black,thick] (a6) to (a10);
\draw [post,black!30,thin] (a6) to (a11);
\draw [post,black,dashed] (a7) to (a9);
\draw [post,black,dashed] (a7) to (a10);
\draw [post,black!30,thin] (a7) to (a11);
\draw [post,black,curly] (a8) to (a9);
\draw [post,black,curly] (a8) to (a10);
\draw [post,black,ultra thick] (a8) to (a11);
\draw [post,black,dashed] (a9) to (a12);
\draw [post,black,curly] (a9) to (a13);
\draw [post,black,dashed] (a10) to (a12);
\draw [post,black,curly] (a10) to (a14);
\draw [post,black,ultra thick] (a11) to (a13);
\draw [post,black,thick] (a11) to (a14);
\draw [post,black!30,thin] (a12) to (a15);
\draw [post,black!30,thin] (a12) to (a16);
\draw [post,black!30,thin] (a12) to (a17);
\draw [post,black,thick] (a13) to (a15);
\draw [post,black,curly] (a13) to (a16);
\draw [post,black,dashed] (a13) to (a17);
\draw [post,black,thick] (a14) to (a15);
\draw [post,black,curly] (a14) to (a16);
\draw [post,black,dashed] (a14) to (a17);
\draw [post,black!30,dotted] (a0) to (a3);
\draw [post,black!30,dotted] (a15) to (a18);
\draw [post,black!30,dotted] (a1) to (a4);
\draw [post,black!30,dotted] (a16) to (a19);
\draw [post,black!30,dotted] (a2) to (a5);
\draw [post,black!30,dotted] (a17) to (a20);
\end{tikzpicture}
}

\def\TubeAOr{
\begin{tikzpicture}
[scale=0.50,
inner sep=1.5,font=\ttfamily,
movi/.style={rectangle,draw=black,minimum size=4mm},
vertex/.style={circle,draw=black!50,fill=black!20, very thick},
parta/.style={circle,draw=black!50,fill=black!0, very thick},
partb/.style={circle,draw=black!50,fill=black!100, very thick},
edgevertex/.style={circle,draw=blue!80,fill=blue!40, very thick},
curly/.style={decorate,decoration={snake,amplitude=0.2mm,segment length=1mm,post length=1mm}},
normal/.style=,
post/.style={->, >=stealth'},
pre/.style={<-, >=stealth'}]
\node(a0) at (1.000000,0.000000) {\footnotesize };
\node(a1) at (1.000000,1.000000) {\footnotesize };
\node(a2) at (1.000000,2.000000) {\footnotesize };
\node(a3) at (2.300000,0.000000) {\footnotesize +};
\node(a4) at (2.300000,1.000000) {\footnotesize o};
\node(a5) at (2.300000,2.000000) {\footnotesize -};
\node(a6) at (4.600000,0.000000) {\footnotesize +};
\node(a7) at (4.600000,1.000000) {\footnotesize o};
\node(a8) at (4.600000,2.000000) {\footnotesize -};
\node(a9) at (6.900000,1.000000) [movi] {\large\textbf{+}};
\node(a10) at (6.900000,0.000000) [movi] {\large\textbf{o}};
\node(a11) at (6.900000,2.000000) {\footnotesize -};
\node(a12) at (9.200000,0.000000) {\footnotesize +};
\node(a13) at (9.200000,2.000000) [movi] {\large\textbf{o}};
\node(a14) at (9.200000,1.000000) [movi] {\large\textbf{-}};
\node(a15) at (11.500000,0.000000) {\footnotesize +};
\node(a16) at (11.500000,1.000000) {\footnotesize o};
\node(a17) at (11.500000,2.000000) {\footnotesize -};
\node(a18) at (12.800000,0.000000) {\footnotesize };
\node(a19) at (12.800000,1.000000) {\footnotesize };
\node(a20) at (12.800000,2.000000) {\footnotesize };
\node(a21) at (2.300000,-1.000000) {\footnotesize -2};
\node(a22) at (4.600000,-1.000000) {\footnotesize -1};
\node(a23) at (6.900000,-1.000000) {\footnotesize  0};
\node(a24) at (9.200000,-1.000000) {\footnotesize  1};
\node(a25) at (11.500000,-1.000000) {\footnotesize  2};
\draw [post,black!30,thin] (a3) to (a6);
\draw [post,black!30,thin] (a3) to (a7);
\draw [post,black!30,thin] (a4) to (a6);
\draw [post,black!30,thin] (a4) to (a8);
\draw [post,black!30,thin] (a5) to (a7);
\draw [post,black!30,thin] (a5) to (a8);
\draw [post,black,thick] (a6) to (a9);
\draw [post,black,thick] (a6) to (a10);
\draw [post,black!30,thin] (a6) to (a11);
\draw [post,black,dashed] (a7) to (a9);
\draw [post,black,dashed] (a7) to (a10);
\draw [post,black!30,thin] (a7) to (a11);
\draw [post,black,curly] (a8) to (a9);
\draw [post,black,curly] (a8) to (a10);
\draw [post,black,ultra thick] (a8) to (a11);
\draw [post,black,dashed] (a9) to (a12);
\draw [post,black,curly] (a9) to (a13);
\draw [post,black,dashed] (a10) to (a12);
\draw [post,black,curly] (a10) to (a14);
\draw [post,black,ultra thick] (a11) to (a13);
\draw [post,black,thick] (a11) to (a14);
\draw [post,black!30,thin] (a12) to (a15);
\draw [post,black!30,thin] (a12) to (a16);
\draw [post,black!30,thin] (a12) to (a17);
\draw [post,black,thick] (a13) to (a15);
\draw [post,black,curly] (a13) to (a16);
\draw [post,black,dashed] (a13) to (a17);
\draw [post,black,thick] (a14) to (a15);
\draw [post,black,curly] (a14) to (a16);
\draw [post,black,dashed] (a14) to (a17);
\draw [post,black!30,dotted] (a0) to (a3);
\draw [post,black!30,dotted] (a15) to (a18);
\draw [post,black!30,dotted] (a1) to (a4);
\draw [post,black!30,dotted] (a16) to (a19);
\draw [post,black!30,dotted] (a2) to (a5);
\draw [post,black!30,dotted] (a17) to (a20);
\end{tikzpicture}
}

\def\DOTS[#1]{
\begin{tikzpicture}[scale=#1]
\draw [dotted] (0,0) to (1,0);
\end{tikzpicture}
}

\def\MATCHvier[#1]{
\begin{tikzpicture}[scale=#1]
\draw (0.05,-1) to (0.95,-1);
\draw (0.05,-0.5) to (0.95,-0.5);
\draw [dotted] (0.5,-0.3) to (0.5,0.3);
\draw (0.05,0.5) to (0.95,0.5);
\draw (0.05,1) to (0.95,1);
\end{tikzpicture}
}
\def\MATCHsechs[#1]{
\begin{tikzpicture}[scale=#1]
\draw (0.05,-1.5) to (0.95,-1.5);
\draw (0.05,-1) to (0.95,-1);
\draw [dotted] (0.5,-0.6) to (0.5,0.6);
\draw (0.05,1) to (0.95,1);
\draw (0.05,1.5) to (0.95,1.5);
\end{tikzpicture}
}

\def\GGvier[#1,#2]{
\begin{tikzpicture}[scale=#1]
\draw (0.05,-1) to (0.95,-1);
\draw (0.95,-1) to (0.95,1);
\draw (0.05,1) to (0.95,1);
\draw (0.05,-1) to (0.05,1);
\node at (0.5,0) {#2};
\end{tikzpicture}
}
\def\GG[#1,#2,#3,#4]{
\begin{tikzpicture}[scale=#1]
\draw (0.05,-#3) to (0.95,-#4);
\draw (0.95,-#4) to (0.95,#4);
\draw (0.05,#3) to (0.95,#4);
\draw (0.05,-#3) to (0.05,#3);
\node at (0.5,0) {#2};
\end{tikzpicture}
}

\def\McKAY[#1]{
\begin{tikzpicture}[scale=#1]
\node at (0,0) {\DOTS[#1]};
\node at (1,0) {\MATCHvier[#1]};
\node at (2,0) {\GGvier[#1,$K_{n,n}$]};
\node at (3,0) {\MATCHvier[#1]};
\node at (4,0) {\DOTS[#1]};
\node at (5,0) {\MATCHvier[#1]};
\node at (6,0) {\GGvier[#1,$K_{n,n}$]};
\node at (7,0) {\MATCHvier[#1]};
\node at (8,0) {\DOTS[#1]};
\end{tikzpicture}
}

\def\LITTLE[#1]{
\begin{tikzpicture}[scale=#1]
\node at (0,0) {\DOTS[#1]};
\node at (1,0) {\GGvier[#1,$I$]};
\node at (2,0) {\GGvier[#1,$K_{n,n}$]};
\node at (3,0) {\GGvier[#1,$I$]};
\node at (4,0) {\GGvier[#1,$K_{n,n}$]};
\node at (5,0) {\GGvier[#1,$I$]};
\node at (6,0) {\GGvier[#1,$K_{n,n}$]};
\node at (7,0) {\DOTS[#1]};
\end{tikzpicture}
}

\def\MOREGENERAL[#1]{
\begin{tikzpicture}[scale=#1]
\node at (0,0) {\DOTS[#1]};
\node at (1,0) {\GGvier[#1,$K_{n,n}$]};
\node at (2,0) {\MATCHvier[#1]};
\node at (3,0) {\GGvier[#1,$I_1$]};
\node at (4,0) {\GGvier[#1,$K_{n,n}$]};
\node at (5,0) {\GGvier[#1,$K_{n,n}$]};
\node at (6,0) {\MATCHvier[#1]};
\node at (7,0) {\GGvier[#1,$I_2$]};
\node at (8,0) {\MATCHvier[#1]};
\node at (9,0) {\GGvier[#1,$K_{n,n}$]};
\node at (10,0) {\DOTS[#1]};
\end{tikzpicture}
}

\def\GENERAL[#1]{
\begin{tikzpicture}[scale=#1]
\node at (0,0) {\DOTS[#1]};
\node at (1,0) {\GG[#1,$K_{n,m}$,1,1.5]};
\node at (2,0) {\MATCHsechs[#1]};
\node at (3,0) {\GG[#1,$I_1$,1.5,0.5]};
\node at (4,0) {\GG[#1,$K_{r,r}$,0.5,0.5]};
\node at (5,0) {\GG[#1,$K_{r,s}$,0.5,1]};
\node at (6,0) {\MATCHvier[#1]};
\node at (7,0) {\GG[#1,$I_2$,1,1]};
\node at (8,0) {\MATCHvier[#1]};
\node at (9,0) {\GG[#1,$K_{n,m}$,1,1.5]};
\node at (10,0) {\DOTS[#1]};
\end{tikzpicture}
}

\def\LOCAL[#1]{
\begin{tikzpicture}[scale=#1]
\node at (0,0) {\GG[#1,$K_{n,m}$,0.5,1]};
\node at (1,0) {\MATCHvier[#1]};
\node at (2,0) {\DOTS[#1]};
\node at (3,0) {\MATCHvier[#1]};
\node at (4,0) {\GG[#1,$I_1$,1,1.5]};
\node at (5,0) {\MATCHsechs[#1]};
\node at (6,0) {\DOTS[#1]};
\node at (7,0) {\MATCHsechs[#1]};
\node at (8,0) {\GG[#1,$K_{r,s}$,1.5,0.5]};
\end{tikzpicture}
}

\def\MULTI[#1]{
\begin{tikzpicture}[scale=#1]
\node at (0,0) {\GG[#1,$K_{n,m}$,0.5,1]};
\node at (1,0) {\MATCHvier[#1]};
\node at (2,0) {\DOTS[#1]};
\node at (3,0) {\MATCHvier[#1]};
\node at (4,0) {\GG[#1,$I_1$,1,1.5]};
\node at (5,0) {\MATCHsechs[#1]};
\node at (6,0) {\DOTS[#1]};
\node at (7,0) {\MATCHvier[#1]};
\node at (8,0) {\GG[#1,$I_t$,1,1]};
\node at (9,0) {\MATCHvier[#1]};
\node at (10,0) {\DOTS[#1]};
\node at (11,0) {\MATCHvier[#1]};
\node at (12,0) {\GG[#1,$K_{r,s}$,1,1.5]};
\end{tikzpicture}
}

\def\INF[#1]{
\begin{tikzpicture}[scale=#1]
\node at (0,0) {\DOTS[#1]};
\node at (1,0) {\MATCHvier[#1]};
\node at (2,0) {\MATCHvier[#1]};
\node at (3,0) {\GG[#1,$I$,1,1.5]};
\node at (4,0) {\MATCHsechs[#1]};
\node at (5,0) {\MATCHsechs[#1]};
\node at (6,0) {\DOTS[#1]};
\end{tikzpicture}
}

\begin{document}

\newbox\Adr
\setbox\Adr\vbox{
\centerline{\sc Christoph Neumann}
\vskip18pt
\centerline{Faculty of Mathematics,
University of Vienna}
}

\title[Layerwise Constructions]{Constructing highly arc transitive digraphs using a layerwise direct product}
\author[Christoph Neumann]{\box\Adr}
\address{Faculty of Mathematics, University of Vienna, Nordbergstrasse~15, A-1090 Vienna, Austria}

\begin{abstract}
We introduce a construction of highly arc transitive digraphs using a layerwise direct product. This product generalizes some known classes of highly arc transitive digraphs but also allows to construct new such. We use the product to obtain counterexamples to a conjecture by Cameron, Praeger and Wormald on the structure of certain highly arc transitive digraphs.
\end{abstract}
\maketitle

\secA{Introduction}\label{sec.intro}

A digraph is highly arc transitive if its automorphism group acts transitively on the set of its $n$-arcs for every $n$. We present a way to construct various highly arc transitive digraphs as a layerwise product. Doing so, we unify different constructions of highly arc transitive digraphs presented in \cite{1}, \cite{9} and \cite{7} and obtain new highly arc transitive digraphs. Depending on the number and structure of the factors, the produced digraphs can have one or two ends.

\names{Cameron}, \names{Praeger} and \names{Wormald} studied highly arc transitive digraphs in \cite{1} and made a conjecture about the subclass of digraphs which have an epimorphism onto the integer line such that the preimages are finite. \names{Möller} \cite{2} constructed a digraph, that was believed to be a counterexample, but later \names{\v{S}parl} \cite{5} showed that it is not transitive. We give digraphs that are layerwise direct products and qualify as counterexamples. Very recently \names{DeVos}, \names{Mohar} and \names{\v{S}\'amal} studied highly arc transitive digraphs in \cite{9}. They answer a question from \cite{1} and clarify the mentioned conjecture using the same counterexamples as we do, but obtained them with a different approach. Moreover they study the structure of highly arc transitive digraphs with two ends. In Section~\ref{sec.art} we explain the relation between the present work and \cite{9}.

Our paper is organized as follows. In Section~\ref{sec.hat} we recall some notions from \cite{1}. In Section~\ref{sec.prod} we introduce the layerwise direct product. In Section~\ref{sec.counter} we use this product to find a counterexample to the mentioned conjecture.  in Section~\ref{sec.rem} we make some remarks on our overlap with \cite{9}, \names{Möller}'s digraph and the automorphisms we needed for our construction in Section~\ref{sec.counter}. In Section~\ref{sec.more} we present further highly arc transitive digraphs that can be constructed using the layerwise direct product. This includes representations of some know highly arc transitive digraphs, generalizations of the factors and the thereby obtained new highly arc transitive digraphs.

\secB{Highly arc transitive digraphs.}\label{sec.hat}
An \defstil{$n$-arc} in a digraph $D$ is a series $(e_i)_{i=0,\dots,n-1}$ of $n$ directed edges $e_i=(x_{e_i},y_{e_i})$ such that $y_{e_i}=x_{e_{i+1}}$ for all $i$ in the range. A digraph is called \defstil{$n$-arc transitive} if his automorphism group acts transitively on the nonempty set of its $n$-arcs. A digraph is called \defstil{highly arc transitive (HAT)} if it is $n$-arc transitive for all $n\in\N$. An \defstil{alternating walk} is a series $(e_i)_{i=0,\dots,n-1}$ of $n$ directed edges $e_i=(x_{e_i},y_{e_i})$ such that consecutive edges agree alternatingly on either the initial or terminal vertex, i.e., for every even $i$ in the range $x_{e_i}=x_{e_{i+1}}$ and $y_{e_i}=y_{e_{i+1}}$ for every odd $i$ (or vice verse). Two edges are \defstil{reachable} from each other if there exists an alternating walk that contains them both. Being reachable from each other is an equivalence relation on the edgeset of a digraph which we call \defstil{reachability relation}. Let $D=(V,E)$ be a HAT digraph and $e\in E$. Let $\Delta(e)$ be the subgraph of $D$ that is spanned by the equivalence class of $e$ with respect to the reachability relation. $\Delta(e)$ is independent from $e$ i.e. $\forall e_1,e_2\in E:\ \Delta(e_1)\cong \Delta(e_2)$. Thus, we can speak of it as $\Delta(D)$ and call it the \defstil{associated digraph} of $D$. If $\Delta(D)\neq D$ then $\Delta(D)$ is bipartite. Let $Z=(\Z,\{(i,i+1)\mid i\in\Z\})$ be the digraph representing the integer line. A digraph $G$ has \defstil{Property $Z$} if there is an epimorphism $\varphi:G\to Z$. Note that if $G$ is HAT all edges of $\Delta(e)$ are mapped to the same edge $\varphi(e)$ by $\varphi$. The inverse image $\varphi^{-1}((i,i+1))$ of every edge of $Z$ induces a subgraph of $G$ which we call \defstil{layer}. If $G$ is HAT all components of every layer are isomorphic to $\Delta(G)$. For more details on these notions we refer to \cite{1} or \cite{4}.

In the paper \cite{1} on HAT digraphs, \names{Cameron}, \names{Praeger} and \names{Wormald} stated the following conjecture.

\begin{conj}[\textsc{Cameron, Praeger, Wormald}]\label{conj}
Let $D$ be a connected HAT digraph with Property $Z$ and finite fibers $\varphi^{-1}(x)$. Then the associated digraph $\Delta(D)$ is complete bipartite.
\end{conj}

\secB{The layerwise direct product}\label{sec.prod}

The layerwise direct product appeared earlier e.g. in the construction of the Diestel--Leader digraph (also known as broom-digraph). It was generalized to the horocyclic product which is defined in \cite{8}. This generalization is related, but does not agree with the one that we are going to use. Our layerwise direct product is a proper subgraph of the direct product. It is induced by a subset of the vertexset that is gained in the following way: Given Property $Z$ for both factors, one restricts to vertices agreeing on their image.

\begin{defn}[\textsc{Layerwise direct product}]
Let $G^1=(V^1,E^1)$ and $G^2=(V^2,E^2)$ be digraphs with Property $Z$ and let $\varphi^1:V^1\to Z$ and $\varphi^2:V^2\to Z$ be the arising epimorphisms. The \defstil{layerwise direct product} $G^1 \,{}_{\varphi^1}\!\!\times_{\varphi^2} G^2$ is the digraph $(V,E)$ with
\begin{eqnarray*}
V &=& \{(x,y)\mid x\in V^1,\, y\in V^2,\, \varphi^1(x)=\varphi^2(y)\},\\
E &=& \{((a,b),(x,y))\mid (a,x)\in E^1,\, (b,y)\in E^2\}.
\end{eqnarray*}
\end{defn}

If the factors are both connected HAT digraphs with Property $Z$, this layerwise direct product gives exactly a connected component of the ordinary direct product -- in this case these components are all isomorphic. For general factors, there is no need for the components to be isomorphic. So, in the general case, the layerwise direct product picks the ``central'' component.

It is often convenient to denote the vertices as $(n,x,y)$ where $n=\varphi^1(x)=\varphi^2(y)$. In the situation $x=(i,a)$ and $y=(j,b)$ and $\varphi^{1,2}$ depend only on $i$ and $j$ respectively, it is very convenient to denote the vertices as $(n,a,b)$.

If we need to consider more factors, we generalize our definition of the layerwise direct product as follows: Let $I$ be some set of indices and $G_i=(V^i,E^i)$ for $i\in I$ digraphs with Property $Z$ and let $\Phi\coloneqq\{\varphi^i:V^i\to Z\mid i\in I\}$ be the set of arising epimorphisms. The \defstil{layerwise direct product} $\prod_I^\Phi\coloneqq G^i$ is the digraph $(V,E)$ with
\begin{eqnarray*}
V &=& \{(x_i)_{i\in I}\mid x_i\in V^i,\ \forall j,k\in I:\ \varphi^j(x_j)=\varphi^k(x_k)\},\\
E &=& \{((x_i)_{i\in I},(y_i)_{i\in I})\mid \forall j\in I:\ (x_j,y_j)\in E^j\}.
\end{eqnarray*}

Again we denote the vertices as $(n,(x_i)_{i\in I})$ or $(n,(a_i)_{i\in I})$ respectively.

\secA{The counterexample}\label{sec.counter}

\begin{figure}[htdg]
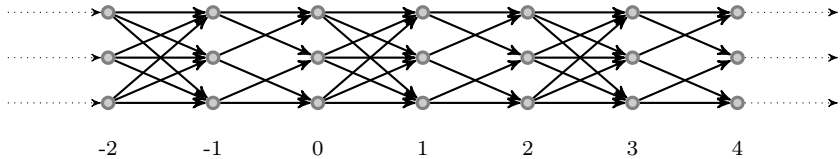

\centering
\TubeL
\caption{The digraph $L$}\label{fig.L}
\end{figure}

Our counterexample for Conjecture~\ref{conj} will be a layerwise direct product of two factors which we have to define first. Let $L=(V^L,E^L)$ be the digraph with vertexset $V^L=\Z\times\{-1,0,1\}$ and $e=((i,x),(i+1,y))$ be an edge if $i$ is even or if $i$ is odd and $x+y\neq 0$ (the latter is a nice way to encode an alternating $6$-cycle). Thus the layers of $L$ are alternating a $K_{3,3}$ and an alternating $6$-cycle (see Figure~\ref{fig.L}). 

$L$ has Property $Z$ with $\varphi^1((i,x))\coloneqq i$ but $\varphi^2((i,x))\coloneqq i+1$ is an epimorphism onto $Z$ too. We can now define our counterexample:
\begin{equation}
D\coloneqq L\,{}_{\varphi^1}\!\!\times_{\varphi^2}L.\label{eq.D}
\end{equation}
For convenience, we denote $-1$/$0$/$1$ with \minus/\oo/\plus. As mentioned above, we extract the $\Z$ coordinate and denote the vertices of $D$ with $(n,x,y)$ rather than with $((n,x),(n+1,y))$ as if $D$ had the vertexset $\Z\times\{$\minus$,$\oo$,$\plus$\}\times\{$\minus$,$\oo$,$\plus$\}$.

\begin{figure}[htdg]
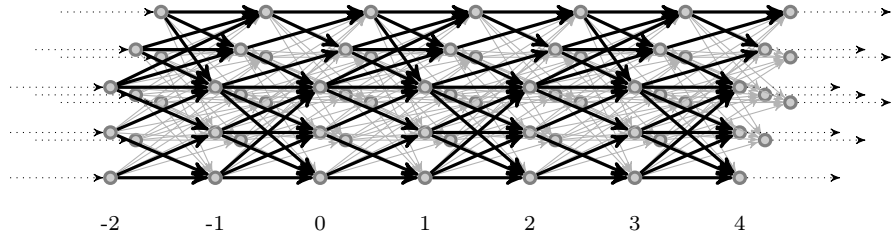

\centering
\TubeD
\caption{The digraph $D$}\label{fig.D}
\end{figure}

\begin{lemma}\label{lemma}
Let $AC_6$ be the alternating $6$-cycle and $D=(V(D),E(D))$ the above digraph. Let $\psi^1:V(K_{3,3})\to\{0,1\}$ and $\psi^2:AC_6\to \{0,1\}$ both map the initial vertices to $0$ and the terminal vertices to $1$. Then $\Delta(e)$ is independent of $e$ thus $\Delta(D)$ is welldefined and
\begin{equation}
\Delta(D)=K_{3,3}\,{}_{\psi^1}\!\!\times_{\psi^2}AC_6 \label{eq.delta}
\end{equation}
is connected, $1$-arc transitive, bipartite but not complete bipartite.
\end{lemma}

\begin{figure}[htdg]
\centering
\TubeDeltaAAA\quad\TubeDeltaABA\quad\TubeDeltaACA\\
\TubeDeltaBAA\quad\TubeDeltaBBA\quad\TubeDeltaBCA\\
\TubeDeltaCAA\quad\TubeDeltaCBA\quad\TubeDeltaCCA\\
\TubeDeltaAAB\quad\TubeDeltaABB\quad\TubeDeltaACB\\
\TubeDeltaBAB\quad\TubeDeltaBBB\quad\TubeDeltaBCB\\
\TubeDeltaCAB\quad\TubeDeltaCBB\quad\TubeDeltaCCB
\caption{$\Delta(D)$}\label{fig.delta}
\end{figure}

\begin{proof}
For $x=(i,x^1,x^2)\in V(D)$ define $\varphi(x)\coloneqq i=\varphi^1(x^1)=\varphi^2(x^2)$. Thus, the layerwise direct product preserves Property $Z$. Thus no alternating walk can leave the layer in which it started. Thus $\Delta(D)$ (if it exists) must be isomorphic to a subgraph of a layer thus bipartite. By a flip of the coordinates the layerwise direct product is commutative. Thus, we do not need to distinguish between even and odd layers in the following. By the construction of $D$, every layer is of the form of the right side of $(\ref{eq.delta})$. If we consider the subgraphs of a layer that are spanned by the vertices that agree on the second coordinate (note that that means ignoring the $\Z$-coordinate) we get three alternating $6$--cycles. If we consider the subgraph spanned by the vertices that have a \minus\ in their third coordinate we get a $K_{3,3}$ that connects the three cycles. Thus, the layer is connected and $\Delta(D)$ exists since all layers are isomorphic. Thus, $(\ref{eq.delta})$ holds. The $1$-arc transitivity follows from the $1$-arc transitivity of the factors in the following way: Every edge in one of the factors corresponds to a copy of the other factor in the product. The edgesets of the collection of these copies form a partition of the edgeset of the product (compare Figure~\ref{fig.delta}). By the $1$-arc transitivity of each one factor its automorphism group acts transitively on these copies and can be embedded in the automorphism group of $\Delta(D)$. Thus, one can adjust in the first step the second coordinates of the endpoints of an edge and in the second step the third coordinates -- only restricted by the existence of edges. This is also the main idea of the construction of $D$ and will be used in pretty much every step in the proof of Proposition~\ref{prop}. Finally the edge $((0,$\oo$,$\oo$),(1,$\oo$,$\oo$))$ does not exist in $\Delta(D)$ and thus it is not complete bipartite.
\end{proof}

\begin{prop}\label{prop}
The digraph $D$ defined in $(\ref{eq.D})$ is HAT.
\end{prop}

\begin{proof}
We are going to prove very carefully that for every $n$-arc $a$ (assume initial vertex $x$ with $\varphi(x)=-m$) in $D$ there is an automorphism $\psi\in\Aut(D)$ such that $\psi(a)=z$ where
\begin{equation*}
z=(((i,\text{\minus},\text{\minus}),(i+1,\text{\minus},\text{\minus})))_{i=-n,\dots,-1}.
\end{equation*}

\begin{figure}[htdg]
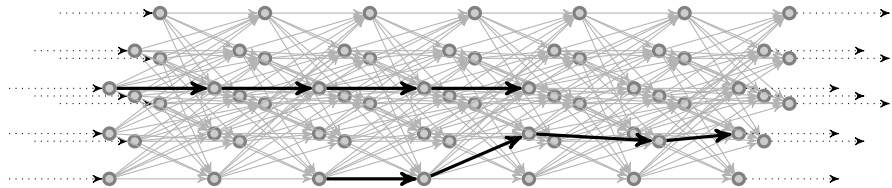

\centering
\TubeAZ
\caption{Two $4$-arcs in $D$}
\end{figure}

Therefore, we define a shift-automorphism $\sigma$, transitivity-automorphisms $\theta_\text{\oo}$, $\theta_\text{\plus}$, $\theta^\text{\oo}$ and $\theta^\text{\plus}$ and arranging-automorphisms $\alpha_\text{\oo}$, $\alpha_\text{\plus}$ and $\alpha^\text{\oo}$. For convenience, we denote $\theta_\text{\minus}=\theta^\text{\minus}=\alpha_\text{\minus}=\alpha^\text{\minus}=\alpha^\text{\plus}=\id$. Our automorphism $\psi$ will then have the form
\begin{equation}\label{eq.psi}
\psi\coloneqq(\sigma^{-1}\circ\alpha_{i_n}\circ\alpha^{i^n})\circ\dots\circ(\sigma^{-1}\circ\alpha_{i_1}\circ\alpha^{i^1})\circ\theta_{i_0}\circ\theta^{i^0}\circ\sigma^m,
\end{equation}
where $m\in\Z$ is of course meant to be a power rather than an index and $i_j,i^j\in\{$\minus$,$\oo$,$\plus$\}$ where $i^0$ cannot take the value \plus\ because there are no edges $((0,x,$\minus$),(1,x,$\plus$))$. We are now going to explain the tasks of the above automorphisms and afterwards check their existence and correctness one by one. In the following, we will refer to the doubleray with the vertices $(i,$\minus$,$\minus$)$ as \defstil{baseline}.
\begin{itemize}
\item $\sigma$ shifts $D$ up by one layer i.e. $\varphi+1\equiv\varphi\circ\sigma$. Thus $\sigma^m$ shifts the initial vertex of $a$ to the zero-layer i.e. $\varphi\circ\sigma^m(x)=0$ (vertex-layers can be defined analogously to the edge-layers we used above, since it is clear from the context which layers we mean, we are not wasting more notation on them). Moreover $\sigma$ stabilizes the baseline setwise.
\item The transitivity-automorphisms guarantee the transitivity of $D$. Their task is to map $\sigma^m(x)$ to the baseline. $\theta_\text{\plus}$ and $\theta_\text{\oo}$ map $(0,$\plus$,y)$ and $(0,$\oo$,y)$ respectively to $(0,$\minus$,y)$, thereby fixing $y$. $\theta^\text{\plus}$ and $\theta^\text{\oo}$ map $(0,x,$\oo$)$ and $(0,x,$\plus$)$ respectively to $(0,x,$\minus$)$, thereby fixing $x$. Thus, they can be combined in a way to map any vertex in the zero-layer to $(0,$\minus$,$\minus$)$. There is no need for the transitivity-automorphisms to stabilize the baseline.
\item The arranging-automorphisms guarantee that $D$ is $(s+1)$-arc transitive if it only is $s$-arc transitive (where we can understand the transitivity gained from the transitivity-automorphisms as $0$-arc transitivity). Therefore, we must be able to map the arcs $((0,$\minus$,$\minus$),(1,x,y))$ onto $((0,$\minus$,$\minus$),(1,$\minus$,$\minus$))$ and at the same time stabilize the negative half of the baseline. $\alpha^\text{\oo}$ maps $((0,$\minus$,$\minus$),(1,$\oo$,y))\mapsto((0,$\minus$,$\minus$),(1,$\minus$,y))$ again fixing $y$. $\alpha_\text{\oo}$ and $\alpha_\text{\minus}$ respectively map $((0,$\minus$,$\minus$),(1,x,$\oo$))$ and $((0,$\minus$,$\minus$),(1,x,$\plus$))$ respectively to $((0,$\minus$,$\minus$),(1,x,$\minus$))$ jet again fixing $x$.
\item The brackets in $(\ref{eq.psi})$ therefore map the arc with initial vertex $(0,$\minus$,$\minus$)$ on the baseline and shift $D$ one step to the left to keep the working-layer the same.
\end{itemize}
We are now going to have a closer look at every single of these automorphisms:
\begin{enumerate}
\item[$\sigma$:] We start with the shift-automorphism. It will increase the $\Z$ coordinate by one and flip the other coordinates.
\begin{equation*}
\sigma:(i,x,y)\mapsto (i+1,y,x)
\end{equation*}
We recognize that we can describe the edgeset of $D$ by joining the conditions of the factors. Thus
\begin{equation*}
((i,x_1,y_1),(i+1,x_2,y_2))\in E(D)\iff (1)\text{ or } (2)
\end{equation*}
with the conditions
\begin{enumerate}
\item[$(1)$] $i$ is even and $y_1+y_2\neq 0$
\item[$(2)$] $i$ is odd and $x_1+x_2\neq 0$.
\end{enumerate}
We have to prove that the image of every edge is again an edge. We consider
\begin{eqnarray*}
\sigma(((i,x_1,y_1),(i+1,x_2,y_2)))&=&(\sigma((i,x_1,y_1)),\sigma((i+1,x_2,y_2)))\\
&=&((i+1,y_1,x_1),(i+2,y_2,x_2))
\end{eqnarray*}
and the two cases:
\begin{enumerate}
\item $i$ is even. Then we must have $y_1+y_2\neq 0$. Obviously, $i+1$ is odd and by condition $(2)$ from above $((i+1,y_1,x_1),(i+2,y_2,x_2))$ is an edge.
\item $i$ is odd. Analogously, we have $x_1+x_2\neq 0$ and $i+1$ is even and by condition $(1)$ from above $((i+1,y_1,x_1),(i+2,y_2,x_2))$ is an edge.
\end{enumerate}
Since $\sigma((i,$\minus$,$\minus$))=(i+1,$\minus$,$\minus$)$, the baseline is stabilized setwise. Also 
\begin{equation*}
\varphi(\sigma((i,x,y)))=\varphi((i+1,y,x))=i+1=\varphi((i,x,y))+1
\end{equation*}
holds.

\item[$\theta_\text{\oo}$:] The automorphism $\theta_\text{\oo}$ must map
\begin{equation*}
\theta_\text{\oo}:(0,\text{\oo},y)\mapsto(0,\text{\minus},y).
\end{equation*}
We choose the automorphism that exchanges the vertices
\begin{eqnarray*}
(0,\text{\oo},y)&\leftrightarrows&(0,\text{\minus},y)\quad\text{and}\\
(-1,\text{\oo},y)&\leftrightarrows&(-1,\text{\plus},y).
\end{eqnarray*}
This is illustrated in Figure~\ref{fig.t_o}. Note that this figure shows $D$  projected along the third coordinate. Thus, two more vertices hide behind every vertex and a $K_{3,3}$ ($AC_6$ respectively -- depending on the layer) hides behind every edge. It is enough to find an automorphism in that view. The third coordinate (that is being projected along) cannot cause any problem since inside a layer all the edges represent the same bipartite digraph (either $K_{3,3}$ or $AC_6$) and the third coordinate is stabilized. Within a layer, the edges that map onto each other are highlighted with the same style. Note that this automorphism (like all the following ones) involves 12 vertices and 126 edges. These are:
\begin{eqnarray*}
((-2,\text{\minus},y),(-1,\text{\oo},y))&\leftrightarrows&((-2,\text{\minus},y),(-1,\text{\plus},y))\\
((-2,\text{\oo},y),(-1,\text{\oo},y))&\leftrightarrows&((-2,\text{\oo},y),(-1,\text{\plus},y))\\
((-2,\text{\plus},y),(-1,\text{\oo},y))&\leftrightarrows&((-2,\text{\plus},y),(-1,\text{\plus},y))\\
((-1,\text{\minus},y),(0,\text{\minus},y))&\leftrightarrows&((-1,\text{\minus},y),(0,\text{\oo},y))\\
((-1,\text{\oo},y),(0,\text{\plus},y))&\leftrightarrows&((-1,\text{\plus},y),(0,\text{\plus},y))\\
((-1,\text{\oo},y),(0,\text{\minus},y))&\leftrightarrows&((-1,\text{\plus},y),(0,\text{\oo},y))\\
((0,\text{\minus},y),(1,\text{\minus},y))&\leftrightarrows&((0,\text{\oo},y),(1,\text{\minus},y))\\
((0,\text{\minus},y),(1,\text{\oo},y))&\leftrightarrows&((0,\text{\oo},y),(1,\text{\oo},y))\\
((0,\text{\minus},y),(1,\text{\plus},y))&\leftrightarrows&((0,\text{\oo},y),(1,\text{\plus},y))
\end{eqnarray*}

\begin{figure}[htdg]
\centering
\TubeTol\quad\TubeTor
\caption{$\theta_\text{\oo}$}\label{fig.t_o}
\end{figure}

\item[$\theta_\text{\plus}$:] The automorphism $\theta_\text{\plus}$ must map
\begin{equation*}
\theta_\text{\plus}:(0,\text{\plus},y)\mapsto(0,\text{\minus},y).
\end{equation*}
We choose the automorphism that exchanges the vertices
\begin{eqnarray*}
(0,\text{\plus},y)&\leftrightarrows&(0,\text{\minus},y)\quad\text{and}\\
(-1,\text{\plus},y)&\leftrightarrows&(-1,\text{\minus},y).
\end{eqnarray*}
 This is illustrated in Figure~\ref{fig.t_p}. For the rest of the proof we will not list the actions on the edges since they are obvious from the figures. 

\begin{figure}[htdg]
\centering
\TubeTpl\quad\TubeTpr
\caption{$\theta_\text{\plus}$}\label{fig.t_p}
\end{figure}

\item[$\theta^\text{\oo}$:] The automorphism $\theta^\text{o}$ must map
\begin{equation*}
\theta^\text{\oo}:(0,x,\text{\oo})\mapsto(0,x,\text{\minus}).
\end{equation*}
We choose the automorphism that exchanges the vertices
\begin{eqnarray*}
(0,x,\text{\oo})&\leftrightarrows&(0,x,\text{\minus})\quad\text{and}\\
(1,x,\text{\oo})&\leftrightarrows&(1,x,\text{\plus}).
\end{eqnarray*}
This is illustrated in Figure~\ref{fig.t^o}.

\begin{figure}[htdg]
\centering
\TubeTOl\quad\TubeTOr
\caption{$\theta^\text{\oo}$}\label{fig.t^o}
\end{figure}

\item[$\theta^\text{\plus}$:] The automorphism $\theta^\text{\plus}$ must map
\begin{equation*}
\theta^\text{\plus}:(0,x,\text{\plus})\mapsto(0,x,\text{\minus}).
\end{equation*}
We choose the automorphism that exchanges the vertices
\begin{eqnarray*}
(0,x,\text{\plus})&\leftrightarrows&(0,x,\text{\minus})\quad\text{and}\\
(1,x,\text{\plus})&\leftrightarrows&(1,x,\text{\minus}).
\end{eqnarray*}
This is illustrated in Figure~\ref{fig.t^p}.
\begin{figure}[htdg]
\centering
\TubeTPl\quad\TubeTPr
\caption{$\theta^\text{\plus}$}\label{fig.t^p}
\end{figure}

\item[$\alpha_\text{\oo}$:] The automorphism $\alpha_\text{\oo}$ must map the edge
\begin{equation*}
((0,\text{\minus},\text{\minus}),(1,\text{\oo},y))\mapsto((0,\text{\minus},\text{\minus}),(1,\text{\minus},y)).
\end{equation*}
We choose the automorphism that exchanges the vertices 
\begin{eqnarray*}
(1,\text{\oo},y)&\leftrightarrows&(1,\text{\minus},y)\quad\text{and}\\ (2,\text{\plus},y)&\leftrightarrows&(2,\text{\oo},y).
\end{eqnarray*}
 This is illustrated in Figure~\ref{fig.a_o}.

\begin{figure}[htdg]
\centering
\TubeAol\quad\TubeAor
\caption{$\alpha_\text{\oo}$}\label{fig.a_o}
\end{figure}

\item[$\alpha_\text{\plus}$:] The automorphism $\alpha_\text{\plus}$ must map the edge
\begin{equation*}
((0,\text{\minus},\text{\minus}),(1,\text{\plus},y))\mapsto((0,\text{\minus},\text{\minus}),(1,\text{\minus},y)).
\end{equation*}
We choose the automorphism that exchanges the vertices
\begin{eqnarray*}
(1,\text{\plus},y)&\leftrightarrows&(1,\text{\minus},y)\quad\text{and}\\ (2,\text{\plus},y)&\leftrightarrows&(2,\text{\minus},y).
\end{eqnarray*}
This is illustrated in Figure~\ref{fig.a_p}.

\begin{figure}[htdg]
\centering
\TubeApl\quad\TubeApr
\caption{$\alpha_\text{\plus}$}\label{fig.a_p}
\end{figure}

\item[$\alpha^\text{\oo}$:] Finally, the automorphism $\alpha^\text{\oo}$ must map the edge 
\begin{equation*}
((0,\text{\minus},\text{\minus}),(1,x,\text{\oo}))\mapsto((0,\text{\minus},\text{\minus}),(1,x,\text{\minus})).
\end{equation*}
We choose the automorphism that exchanges the vertices
\begin{eqnarray*}
(0,x,\text{\oo})&\leftrightarrows&(0,x,\text{\plus})\quad\text{and}\\
(1,x,\text{\minus})&\leftrightarrows&(1,x,\text{\oo}).
\end{eqnarray*}
This is illustrated in Figure~\ref{fig.a^o}.

\begin{figure}[htdg]
\centering
\TubeAOl\quad\TubeAOr
\caption{$\alpha^\text{\oo}$}\label{fig.a^o}
\end{figure}

\end{enumerate}
\end{proof}

\begin{thm}
There is a connected HAT digraph with Property $Z$ and finite fibers which has an associated digraph that is {\em not} complete bipartite, i.e., Conjecture~\ref{conj} does not hold.
\end{thm}

\begin{proof}
The digraph $D$ defined in $(\ref{eq.D})$ has fibersize $9$. By Lemma~\ref{lemma} its associated digraph is not complete bipartite. It follows from Lemma~\ref{lemma} as well that $D$ is connected. By Proposition~\ref{prop}, $D$ is HAT.
\end{proof}

\secA{Remarks}\label{sec.rem}

\secB{On the overlap with the work of \names{DeVos}, \names{Mohar} and \names{\v{S}\'amal}}\label{sec.art}

Our digraph $D$ was recently constructed by \names{DeVos}, \names{Mohar} and \names{\v{S}\'amal} \cite[Construction 2]{9} without the use of the layerwise direct product. They basically gave the joined condition from the above proof directly and mixed the coordinates in a way such that the shift does not need to flip coordinates. Our constructions from Section~\ref{sec.other} correstponds in a similar way to \cite[Construction 3]{9}.

\secB{Möller's digraph} \label{sec.moeller}

The digraph \names{Möller} constructed in \cite{2} looked promising because the arranging-automorphisms work. Unfortunately, he missed to check transitivity. His approach was to use an alternating eight cycle of $K_{2,2}$'s as $\Delta$ and concatenate them cleverly. The easiest way to see that his digraph is not HAT is probably to see that an edge cannot be mapped to an edge that agrees on the terminal vertex  and lies in the same $K_{2,2}$ (the contradiction arises soon in the outward direction).

\secB{On the automorphisms} \label{sec.auto}

In the proof of Proposition~\ref{prop} we chose the $n$-arc
\begin{equation*}
z=(((i,\text{\minus},\text{\minus}),(i+1,\text{\minus},\text{\minus})))_{i=-n,\dots,-1}
\end{equation*}
and the automorphism
\begin{equation*}
\psi\coloneqq(\sigma^{-1}\circ\alpha_{i_n}\circ\alpha^{i^n})\circ\dots\circ(\sigma^{-1}\circ\alpha_{i_1}\circ\alpha^{i^1})\circ\theta_{i_0}\circ\theta^{i^0}\circ\sigma^m.
\end{equation*}
We could as well have chosen the $n$-arc
\begin{equation*}
z'=(((i,\text{\minus},\text{\minus}),(i+1,\text{\minus},\text{\minus})))_{i=0,\dots,n-1}
\end{equation*}
in which case we would not need the shifts $\sigma^{-1}$. Moreover, we would not need that the shift stabilizes the baseline. But we would have needed an argument that the arranging-automorphisms work in every layer (this is not difficult, since all the layers look the same, i.e., just differ by flips of the coordinates). The automorphism would then appear as
\begin{equation*}
\psi'\coloneqq({}_n\alpha_{i_n}\circ{}_n\alpha^{i^n})\circ\dots\circ({}_1\alpha_{i_1}\circ{}_1\alpha^{i^1})\circ\theta_{i_0}\circ\theta^{i^0}\circ\sigma^m,
\end{equation*}
where the lower left indices denote the working-layer. More importantly, we notice that the different $\alpha$- and $\theta$-automorphisms act only on one of the coordinates and thus we need not perform them one after the other but can apply them simultaneously in every layer. We can therefore denote
\begin{equation*}
\Theta\coloneqq\left(\begin{array}{c}\id\\\theta_{i_0}\\\theta^{i^0}\end{array}\right),\quad
A_j\coloneqq\left(\begin{array}{c}\id\\{}_j\alpha_{i_j}\\{}_j\alpha^{i^j}\end{array}\right),
\end{equation*}
and our automorphism reads
\begin{equation*}
\psi'=A_n\circ\dots\circ A_1\circ\Theta\circ\sigma^m.
\end{equation*}

\secA{More constructions}\label{sec.more}

\secB{The example of McKay and Praeger}\label{sec.mckay}

\begin{figure}[htdg]
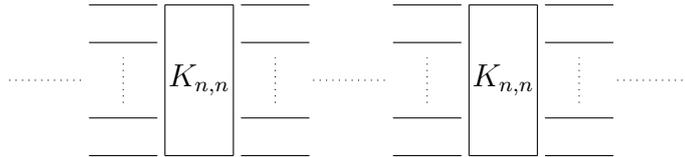

\centering
\McKAY[1]
\caption{Factor of the digraph of Praeger and McKay}\label{fig.McKay}
\end{figure}

In \cite{1} a nontrivial HAT digraph is constructed that has Property $Z$, finite fibers and complete bipartite digraphs as associated digraph. This example can be realized using the layerwise direct product. One uses $K_{n,n}$'s instead of $K_{3,3}$'s and reaplaces $AC_6$ with a matching ($n$ horizontal edges). One places $m$ matchings in between the $K_{n,n}$s and multiplies the digraph with its $m$ shifts (see Figure~\ref{fig.McKay}).

\secB{Other factors}\label{sec.other}

\begin{figure}[htdg]
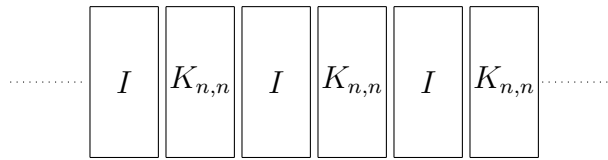

\centering
\LITTLE[1]
\caption{Little variation of the factors from $D$}\label{fig.little}
\end{figure}

Considering the automorphisms in the proof, every $1$-arc transitive, connected, bipartite, noncomplete bipartite digraph (with equal partition sizes) concatenated with a suitable $K_{n,n}$ would have done the job. (Connectedness is not necessary in every situation.) From now on, we will refer to the non-complete bipartite layers as \defstil{involvers} (unless they have outvalency smaller $2$). Figure~\ref{fig.little} shows this situation where $I$ stands for the Involver.

Like in Section~\ref{sec.mckay}, we are not restricted to two factors. We are free to place arbitrary many matchings (and some more $K_{n,n}$'s) between the $K_{n,n}$'s and $1$-arc transitive digraphs as long as the factor stays periodic (say with the length $l$ of the period). Then the product of the $l$ different shifts is again HAT.

\begin{figure}[htdg]
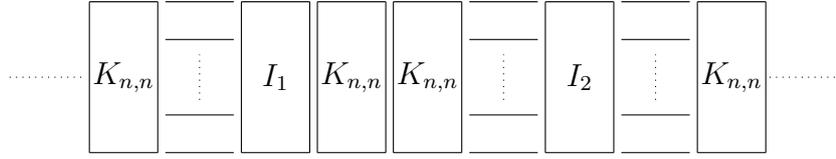

\centering
\MOREGENERAL[1]
\caption{More general factors}\label{fig.moregen}
\end{figure}

Moreover, we can use different $1$-arc transitive, bipartite digraphs as involvers. We just have to ensure that there is a $K_{n,n}$ between each two of them. Again multiplying the $l$ shifts yield a HAT digraph (see Figure~\ref{fig.moregen}).

Indeed we do not even need to keep the fibersizes of the factor constant but can alter them periodically.

\begin{figure}[htdg]
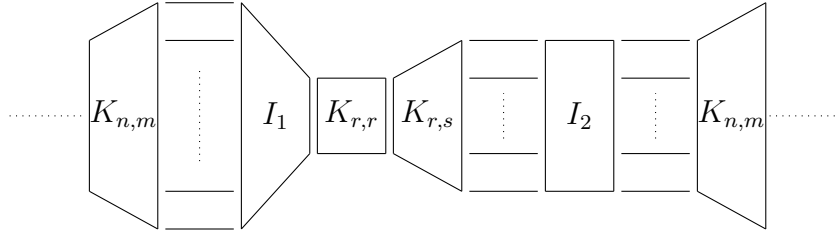

\centering
\GENERAL[1]
\caption{Even more general factors}\label{fig.gen}
\end{figure}

In this general situation we need to think about the shift globally and the transitivity- and arranging-automorphisms locally in one factor from the $K_{n,m}$ preceding an involver to the $K_{r,s}$ succeeding it.

\begin{figure}[htdg]
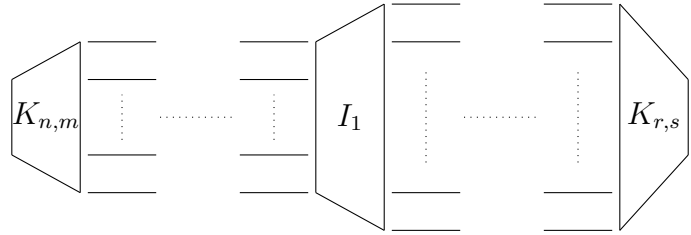

\centering
\LOCAL[1]
\caption{Local view of a factor}\label{fig.local}
\end{figure}

The shift will simply increase the $\Z$ coordinate and circularly shift the other coordinates:
\begin{equation*}
\sigma:(i,x_1,\dots,x_n)\mapsto(i+1,x_n,x_1,\dots,x_{n-1})
\end{equation*}
We must consider the $\alpha$- and $\theta$-automorphisms in all the layers of the local view (we can skip the $\theta$ for the very first and the $\theta$ and $\alpha$ for the very last layer because they actually belong to the preceding or succeeding local view).

\begin{itemize}
\item The $\alpha$-automorphisms for the $K_{n,m}$ map an arbitrary edge $e_1$ onto an other $e_0$. The corresponding map for the terminal vertices $t_{e_1}\mapsto t_{e_0}$ will be transported by the matchings to the involver. There it can be realized because the $1$-arc transitivity of the involver guarantees transitivity on the initial vertices. But it will produce some permutation on its initial vertices (which is transported back by the matchings and stopped at the $K_{n,m}$) and some permutation on its terminal vertices (which is transported by the matchings to the $K_{r,s}$ where it is stopped).
\item The $\theta$-automorphisms for the terminal vertices of the $K_{n,m}$ works in exactly the same way.
\item In the layers between the $K_{n,m}$ and the involver there are no branchings, thus we have $\alpha=\id$.
\item The $\theta$-automorphisms in these layers must again be considered first in the forward direction where they have some effect on the initial and terminal vertices of the involver. These effects are again transported by the matchings to the complete bipartite layers, where they are stopped.
\item The same happens for the $\theta$-automorphisms in the initial layer of the involver -- only the forward transportation to the involver is now trivial.
\item For the $\alpha$-automorphisms in the involver we need its $1$-arc transitivity (which is more than the transitivity on the initial vertices). But then we get the transportation and stoppage as above.
\item The rest of the layers can be dealt analogously by first considering the backward direction.
\end{itemize}

Considering this situation, we find, that we can place more than one involver between two complete bipartite layers, if they are compatible in the way, that there is a set of permutations of the vertices of the terminal layer of the one involver, that realizes all the needed automorphisms and can be realized by the next involver on its initial layer (where the permutations get transported by the matchings in between).

\begin{figure}[htdg]
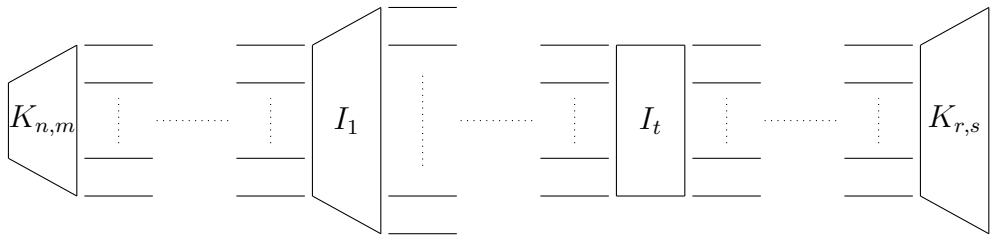

\centering
\MULTI[1]
\caption{Compatible involvers}\label{fig.multi}
\end{figure}

\secB{Factors with infinite fibers} \label{sec.inffib}

There is no need to restrict to finite bipartite digraphs. We could build up our factors from arbitrary $1$-arc transitive, bipartite involvers with finite and infinite partition sizes concatenated with complete bipartite digraphs. Doing this, we loose the local finiteness. Moreover, our HAT digraphs will have only one end.

\secB{Infinitely many factors}\label{sec.inffac}

Consider the factor from Figure~\ref{fig.gen} to be non-periodic. We still could use it to construct a HAT digraph if we multiply all its shifts ($\Z$-many). The vertices of such a product can be understood as a $\Z$ coordinate together with a two way infinite sequence with entries at the $i$-th position from the $i$-th layer of the factor (which is different for every shift).

This construction goes along with questions of connectedness and local finiteness. Let us first consider the latter. If there are finitely many involvers we are free to set all layers to the right of the rightmost involver and to the left of the leftmost involver to matchings. If we moreover use only locally finite layers, the product will be locally finite.

\begin{figure}[htdg]
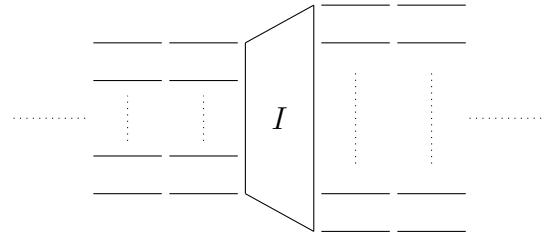

\centering
\INF[1]
\caption{A factor realizing a construction from \cite{1}}\label{fig.seq}
\end{figure}

Consider the factor that has an involver (or $K_{n,m}$) in the $(0,1)$ layer and matchings in all the other layers (see Figure~\ref{fig.seq}). Given every layer contains vertices $(i,0)$ and ($i,1)$ then there will be no finite walk connecting $(0,(\dots,0,0,0,\dots))$ and $(0,(\dots,1,1,1,\dots))$ in the product $G$. Thus,  $G$ is not connected -- indeed its components are formed by the vertices whose sequences differ in at most finitely many entries (this is modulo the shifts by which the sequences differ from layer to layer). Since the entire digraph $G$ is HAT, all its components $C$ are isomorphic and HAT. Now $C$ is isomorphic to the digraphs constructed in \cite[Theorem 4.8]{1}. These digraphs are direct products of a so called sequences digraph with the integer line.

The question for connectedness of such an infinite product is somewhat more involved and we do not address it here.

\secB{Nonisomorphic factors}\label{sec.nonisofac}

Up to now we considered non-HAT factors and created a HAT product. Finally, we want to build new HAT digraphs from existing HAT digraphs. From \cite[Lemma 4.3, (a)]{1} it follows that the layerwise direct product of two HAT digraphs with Property $Z$ is again HAT with Property $Z$. If the factors had finite fibers so has their product. If one of the factors has a noncomplete bipartite associated digraph so has the product. Thus, we can obtain numerous HAT digraphs with Property $Z$ by multiplying some of the above digraphs, thereby gaining even more counterexamples of Conjecture~\ref{conj}.

Building the layerwise direct product of the $1$-invalent-$2$-outvalent tree with the $2$-invalent-$1$-outvalent tree, one gets a Cayley graph of the Lamplighter group. It has $K_{2,2}$ as associated digraph. If one replaces the second tree by the $3$-invalent-$1$-outvalent tree, one gets a digraph that is not quasi-isomorphic to any Cayley graph, it has $K_{3,2}$ as associated digraph. Such digraphs are called Diestel--Leader graphs, and they got to fame for these mentioned properties. Any other choice of regular trees will also yield a HAT product.

We close with a remark on the digraphs we mentioned above to realize the construction from \cite[Theorem 4.8]{1}. Their layerwise direct product with some HAT digraph with Property $Z$ will realize the digraphs from \cite[Corollary 4.9]{1}.

\end{document}